\documentclass[11pt]{amsart}
\usepackage[T1]{fontenc}
\usepackage{lmodern}
\usepackage[utf8]{inputenc}
\usepackage[english]{babel}
\usepackage{amsmath,amssymb,amsthm,mathtools,longtable,array}
\usepackage[margin=1.08in]{geometry}
\usepackage{microtype}
\usepackage[hidelinks]{hyperref}

\newcommand{\kk}{\Bbbk}

\newcommand{\ord}{\operatorname{ord}}
\newcommand{\PIdeg}{\operatorname{PIdeg}}
\newcommand{\End}{\operatorname{End}}

\newcommand{\Ad}{\operatorname{Ad}}
\newcommand{\mi}{\mu_\infty}
\theoremstyle{plain}
\newtheorem{theorem}{Theorem}[section]
\newtheorem{proposition}[theorem]{Proposition}
\newtheorem{lemma}[theorem]{Lemma}
\newtheorem{corollary}[theorem]{Corollary}
\theoremstyle{definition}
\newtheorem{definition}[theorem]{Definition}
\newtheorem{question}[theorem]{Question}
\theoremstyle{remark}
\newtheorem{remark}[theorem]{Remark}

\numberwithin{equation}{section}

\begin{document}
\begin{center}
{\LARGE\bfseries  Polynomial identities and primitivity in double extension regular algebras of type (14641)\par}
\vspace{1.2em}
{\large James G\'omez$^{1}$ $\cdot$ Helbert Venegas$^{2}$\par}
\vspace{.6em}
{\small $^{1}$Grupo de \'Algebra y An\'alisis, Universidad Pedag\'ogica y Tecnol\'ogica de Colombia, Tunja, Colombia\\
\texttt{james.gomez@uptc.edu.co}\par}
\vspace{.4em}
{\small $^{2}$Departamento de Matem\'aticas, Universidad Militar Nueva Granada, Bogot\'a, Colombia\\
\texttt{helbert.venegas@unimilitar.edu.co}\par}
\end{center}

\vspace{1em}
\noindent\textbf{Abstract}\quad
We determine which members of the twenty-six families of double extension regular algebras of type $(14641)$ in the Zhang--Zhang classification satisfy a polynomial identity. Retaining the scalar twist parameter, we compute the PI degree for twenty-one of the twenty-two families with PI members and give three possibilities for the remaining family. We establish exact primitivity criteria for the nondegenerate members of eight families and show that sixteen families have no primitive members. The proofs use monomial data, normal subalgebras and crossed products of quantum tori. Twisting by a root of unity preserves the PI property, whereas finite order of the double-extension determinant is necessary but not sufficient for an algebra in the list to be PI.
\medskip

\noindent\textbf{Keywords}\quad Double Ore extension $\cdot$ Artin--Schelter regular algebra $\cdot$ Polynomial identity $\cdot$ PI degree $\cdot$ Primitive ring $\cdot$ Zhang twist

\medskip

\noindent\textbf{Mathematics Subject Classification}\quad Primary 16R20 $\cdot$ 16S36 $\cdot$ Secondary 16S35 $\cdot$ 16D60 $\cdot$ 16E65

\section{Introduction}

Zhang and Zhang introduced double Ore extensions in \cite{ZZ08} and used them in \cite{ZZ09} to build Artin--Schelter regular algebras of global dimension four. A trimmed double extension of type $(14641)$ that is not an Ore extension of a three-dimensional regular algebra lies in one of twenty-six families, $\mathbb A$ to $\mathbb Z$ \cite[Section 4]{ZZ09}; some members are still iterated Ore extensions \cite{RR,CLM}. The introduction of \cite{ZZ09} ends with questions about each algebra $B$ in the list, the first being
\begin{quote}
(a) Is $B$ primitive? Does $B$ satisfy a polynomial identity? What is $\operatorname{Spec}B$?
\end{quote}
For the algebra $B(h)$ of \cite{ZZ08}, which is family $\mathbb R$, they had already shown that it is PI exactly when $h$ is a root of unity \cite[Corollary 4.8]{ZZ08}.

Our earlier studies of the families $\mathbb G$ and $\mathbb K$ \cite{GV-G,GV-K} determined their PI criteria and PI degrees, alongside their centers and ozone groups.\footnote{The computations of \cite{GV-G} were also checked by computer \cite{GV-codeG}.} Here we settle the PI question for all twenty-six families, give the PI degree for all but one family with PI members, and decide primitivity in most cases. Rubiano \cite{Rub-ozone} also obtained PI criteria for $\mathbb K$ and $\mathbb O$, and PI degrees in the parameter regimes specified below. Earlier work concerns centers \cite{Rub-center}, Gr\"obner--Shirshov bases \cite{HHR}, differential smoothness \cite{RR}, homological invariants \cite{Rub-homol} and normal ordering \cite{Rub-normal}. For single Ore extensions see \cite{PV,DS,Jondrup,LM,LM2,Haynal}; for other algebras of polynomial type see \cite{GG-poly,GG-PBW}, and for related work on Gelfand--Kirillov dimension and centers \cite{LV-GK,LV-center,LV-fractions}.

In \cite{HHR,Rub-center} the twist parameter is set to $h=1$. That is harmless for PBW bases, not here: by Table~\ref{tab:main}, no member of $\mathbb C$, $\mathbb E$, $\mathbb F$, $\mathbb I$, $\mathbb J$, $\mathbb S$, $\mathbb T$, $\mathbb U$ with $h=1$ is PI, though all eight families have PI members.

Write $\mi$ for the group of roots of unity in $\kk$, and $s$ for a square root of $f$ whenever a family has a parameter $f$; the conditions below do not depend on the choice of $s$, and $i$ is a fixed square root of $-1$.

\begin{theorem}\label{thm:A}
Let $\kk$ be algebraically closed of characteristic zero, and let $B$ be an algebra in one of the twenty-six families of \cite[Section 4]{ZZ09}, with admissible parameters. Then $B$ is PI if and only if the condition in Table~\ref{tab:main} holds.
\end{theorem}

\begin{longtable}{|c|c|c|}
\caption{The PI criterion. The second column lists the parameters as in \cite[Section 4]{ZZ09}; $h\in\kk^\times$ is the twist parameter.}\label{tab:main}\\
\hline
Family & Parameters & $B$ is PI if and only if\\\hline
\endfirsthead
\hline
Family & Parameters & $B$ is PI if and only if\\\hline
\endhead
$\mathbb A$ & $P=(1,1)$ & never\\
$\mathbb B$ & $p^2=-1$ & $h\in\mi$\\
$\mathbb C$ & $p^2+p+1=0$ & $3h^2\in\mi$\\
$\mathbb D$ & $p\ne0$ & $p,h\in\mi$ and $p\ne\pm1$\\
$\mathbb E$, $\mathbb J$ & $p^2=-1$, resp.\ $q^2=-1$ & $2h^2\in\mi$\\
$\mathbb F$, $\mathbb I$ & $p^2=-1$, resp.\ $q^2=-1$ & $2h^2\in\mi$\\
$\mathbb G$ & $p\ne0,\pm1$, $f\ne0$ & $p,h\in\mi$\\
$\mathbb H$ & $Q=(1,1)$ & never\\
$\mathbb K$, $\mathbb L$ & $q\ne0$, $f\ne0$ & $q,f,h\in\mi$\\
$\mathbb M$ & $f\ne1$ & $f\ne0$ and $h(s\pm1)\in\mi$\\
$\mathbb N$ & $f^2\ne g^2$ & $f\pm g\in\mi$\\
$\mathbb O$, $\mathbb P$ & $f\ne1$ & $f\ne0$ and $h(1\pm s)\in\mi$\\
$\mathbb Q$, $\mathbb V$ & --- & $h\in\mi$\\
$\mathbb R$ & --- & $h\in\mi$\\
$\mathbb S$, $\mathbb T$, $\mathbb U$ & --- & $2h\in\mi$\\
$\mathbb W$ & $f\ne-1$ & $f\ne0$ and $h(s\pm i)\in\mi$\\
$\mathbb X$ & --- & never\\
$\mathbb Y$ & $f\in\kk$ & never\\
$\mathbb Z$ & $f(1+f)\ne0$ & $h(s\pm i)\in\mi$\\
\hline
\end{longtable}

The families $\mathbb A$, $\mathbb H$, $\mathbb X$, $\mathbb Y$ are treated in Section~\ref{sec:noPI}, the fourteen monomial families in Section~\ref{sec:mono}, $\mathbb D$, $\mathbb G$, $\mathbb Q$, $\mathbb V$ in Section~\ref{sec:normal}, and $\mathbb C$, $\mathbb F$, $\mathbb I$ in Section~\ref{sec:CFI}; for $\mathbb R$ see \cite[Corollary 4.8]{ZZ08}. The families $\mathbb M$, $\mathbb O$, $\mathbb P$, $\mathbb W$ allow $f=0$, where the monomial description breaks down. Those members are never PI, like $\mathbb D(\pm1,h)$; in fact $\mathbb O(0,h)=\mathbb D(-1,h)$. The twist only enters through its class in $\kk^\times/\mi$:

\begin{theorem}\label{thm:B}
Let $B(h)$ be the double extension with mixed matrix $h\Sigma_0$, for a fixed $\Sigma_0$ of type $(14641)$. If $\zeta\in\mi$, then $B(h)$ is PI if and only if $B(h\zeta)$ is PI.
\end{theorem}

This is Theorem~\ref{thm:twist}; a general Zhang twist does not preserve the PI property, since $\kk_q[x,y]$ is a twist of $\kk[x,y]$.

For the PI degree we write $Q(B)$ as a crossed product of a quantum torus by a group of order at most eight. This gives closed formulas for twenty-one of the twenty-two families with PI members, collected in Table~\ref{tab:degrees}; for instance
\begin{align*}
\PIdeg\mathbb S(h)&=8\ord(16h^4), & \PIdeg\mathbb E(h)&=8\ord(4h^4), & \PIdeg\mathbb B(p,h)&=8\ord(h^4),\\
\PIdeg\mathbb R(h)&=4\ord(h^2), & \PIdeg\mathbb F(p,h)&=8\ord(-4h^4) .
\end{align*}
For $\mathbb C$ we only show that the PI degree is $n$, $3n$ or $9n$, with $n$ explicit (Proposition~\ref{prop:bounds}).

No PI member is primitive (Proposition~\ref{prop:notprimitive}). For the non-PI members, in sixteen families the quotient division ring has a central element that does not depend on the twist, such as $(x_1^4+x_2^4)^2/(x_1^4x_2^4)$ in $\mathbb B$, so no member is primitive. In $\mathbb M,\mathbb N,\mathbb O,\mathbb P,\mathbb W,\mathbb Z$
(with $f\ne0$ in $\mathbb M,\mathbb O,\mathbb P,\mathbb W$), a member
is primitive exactly when neither $u$ nor $v$ is a root of unity
(Theorem~\ref{thm:primmono}(ii)). In $\mathbb D$ with $p^2\ne1$
and in $\mathbb G$, a non-PI member is primitive exactly when neither
$p$ nor $hp$ is a root of unity (Theorem~\ref{thm:DGprim}). What remains open is $\mathbb A$, $\mathbb Y$ and five degenerate specializations, for $h\notin\mi$ (Question~\ref{q:remainingprim}).

Zhang and Zhang showed that $B(h)$ is PI exactly when $\det\sigma$ has finite order \cite[Corollary 4.8]{ZZ08}; on the whole list this fails (Remark~\ref{rem:detsigma}). Except for two scalars taken from \cite[Proposition~5.17]{ZZ09}, every computation is done here by hand.

\section{Preliminaries}\label{sec:prelim}

Throughout, $\kk$ is algebraically closed of characteristic zero.

\subsection{The families}

We follow \cite[Sections 3 and 4]{ZZ09}. A trimmed double extension $B=(\kk_Q[x_1,x_2])_P[y_1,y_2;\sigma]$ is generated by $x_1,x_2,y_1,y_2$ with
\begin{align}
x_2x_1&=q_{12}x_1x_2+q_{11}x_1^2, & y_2y_1&=p_{12}y_1y_2+p_{11}y_1^2,\notag
\end{align}
where $Q=(q_{12},q_{11})$ and $P=(p_{12},p_{11})$, and four mixed relations
\begin{equation}\notag
\begin{pmatrix}y_1x_1\\ y_1x_2\\ y_2x_1\\ y_2x_2\end{pmatrix}=\Sigma\begin{pmatrix}x_1y_1\\ x_2y_1\\ x_1y_2\\ x_2y_2\end{pmatrix}.
\end{equation}
In every family except $\mathbb N$, $\Sigma=h\Sigma_0$ with $\Sigma_0$ the matrix of \cite[Section 4]{ZZ09} and $h\in\kk^\times$; in $\mathbb N$ the matrix depends on $f,g$ and a twist only rescales them. We always read the relations off these matrices. For $\mathbb Z$ this matters: the expanded list of relations at the end of \cite[Subcase 4.4.4]{ZZ09} has $x_1y_2$ in place of $x_1y_1$ in the relations for $y_1x_1$ and $y_2x_2$, whereas the matrix $\Sigma_{\mathbb Z}$, and the computation in the proof of \cite[Proposition 5.16]{ZZ09}, give
\[
y_1x_1=h(x_1y_1+x_2y_2),\qquad y_2x_2=h(fx_1y_1-x_2y_2).
\] We write $A_x=\kk_Q[x_1,x_2]$ and $A_y=\kk_P[y_1,y_2]$ for the two planes, and $\kk_q[x_1,x_2]$ for the quantum plane with $x_2x_1=qx_1x_2$.

By \cite[Theorem 0.1]{ZZ09} every $B$ in the list is a strongly noetherian, Auslander regular and Cohen--Macaulay domain, AS-regular of global dimension four, and the ordered monomials $x_1^ax_2^by_1^cy_2^d$ form a basis. It is $\mathbb Z^2$-graded with $\deg x_i=(1,0)$ and $\deg y_i=(0,1)$. Multiplication gives an isomorphism of vector spaces $A_x\otimes A_y\to B$; in particular, if $a_1,\dots,a_r\in A_x$ and $b_1,\dots,b_t\in A_y$ are linearly independent, so are the products $a_jb_l$. Since $B$ is a double extension in the sense of \cite{ZZ08}, it is also a free right $A_x$-module on a basis of $A_y$, so multiplication $A_y\otimes A_x\to B$ is bijective as well. We use this mostly in the following form. The planes are noetherian domains, so $Q(A_x)$ and $Q(A_y)$ are subrings of $Q(B)$.

\begin{lemma}\label{lem:fracind}
If $\alpha_1,\dots,\alpha_r\in Q(A_x)$ are linearly independent over $\kk$, then they are right linearly independent over $Q(A_y)$. In particular, if $\beta_1,\dots,\beta_t\in Q(A_y)$ are also linearly independent, so are the products $\alpha_j\beta_l$ in $Q(B)$. The same holds with $A_x$ and $A_y$, and left and right, exchanged.
\end{lemma}

\begin{proof}
Let $\sum_j\alpha_j\beta_j=0$ with $\beta_j\in Q(A_y)$. Take a common left denominator in $A_x$ and a common right denominator in $A_y$: $\alpha_j=g^{-1}a_j$ and $\beta_j=b_jg'^{-1}$ with $0\ne g,a_j\in A_x$ and $0\ne g'$, $b_j\in A_y$. The $a_j$ are linearly independent, since left multiplication by $g$ is injective, and $\sum_ja_jb_j=g\bigl(\sum_j\alpha_j\beta_j\bigr)g'=0$ in $B$. As $B$ is a free right $A_y$-module on any basis of $A_x$, every $b_j$ is zero. For the second claim apply the first to $\sum_j\alpha_j\bigl(\sum_lc_{jl}\beta_l\bigr)$.
\end{proof}

\begin{lemma}\label{lem:qplane}
Let $a\in A_x$ and $b\in A_y$ be nonzero and homogeneous of positive degree, with $ba=\mu\,ab$ for some $\mu\in\kk^\times$. If $B$ is PI, then $\mu\in\mi$.
\end{lemma}

\begin{proof}
The powers of $a$ have different degrees, so they are linearly independent, and the same holds for $b$. So the monomials $a^jb^l$ are independent and $\kk\langle a,b\rangle$ is the quantum plane $\kk_\mu[a,b]$, which is PI only when $\mu\in\mi$.
\end{proof}

\begin{lemma}\label{lem:finite}
Let $R$ be a domain and $S\subseteq R$ a subring which is a right noetherian PI domain, with $R$ finitely generated as a right $S$-module. Then $R$ is PI. If $d=\dim_{Q(S)}R\otimes_SQ(S)$, then $\PIdeg R\le d\,\PIdeg S$. The same holds with right replaced by left.
\end{lemma}

\begin{proof}
$S$ is a right Ore domain, and by Posner's theorem $Q(S)$ is a division algebra of dimension $m^2$ over its center $Z$, with $m=\PIdeg S$ \cite[Theorem 13.6.5]{MR}. Since $R$ is a domain it is torsion-free over $S$, so $R\to V=R\otimes_SQ(S)$ is injective, and $V$ is a right $Q(S)$-space of dimension $d<\infty$. Left multiplication gives a ring homomorphism $R\to\End_{Q(S)}(V)\cong M_d(Q(S))$, injective because $b\otimes1\ne0$ for $b\ne0$. If $L$ is a splitting field of $Q(S)$ over $Z$, then $M_d(Q(S))\subseteq M_{dm}(L)$, so $R$ satisfies the identities of $M_{dm}(L)$.
\end{proof}

\begin{lemma}\label{lem:shift}
Let $D$ be a division algebra finite-dimensional over its center $F\supseteq\kk$, and let $t,w\in D$, $t\notin\kk$, $w\ne0$, with $wtw^{-1}=\mu t+c$, $\mu\in\kk^\times$, $c\in\kk$. Then the affine map $X\mapsto\mu X+c$ has finite order. In particular, if $c=0$ then $\mu\in\mi$.
\end{lemma}

\begin{proof}
$t$ is algebraic over $F$; let $m\in F[X]$ be its minimal polynomial. Conjugation by $w^n$ fixes $F$, so $m$ vanishes at $w^ntw^{-n}=\mu^nt+c_n$ for every $n$, and these elements lie in the field $F(t)$. If the affine map had infinite order, each nontrivial power would have at most one fixed point, which lies in $\kk$, while $t\notin\kk$. So the elements $\mu^nt+c_n$ would be pairwise distinct, and $m$ would have infinitely many roots in a field.
\end{proof}

\begin{lemma}\label{lem:shiftmat}
Let $D$ be as in Lemma~\ref{lem:shift}, $t\in D\setminus\kk$, and $y_1,y_2\in D$ linearly independent over $\kk$ with
\[
y_jt=a\,ty_j+\sum_{l=1}^2C_{jl}y_l\qquad(j=1,2),\qquad a\in\{1,-1\},\quad C\in M_2(\kk).
\]
Then $C$ has no nonzero eigenvalue if $a=1$, and $C$ does not have two distinct eigenvalues if $a=-1$.
\end{lemma}

\begin{proof}
If $vC=\lambda v$, then $z=v_1y_1+v_2y_2\ne0$ satisfies $zt=(at+\lambda)z$. For $a=1$ and $\lambda\ne0$ this is a translation of infinite order, against Lemma~\ref{lem:shift}. For $a=-1$, take $z_1,z_2$ with eigenvalues $\lambda_1\ne\lambda_2$ and put $t_0=t-\lambda_2/2$, $\delta=(\lambda_1-\lambda_2)/2$. Then $z_2t_0=-t_0z_2$, $z_1t_0=(-t_0+2\delta)z_1$, and $w=z_1z_2$ satisfies $wt_0=(t_0-2\delta)w$, again a translation.
\end{proof}

\subsection{Twisting}

Fix $\Sigma_0$ and write $B(h)$ for the algebra with mixed matrix $h\Sigma_0$. For bihomogeneous $a,b\in B(1)$ of bidegrees $(k_a,j_a)$ and $(k_b,j_b)$ put
\begin{equation}\notag
a*_hb=h^{j_ak_b}\,ab .
\end{equation}
The exponent is bilinear in the bidegrees, so $*_h$ is associative, and the plane relations do not change while $y*_hx=h\,yx$. Hence $(B(1),*_h)\cong B(h)$; compare \cite[Lemma 3.3]{ZZ09}. Two consequences are used often. Let $N,N'\in B(1)$ have bidegrees $(0,d)$ and $(d,0)$. Then
\begin{equation}\label{eq:twistscalar}
\begin{aligned}
Nx_i&=c\,x_iN\ \text{ in }B(1) &&\Longrightarrow && Nx_i&=c\,h^d\,x_iN\ \text{ in }B(h),\\
y_jN'&=c'N'y_j\ \text{ in }B(1) &&\Longrightarrow && y_jN'&=c'h^{d}N'y_j\ \text{ in }B(h).
\end{aligned}
\end{equation} And composing twists, $(B(h),*_\zeta)\cong B(h\zeta)$.

\begin{theorem}\label{thm:twist}
Let $\zeta\in\kk^\times$ be a root of unity of order $m$. Then $B(h)$ is PI if and only if $B(h\zeta)$ is PI.
\end{theorem}

\begin{proof}
Let $S=\bigoplus_{m\mid k}B(h)_{(k,j)}$, the span of the bihomogeneous components whose $x$-degree is a multiple of $m$. If $b\in S$, then $a*_\zeta b=\zeta^{j_ak_b}ab=ab$. So $S$ is a subalgebra of $B(h)$ and of $(B(h),*_\zeta)\cong B(h\zeta)$, with the same multiplication.

Every basis monomial $x_1^ax_2^by_1^cy_2^d$ can be written as $x^{\alpha'}\cdot(x^{\alpha''}y_1^cy_2^d)$, where $x^{\alpha'}$ consists of the first $l$ letters of $x_1^ax_2^b$, $l$ is the residue of $a+b$ modulo $m$, and the second factor lies in $S$. So $B(h)=\sum_{\alpha'}x^{\alpha'}S$, with $x^{\alpha'}$ running over the finitely many monomials of degree less than $m$. Since $x^{\alpha'}$ has $y$-degree zero, $x^{\alpha'}*_\zeta c=x^{\alpha'}c$, and the same finite set generates $B(h\zeta)$ as a right $S$-module.

$S$ is right noetherian. Let $B_l$ be the span of the bihomogeneous components of $B(h)$ whose $x$-degree is congruent to $l$ modulo $m$, so that $B(h)=\bigoplus_{l}B_l$, $B_0=S$ and $SB_l\subseteq B_l$. For a right ideal $I$ of $S$, not necessarily graded, $IB(h)=\bigoplus_lIB_l$, and its part in $S$ is $IS=I$. So $I\mapsto IB(h)$ is injective on right ideals and preserves inclusions. It is a domain, being inside $B(h)$. If $B(h)$ is PI, so is $S$, and Lemma~\ref{lem:finite} shows that the domain $B(h\zeta)$ is PI. The converse is the same argument with $\zeta^{-1}$.
\end{proof}

So in each family the values of $h$ giving PI members form a union of cosets of $\mi$.

\section{Members that are never PI}\label{sec:noPI}

\begin{proposition}\label{prop:jordan}
The families $\mathbb A$ and $\mathbb H$ are never PI.
\end{proposition}

\begin{proof}
For $\mathbb A$, $P=(1,1)$, and for $\mathbb H$, $Q=(1,1)$. So $B$ contains a Jordan plane, which is not PI in characteristic zero.
\end{proof}

\begin{proposition}\label{prop:noPI}
The following algebras are not PI, for any value of the remaining parameters:
\[
\mathbb X(h),\quad\mathbb Y(f,h),\quad\mathbb D(\pm1,h),\quad\mathbb M(0,h),\quad\mathbb O(0,h),\quad\mathbb P(0,h),\quad\mathbb W(0,h).
\]
Moreover $\mathbb O(0,h)=\mathbb D(-1,h)$.
\end{proposition}

\begin{proof}
For $f=0$ the matrices of $\mathbb O$ and of $\mathbb D$ at $p=-1$ are both
\[
h\begin{pmatrix}1&0&0&0\\0&-1&1&0\\0&0&-1&0\\1&0&0&1\end{pmatrix},
\]
with $P=Q=(-1,0)$ in both cases. This proves the last claim.

In each of the six remaining cases one generator $n\in\{x_1,x_2\}$ is normal. Reading off the mixed relations,
\[
\begin{array}{c|c|c}
 & n & \text{relations of }n\text{ with }y_1,y_2\\\hline
\mathbb X & x_1 & y_1x_1=h\,x_1y_2,\ \ y_2x_1=h\,x_1y_1\\
\mathbb Y & x_1 & y_1x_1=h\,x_1y_1,\ \ y_2x_1=h\,x_1y_2\\
\mathbb D(p,h) & x_1 & y_1x_1=-hp\,x_1y_1,\ \ y_2x_1=hp\,x_1y_2\\
\mathbb M(0,h) & x_2 & y_1x_2=-h\,x_2y_2,\ \ y_2x_2=-h\,x_2y_1\\
\mathbb P(0,h) & x_1 & y_1x_1=h\,x_1y_2,\ \ y_2x_1=h\,x_1y_1\\
\mathbb W(0,h) & x_1 & y_1x_1=h\,x_1y_2,\ \ y_2x_1=h\,x_1y_1
\end{array}
\]
and $n$ $q$-commutes with the other $x$-generator. Let $x'$ be that generator and $t=n^{-1}x'$ in the quotient division ring $Q(B)$; $t$ is transcendental over $\kk$. Using the remaining mixed relations, $h$ cancels and
\begin{equation}\notag
\begin{array}{c|c|c|c}
 & y_1t & y_2t & (a,C)\\\hline
\mathbb X & ty_1+y_1 & ty_2+y_2 & (1,\ I)\\
\mathbb Y & -ty_1+fy_1+y_2 & -ty_2+y_1+fy_2 & \bigl(-1,\ \begin{psmallmatrix}f&1\\1&f\end{psmallmatrix}\bigr)\\
\mathbb D(\pm1,h) & p\,ty_1-p^{-1}y_2 & p^{-1}ty_2+p^{-1}y_1 & \bigl(p,\ p\begin{psmallmatrix}0&-1\\1&0\end{psmallmatrix}\bigr)\\
\mathbb M(0,h) & -ty_1+y_2 & -ty_2-y_1 & \bigl(-1,\ \begin{psmallmatrix}0&1\\-1&0\end{psmallmatrix}\bigr)\\
\mathbb P(0,h) & ty_1-y_1 & ty_2+y_2 & \bigl(1,\ \begin{psmallmatrix}-1&0\\0&1\end{psmallmatrix}\bigr)\\
\mathbb W(0,h) & -ty_1+y_2 & -ty_2+y_1 & \bigl(-1,\ \begin{psmallmatrix}0&1\\1&0\end{psmallmatrix}\bigr)
\end{array}
\end{equation}
As an example take $\mathbb W(0,h)$, where $y_1x_2=h(x_1y_1-x_2y_2)$ and $y_2x_2=h(x_1y_2-x_2y_1)$. From $y_2x_1=hx_1y_1$ we get $y_1x_1^{-1}=h^{-1}x_1^{-1}y_2$, so
\begin{align*}
y_1t&=h^{-1}x_1^{-1}y_2x_2=x_1^{-1}(x_1y_2-x_2y_1)=y_2-ty_1,
\end{align*}
and $y_2t=y_1-ty_2$ in the same way. In every row, $C$ has a nonzero eigenvalue when $a=1$, and two distinct eigenvalues when $a=-1$: they are $f\pm1$ for $\mathbb Y$, $\pm ip$ for $\mathbb D$, $\pm i$ for $\mathbb M$ and $\pm1$ for $\mathbb W$. If $B$ were PI, $Q(B)$ would be finite-dimensional over its center by Posner's theorem, and Lemma~\ref{lem:shiftmat} gives a contradiction.
\end{proof}

\begin{remark}\label{rem:LM}
The proof is an instance of the criterion of Leroy and Matczuk for Ore extensions \cite[Theorem 2.10]{LM}, which for automorphisms is due to J\o ndrup \cite{Jondrup}. By the table and Lemma~\ref{lem:fracind}, the subalgebra of $Q(B)$ generated by $A_y$ and $t$ is an Ore extension $A_y[t;\sigma,\delta]$; solving the table for $ty_j$ gives $ty_j=a\,y_jt-a\sum_lC_{jl}y_l$, so $\sigma(y_j)=a\,y_j$, $\delta(y_j)=-a\sum_lC_{jl}y_l$, and $\sigma$ is the identity on the center $Z$ of $A_y$. In characteristic zero such an extension is PI only if, after multiplying by a unit of $Q(A_y)$, $\delta$ becomes an inner derivation; in particular $\delta(Z)=0$. In every row, however, $\delta$ moves a central element of $A_y$. For instance, in $\mathbb W(0,h)$ we have $\sigma(y_1)=-y_1$ and $\delta(y_1)=y_2$, so
\[
\delta(y_1^2)=\delta(y_1)y_1+\sigma(y_1)\delta(y_1)=y_2y_1-y_1y_2=-2y_1y_2\ne0 .
\]
The elements $\vartheta$ of Lemma~\ref{lem:DGmodel} and $\Theta$ of \eqref{eq:Theta} play the opposite role: they are changes of variable that remove the derivation, as in the proof of \cite[Theorem 2.10]{LM} and in \cite{LM2}.
\end{remark}

\section{The monomial criterion}\label{sec:mono}

\subsection{Monomial data}

\begin{definition}\label{def:mono}
A \emph{monomial datum} for $B$ consists of bases $z_1,z_2$ of $\kk y_1\oplus\kk y_2$ and $e_1,e_2$ of $\kk x_1\oplus\kk x_2$, permutations $\sigma,\tau$ of $\{1,2\}$ and scalars $\lambda_{ab}\in\kk^\times$ with
\begin{equation}\label{eq:mono}
z_ae_b=\lambda_{ab}\,e_{\sigma(b)}z_{\tau(a)}\qquad(a,b\in\{1,2\}).
\end{equation}
It is of type $0$, I, II or III when neither, both, only $\tau$, or only $\sigma$ is the transposition. We call it \emph{admissible} if each plane is $\kk_{\pm1}$ (any basis), or the chosen basis of that plane is a coordinate basis of a quantum plane $\kk_q$.
\end{definition}

Put
\begin{align*}
\mathcal I_0&=\{\lambda_{11},\lambda_{12},\lambda_{21},\lambda_{22}\}, & \mathcal I_{\rm I}&=\{\lambda_{11}\lambda_{22},\ \lambda_{12}\lambda_{21}\},\\
\mathcal I_{\rm II}&=\Bigl\{\lambda_{11}\lambda_{21},\ \lambda_{12}\lambda_{22},\ \frac{\lambda_{11}}{\lambda_{12}},\ \frac{\lambda_{21}}{\lambda_{22}}\Bigr\}, & \mathcal I_{\rm III}&=\Bigl\{\lambda_{11}\lambda_{12},\ \lambda_{21}\lambda_{22},\ \frac{\lambda_{11}}{\lambda_{21}},\ \frac{\lambda_{12}}{\lambda_{22}}\Bigr\}.
\end{align*}

\begin{lemma}\label{lem:rescale}
Under $z_a\mapsto d_az_a$, $e_b\mapsto c_be_b$ the scalars become $\lambda'_{ab}=\lambda_{ab}\,d_ac_b/(d_{\tau(a)}c_{\sigma(b)})$, and the set $\mathcal I$ of the corresponding type does not change. One can choose $d_a,c_b\in\kk^\times$ so that every $\lambda'_{ab}$ lies in the subgroup generated by $\mathcal I$; the choice may require square roots, which exist since $\kk$ is algebraically closed.
\end{lemma}

\begin{proof}
The formula follows from \eqref{eq:mono}. In type II put $r=d_1/d_2$, so $\lambda'_{1b}=r\lambda_{1b}$ and $\lambda'_{2b}=r^{-1}\lambda_{2b}$; with $r=\lambda_{11}^{-1}$ the new scalars are $1$, $\lambda_{12}/\lambda_{11}$, $\lambda_{11}\lambda_{21}$ and $\lambda_{12}\lambda_{22}\cdot\lambda_{11}/\lambda_{12}$. In type I put also $u=c_1/c_2$; then $\lambda'_{11}=ru\lambda_{11}$, $\lambda'_{12}=ru^{-1}\lambda_{12}$, $\lambda'_{21}=r^{-1}u\lambda_{21}$, $\lambda'_{22}=(ru)^{-1}\lambda_{22}$, and choosing $ru=\lambda_{11}^{-1}$, $ru^{-1}=\lambda_{12}^{-1}$ (a square root is needed) gives $1,1,\lambda_{12}\lambda_{21},\lambda_{11}\lambda_{22}$. Type III is type II with the roles of the planes exchanged, and type $0$ needs nothing.
\end{proof}

\begin{lemma}\label{lem:plane}
Let $X=\kk_\varepsilon[x_1,x_2]$, with $\varepsilon=\pm1$,
and let $e_1,e_2$ be a basis of $X_1$. Put
\[
u=e_1^2,\qquad v=e_2^2,\qquad
a=e_1e_2,\qquad a'=e_2e_1.
\]
Then $u,v,a+a'$ and $aa'$ are central in $X$, and $a$ and $a'$
commute. For every positive even integer $n$, the algebra $X$
is a finite module over the central subalgebra
\[
S_n=\kk[u^n+v^n,(uv)^n,a^n+(a')^n].
\]
If $\varepsilon=-1$, then
\[
\kk[x_1^2,x_2^2]=\kk[u,v,a+a'].
\]
\end{lemma}

\begin{proof}
If $\varepsilon=1$, the elements $u^n$ and $v^n$ are integral
over $S_n$, since they satisfy
\[
T^2-(u^n+v^n)T+(uv)^n=0.
\]
Consequently, $u,v,e_1,e_2$ are integral over $S_n$,
so $X$ is a finite $S_n$-module.

Now assume $\varepsilon=-1$ and write
$e_b=\alpha_bx_1+\beta_bx_2$. Squares of linear forms are central, $u=\alpha_1^2x_1^2+\beta_1^2x_2^2$, $v=\alpha_2^2x_1^2+\beta_2^2x_2^2$, and $a+a'=2(\alpha_1\alpha_2x_1^2+\beta_1\beta_2x_2^2)$. These three elements span $\kk x_1^2+\kk x_2^2$. Indeed, if the vectors $(\alpha_1^2,\beta_1^2)$, $(\alpha_2^2,\beta_2^2)$, $(\alpha_1\alpha_2,\beta_1\beta_2)$ lay on a line, comparing the third with the first two would give $\alpha_1\beta_1(\alpha_1\beta_2-\alpha_2\beta_1)=\alpha_2\beta_2(\alpha_1\beta_2-\alpha_2\beta_1)=0$, so $\alpha_1\beta_1=\alpha_2\beta_2=0$; then $e_1,e_2$ are multiples of $x_1,x_2$ in some order, and $u,v$ are independent. So $\kk[u,v,a+a']=\kk[x_1^2,x_2^2]$. Next $aa'=e_1v e_1=uv=a'a$. Write $a=m+w$, $a'=m-w$ with $m$ central and $w$ a multiple of $x_1x_2$; then $a^n+(a')^n$ only contains even powers of $w$, and $w^2$ is central. Finally, $u^n,v^n$ are roots of $T^2-(u^n+v^n)T+(uv)^n$, so $u,v$ are integral over $S_n$.
Put $c=a+a'$ and $R=S_n[u,v]$. The power sum $a^n+(a')^n$ is a monic polynomial of degree $n$ in $c$ with coefficients in $\kk[uv]\subseteq R$.
Since $a^n+(a')^n\in S_n$, this gives a monic equation for $c$
over $R$. Thus $R[c]=\kk[u,v,c]=\kk[x_1^2,x_2^2]$ is finite
over $S_n$, by transitivity of finite integral extensions.
As $X$ is free of rank four over $\kk[x_1^2,x_2^2]$,
it is a finite $S_n$-module.
\end{proof}

\begin{theorem}\label{thm:mono}
Let $B$ have an admissible monomial datum. Then $B$ is PI if and only if the parameters of the two planes are roots of unity and $\mathcal I\subseteq\mi$, where $\mathcal I$ is the set of the corresponding type.
\end{theorem}

\begin{proof}
The planes are subalgebras, so their parameters must be roots of unity.

\emph{Necessity, products.} From \eqref{eq:mono}, $z_a^2e_b=\lambda_{ab}\lambda_{a,\sigma(b)}\,e_bz_{\tau(a)}^2$ in types $0$ and I, and $z_a^2e_b=\lambda_{ab}^2\,e_bz_{\tau(a)}^2$ in type II. In type $0$ put $\gamma=1$. In types I and II, where $\tau$ exchanges $z_1$ and $z_2$, let $\gamma$ be the scalar with $z_2^2z_1^2=\gamma\,z_1^2z_2^2$; it is $1$ in $\kk_{\pm1}$, where squares are central, and $q^4$ for a coordinate basis of $\kk_q$, so it is a root of unity. So
\[
(z_1^2z_2^2)\,e_b=\gamma\,\Lambda_b\,e_b\,(z_1^2z_2^2),\qquad \Lambda_b=\begin{cases}\lambda_{11}\lambda_{12}\lambda_{21}\lambda_{22}&\text{type I},\\ (\lambda_{1b}\lambda_{2b})^2&\text{types }0,\ \mathrm{II},\end{cases}
\]
and Lemma~\ref{lem:qplane} gives $\Lambda_b\in\mi$. In type $0$, Lemma~\ref{lem:qplane} applied to $e_b$ and $z_a$ gives $\lambda_{ab}\in\mi$ directly.

\emph{Necessity, ratios.} Work in $Q(B)$. In type II, the relation $z_{\tau(a)}e_1=\lambda_{\tau(a),1}e_1z_a$ gives $z_ae_1^{-1}=\lambda_{\tau(a),1}^{-1}e_1^{-1}z_{\tau(a)}$, so for $t=e_1^{-1}e_2$,
\[
z_atz_a^{-1}=\frac{\lambda_{\tau(a),2}}{\lambda_{\tau(a),1}}\,t .
\]
In type I the same computation with $t'=e_2^{-1}e_1$ gives $z_at'z_a^{-1}=c_a\,t'^{-1}$ with $c_1=\lambda_{21}/\lambda_{22}$ and $c_2=\lambda_{11}/\lambda_{12}$, hence
\[
(z_1z_2)\,t'\,(z_1z_2)^{-1}=\frac{c_2}{c_1}\,t'=\frac{\lambda_{11}\lambda_{22}}{\lambda_{12}\lambda_{21}}\,t' .
\]
Both $t$ and $t'$ lie outside $\kk$, and Lemma~\ref{lem:shift} shows that these ratios are roots of unity. Together with the products this gives $\mathcal I\subseteq\mi$ in types I and II; type III follows by exchanging the roles of the planes, that is, by using the relations $e_cz_d=\lambda_{\tau(d),\sigma(c)}^{-1}z_{\tau(d)}e_{\sigma(c)}$ obtained from \eqref{eq:mono}.

\emph{Sufficiency.} By Lemma~\ref{lem:rescale} we may assume all $\lambda_{ab}\in\mi$. Choose a positive even integer $n$ such that $n/2$ is divisible by the orders of all $\lambda_{ab}$, of all products of two of them, and of the plane parameters. From \eqref{eq:mono},
\[
z_ae_b^n=(\lambda_{ab}\lambda_{\tau(a),b})^{n/2}\,e_{\sigma(b)}^n\,z_a=e_{\sigma(b)}^nz_a,\qquad z_a\,e_ce_d=\kappa_{a,cd}\,e_{\sigma(c)}e_{\sigma(d)}\,z_a,\quad\kappa_{a,cd}=\lambda_{ac}\lambda_{\tau(a),d},
\]
and $\kappa_{a,cd}^n=1$ by the choice of $n$, so that $z_a(e_ce_d)^n=(e_{\sigma(c)}e_{\sigma(d)})^nz_a$ for $\{c,d\}=\{1,2\}$. So conjugation by each $z_a$ permutes $\{e_1^n,e_2^n\}$ and $\{(e_1e_2)^n,(e_2e_1)^n\}$. If $A_x=\kk_{\pm1}$ put
\[
Z_x=\kk\bigl[e_1^{2n}+e_2^{2n},\ (e_1e_2)^{n}+(e_2e_1)^{n},\ e_1^{2n}e_2^{2n}\bigr],
\]
and if $e_1,e_2$ is a coordinate basis of $\kk_q$ put $Z_x=\kk[e_1^n+e_2^n,\ e_1^ne_2^n]$. In both cases $Z_x$ commutes with $A_x$ and with $z_1,z_2$, so it is central in $B$, and $A_x$ is a finite $Z_x$-module, by Lemma~\ref{lem:plane} in the first case and clearly in the second. Rewriting \eqref{eq:mono} as $e_cz_d=\lambda_{\tau(d),\sigma(c)}^{-1}z_{\tau(d)}e_{\sigma(c)}$, the same construction on the other side gives a central $Z_y$ over which $A_y$ is finite. Since $B=A_xA_y$, the algebra $B$ is a finite module over the commutative central subalgebra $Z_xZ_y$, hence PI.
\end{proof}

\subsection{The fourteen monomial families}

\begin{corollary}\label{cor:mono}
The families $\mathbb B$, $\mathbb E$, $\mathbb J$, $\mathbb K$, $\mathbb L$, $\mathbb N$, $\mathbb S$, $\mathbb T$, $\mathbb U$, $\mathbb Z$, and $\mathbb M$, $\mathbb O$, $\mathbb P$, $\mathbb W$ for $f\ne0$, have the admissible monomial data of Table~\ref{tab:mono}, and satisfy the criteria of Table~\ref{tab:main}.
\end{corollary}

{\small
\begin{longtable}{|c|c|c|c|c|c|}
\caption{Monomial data; $s^2=f$, the first sign is index $1$, and the scalars are divided by $h$ (except for $\mathbb N$). For $\mathbb B$, $\mathbb E$, $\mathbb J$ the data are the same for both values of $p$ or $q$.}\label{tab:mono}\\
\hline
 & $z_1,z_2$ & $e_1,e_2$ & type & $(\lambda_{11},\lambda_{12},\lambda_{21},\lambda_{22})$ & $\mathcal I\subseteq\mi$ iff\\\hline
\endfirsthead
\hline
 & $z_1,z_2$ & $e_1,e_2$ & type & $(\lambda_{11},\lambda_{12},\lambda_{21},\lambda_{22})$ & $\mathcal I\subseteq\mi$ iff\\\hline
\endhead
$\mathbb B$ & $y_1,\ y_2$ & $x_1,\ x_2$ & I & $(1,1,-1,1)$ & $h\in\mi$\\
$\mathbb E$ & $y_1,\ y_2$ & $x_1\pm ix_2$ & I & $(1+i,\ 1-i,\ -1+i,\ -1-i)$ & $2h^2\in\mi$\\
$\mathbb J$ & $y_1\pm iy_2$ & $x_1,\ x_2$ & I & $(1+i,\ -1+i,\ 1-i,\ -1-i)$ & $2h^2\in\mi$\\
$\mathbb K$ & $y_1\pm s^{-1}y_2$ & $x_1,\ x_2$ & 0 & $(1,\ s,\ 1,\ -s)$ & $f,h\in\mi$\\
$\mathbb L$ & $y_1,\ y_2$ & $x_1,\ x_2$ & II & $(f,\ 1,\ f,\ 1)$ & $f,h\in\mi$\\
$\mathbb M$ & $y_1\pm iy_2$ & $x_1\pm is^{-1}x_2$ & I & $(i(s+1),\ -i(s-1),\ i(s-1),\ -i(s+1))$ & $h(s\pm1)\in\mi$\\
$\mathbb N$ & $y_1\pm y_2$ & $x_1,\ x_2$ & III & $(f-g,\ f+g,\ -f-g,\ g-f)$ & $f\pm g\in\mi$\\
$\mathbb O$ & $y_1\pm iy_2$ & $x_1\pm isx_2$ & I & $(1-s,\ 1+s,\ 1+s,\ 1-s)$ & $h(1\pm s)\in\mi$\\
$\mathbb P$ & $y_1,\ y_2$ & $x_1\pm sx_2$ & II & $(1+s,\ 1-s,\ 1-s,\ 1+s)$ & $h(1\pm s)\in\mi$\\
$\mathbb S$ & $y_1\pm y_2$ & $x_1\pm x_2$ & 0 & $(2,-2,-2,-2)$ & $2h\in\mi$\\
$\mathbb T$ & $y_1\pm iy_2$ & $x_1\pm x_2$ & II & $(2i,\ -2,\ -2i,\ -2)$ & $2h\in\mi$\\
$\mathbb U$ & $y_1\pm y_2$ & $x_1\pm ix_2$ & III & $(2i,\ -2i,\ -2,\ -2)$ & $2h\in\mi$\\
$\mathbb W$ & $y_1\pm y_2$ & $x_1\pm isx_2$ & III & $(i(s-i),\ -i(s+i),\ i(s+i),\ -i(s-i))$ & $h(s\pm i)\in\mi$\\
$\mathbb Z$ & $y_1\pm is^{-1}y_2$ & $x_1\pm x_2$ & II & $(i(s-i),\ -i(s+i),\ -i(s+i),\ i(s-i))$ & $h(s\pm i)\in\mi$\\
\hline
\end{longtable}
}

\begin{proof}
Each entry is a direct substitution into the mixed relations. For instance, in $\mathbb O$ the relations are $y_1x_1=h(x_1y_1+fx_2y_2)$, $y_1x_2=h(-x_2y_1+x_1y_2)$, $y_2x_1=h(fx_2y_1-x_1y_2)$, $y_2x_2=h(x_1y_1+x_2y_2)$, and with $z_1=y_1+iy_2$, $e_1=x_1+isx_2$,
\begin{align*}
z_1e_1&=y_1x_1+is\,y_1x_2+i\,y_2x_1-s\,y_2x_2\\
&=h\bigl((1-s)x_1y_1+(is^2-is)x_2y_1+(is-i)x_1y_2+(s^2-s)x_2y_2\bigr)\\
&=h(1-s)\,(x_1-isx_2)(y_1-iy_2)=h(1-s)\,e_2z_2 .
\end{align*}
The planes are $\kk_{\pm1}$ except in $\mathbb B$ ($A_x=\kk_p$, $A_y=\kk_p$, coordinate bases), $\mathbb E$ ($A_y=\kk_p$, coordinate), $\mathbb J$ ($A_x=\kk_q$, coordinate) and $\mathbb K$, $\mathbb L$ ($A_x=\kk_q$, coordinate), so the data are admissible. The bases are bases because $s\ne0$, which is where $f\ne0$ is needed, and the $\lambda_{ab}$ are nonzero because $f\ne1$ in $\mathbb M$, $\mathbb O$, $\mathbb P$, $f\ne-1$ in $\mathbb W$, $\mathbb Z$ and $f^2\ne g^2$ in $\mathbb N$. For the last column compute $\mathcal I$: in $\mathbb O$, for example, $\lambda_{11}\lambda_{22}=h^2(1-s)^2$ and $\lambda_{12}\lambda_{21}=h^2(1+s)^2$; in $\mathbb E$, $\lambda_{11}\lambda_{22}=-2ih^2$ and $\lambda_{12}\lambda_{21}=2ih^2$; in $\mathbb P$, $\lambda_{11}\lambda_{21}=h^2(1-f)$ and $\lambda_{11}/\lambda_{12}=(1+s)/(1-s)$. Theorem~\ref{thm:mono} finishes the proof. The plane parameters are $\pm1$ or $\pm i$ except in $\mathbb K$ and $\mathbb L$, where the condition $q\in\mi$ has to be added to the last column.
\end{proof}

For $\mathbb K$ the criterion, with the PI degree, is in \cite{GV-K}, and also in \cite[Proposition~3.1 and Theorem~3.2]{Rub-ozone}. For $\mathbb O$ the criterion is \cite[Proposition~4.1]{Rub-ozone}. The algebras $\mathbb E$ and $\mathbb J$ are $\Sigma$-$M$-dual \cite[Subcase 4.3.2]{ZZ09}: the matrix $M$ of $\mathbb E$ is the matrix $\Sigma$ of $\mathbb J$. By \cite[Proposition 3.8]{ZZ09}, or directly on the four mixed relations, $x_i\mapsto y_i$, $y_i\mapsto x_i$ defines an isomorphism
\begin{equation}\label{eq:EJ}
\mathbb J(p^{-1},h)\cong\mathbb E(p,h)^{\mathrm{op}},
\end{equation}
with the same $h$. Under \eqref{eq:EJ} the data of $\mathbb E$ in Table~\ref{tab:mono} go to those of $\mathbb J$, with $(\lambda_{ab})$ transposed.

Several of the bases of Table~\ref{tab:mono} already appear in \cite[Section 5]{ZZ09}, where they serve to prove the noetherian property: $\mathbb B$, $\mathbb E$, $\mathbb P$ have $\Sigma_{11}=\Sigma_{22}=0$ \cite[Lemma 5.11, Proposition 5.12]{ZZ09}, the elements $x_1\pm x_2$ skew-commute with $\kk y_1+\kk y_2$ in $\mathbb S$, $\mathbb T$, $\mathbb Z$ \cite[Lemma 5.15, Proposition 5.16]{ZZ09}, and $y_1+iy_2$ is used for $\mathbb M$ and $\mathbb O$ in the proof of \cite[Proposition 5.14]{ZZ09}. What we need here is the four relations \eqref{eq:mono} at once, in $B$ itself.

\section{Families with a normal generator}\label{sec:normal}

Each of the families $\mathbb D$, $\mathbb G$, $\mathbb K$, $\mathbb L$, $\mathbb Q$, $\mathbb V$ has a normal generator of degree one, as one checks on the relations of \cite[Section 4]{ZZ09}. The families $\mathbb K$ and $\mathbb L$ were treated in Corollary~\ref{cor:mono}. For the other four, let $n$ be the normal generator, $x'$ the other $x$-generator and $t=n^{-1}x'$. As in Proposition~\ref{prop:noPI}, the twist cancels and
\begin{equation}\label{eq:Mt}
\begin{pmatrix}y_1\\y_2\end{pmatrix}t=M(t)\begin{pmatrix}y_1\\y_2\end{pmatrix},\qquad\text{hence}\qquad \begin{pmatrix}y_1\\y_2\end{pmatrix}F(t)=F\bigl(M(t)\bigr)\begin{pmatrix}y_1\\y_2\end{pmatrix}\quad\text{for all }F\in\kk[t],
\end{equation}
where $M(t)$ is a $2\times2$ matrix with entries in $\kk[t]$. Central elements come from polynomials $F$ with $F(M(t))=F(t)I$.

\subsection{The families \texorpdfstring{$\mathbb D$}{D} and \texorpdfstring{$\mathbb G$}{G}}

The relations of $\mathbb D(p,h)$ are $x_2x_1=-x_1x_2$, $y_2y_1=p\,y_1y_2$ and
\begin{align*}
y_1x_1&=-hp\,x_1y_1, & y_1x_2&=h(-p^2x_2y_1+x_1y_2),\\
y_2x_1&=hp\,x_1y_2, & y_2x_2&=h(x_1y_1+x_2y_2).
\end{align*}
In $\mathbb G(p,f,h)$ we replace the generator $x_2$ of \cite{ZZ09} by $x_2+\frac1{p-1}x_1$, as in \cite{GV-G}; then the relations are $x_2x_1=x_1x_2$, $y_2y_1=p\,y_1y_2$ and
\begin{align*}
y_1x_1&=hp\,x_1y_1, & y_1x_2&=h(p^2x_2y_1+x_1y_2),\\
y_2x_1&=hp\,x_1y_2, & y_2x_2&=h(f\,x_1y_1+x_2y_2).
\end{align*}
In both, $n=x_1$, $t=x_1^{-1}x_2$ and
\begin{equation}\label{eq:MDG}
M_{\mathbb D}(t)=\begin{pmatrix}pt&-p^{-1}\\ p^{-1}&p^{-1}t\end{pmatrix},\qquad M_{\mathbb G}(t)=\begin{pmatrix}pt&p^{-1}\\ fp^{-1}&p^{-1}t\end{pmatrix}.
\end{equation}

\begin{lemma}\label{lem:potY}
Let $y_2y_1=p\,y_1y_2$ and $y_1t=p\,ty_1+c_{12}y_2$, $y_2t=p^{-1}ty_2+c_{21}y_1$ with $c_{12},c_{21}\in\kk$. Then for $k\ge1$
\begin{align*}
y_1^kt&=p^k\,ty_1^k+c_{12}[k]_{p^2}\,y_1^{k-1}y_2, & y_2^kt&=p^{-k}\,ty_2^k+c_{21}\,p^{k-1}[k]_{p^{-2}}\,y_1y_2^{k-1},
\end{align*}
where $[k]_q=1+q+\dots+q^{k-1}$. In particular, if $p^N=1$ and $p^2\ne1$, then $y_1^N$ and $y_2^N$ commute with $t$.
\end{lemma}

\begin{proof}
Induction on $k$. Using $y_2y_1^k=p^ky_1^ky_2$,
\begin{align*}
y_1^{k+1}t&=y_1\bigl(p^kty_1^k+c_{12}[k]_{p^2}y_1^{k-1}y_2\bigr)=p^k(p\,ty_1+c_{12}y_2)y_1^k+c_{12}[k]_{p^2}y_1^ky_2\\
&=p^{k+1}ty_1^{k+1}+c_{12}\bigl(p^{2k}+[k]_{p^2}\bigr)y_1^ky_2,
\end{align*}
and the other formula is proved in the same way. If $p^N=1$ and $p^2\ne1$ then $[N]_{p^{\pm2}}=0$.
\end{proof}

For $\alpha\in\kk$ let $D_n(X,\alpha)\in\kk[X]$ be the Dickson polynomial, determined by $D_n(w+\alpha w^{-1},\alpha)=w^n+\alpha^nw^{-n}$. It is monic of degree $n$ and has the parity of $n$.

\begin{lemma}\label{lem:dickson}
Let $M(t)=\begin{psmallmatrix}pt&c_{12}\\ c_{21}&p^{-1}t\end{psmallmatrix}$ with $p^n=1$, $p^2\ne1$, and $\alpha=-c_{12}c_{21}\,p^2(p^2-1)^{-2}$. Then $D_n(M(t),\alpha)=D_n(t,\alpha)I$.
\end{lemma}

\begin{proof}
Write $t=w+\alpha w^{-1}$ in an extension of $\kk(t)$. The elements $\mu_1=pw+\alpha(pw)^{-1}$ and $\mu_2=p^{-1}w+\alpha p\,w^{-1}$ satisfy
\begin{align*}
\mu_1+\mu_2&=(p+p^{-1})t, & \mu_1\mu_2&=t^2+\alpha(p-p^{-1})^2=t^2-c_{12}c_{21},
\end{align*}
so they are the eigenvalues of $M(t)$. They are distinct because $p^2\ne1$, and $M(t)$ is diagonalizable. Since $p^n=1$, $D_n(\mu_1,\alpha)=(pw)^n+\alpha^n(pw)^{-n}=D_n(t,\alpha)$, and the same holds for $\mu_2$.
\end{proof}

\begin{theorem}\label{thm:DGPI}
$\mathbb D(p,h)$ is PI if and only if $p,h\in\mi$ and $p\ne\pm1$. $\mathbb G(p,f,h)$ is PI if and only if $p,h\in\mi$.
\end{theorem}

\begin{proof}
Necessity. The plane $A_y=\kk_p[y_1,y_2]$ forces $p\in\mi$, Lemma~\ref{lem:qplane} applied to $x_1$ and $y_1$ gives $hp\in\mi$, and $\mathbb D(\pm1,h)$ is excluded by Proposition~\ref{prop:noPI}.

Sufficiency. By Theorem~\ref{thm:twist} we may take $h=1$. Let $N=\ord p$, so $N\ge3$.

In $\mathbb D(p,1)$ let $\ell=\operatorname{lcm}(2,N)$; by \eqref{eq:MDG}, Lemmas~\ref{lem:potY} and~\ref{lem:dickson} apply with $c_{12}c_{21}=-p^{-2}$, so $\alpha=(p^2-1)^{-2}$. The elements $y_1^\ell$ and $y_2^N$ commute with $t$ by Lemma~\ref{lem:potY}, with $x_1$ because $(-p)^\ell=p^\ell=1$, hence with $x_2=x_1t$, and they are central in $A_y$. The element $x_1^\ell$ commutes with $y_1,y_2$ and, $\ell$ being even, with $x_2$. Finally let
\[
\Phi=x_1^\ell\,D_N(t,\alpha)^{\ell/N}.
\]
The polynomial $D_N(t,\alpha)^{\ell/N}$ is even in $t$, so it commutes with $x_1$, which anticommutes with $t$; it commutes with $y_1,y_2$ by \eqref{eq:Mt} and Lemma~\ref{lem:dickson}. So $\Phi$ is central. It lies in $\kk[x_1^2,x_2^2]$ and has the form $\pm x_2^\ell+{}$terms of lower degree in $x_2$. Hence $A_x$ is a finite module over $\kk[x_1^\ell,\Phi]$, $A_y$ is finite over $\kk[y_1^\ell,y_2^N]$, and $B=A_xA_y$ is a finite module over a central subalgebra.

In $\mathbb G(p,f,1)$, \eqref{eq:MDG} gives $c_{12}c_{21}=fp^{-2}$ and $\alpha=-f(p^2-1)^{-2}$. In the same way $x_1^N$, $y_1^N$, $y_2^N$ and $x_1^ND_N(t,\alpha)=D_N(x_2,\alpha x_1^2)$ are central, and $B$ is a finite module over the subalgebra they generate.
\end{proof}

The PI degrees of $\mathbb D$ and $\mathbb G$ are computed in Proposition~\ref{prop:DGdeg}, from a quantum torus with the same quotient division ring (Lemma~\ref{lem:DGmodel}); for $\mathbb G$ see also \cite{GV-G}. For $\mathbb K$, a crossed product of index two over a quantum torus, see \cite[Theorem 2.5]{GV-K}.

\subsection{The families \texorpdfstring{$\mathbb Q$}{Q} and \texorpdfstring{$\mathbb V$}{V}}

For $h=1$ the relations are
\begin{align*}
\mathbb Q:\quad & y_1x_1=x_1y_2, && y_1x_2=x_1y_1+x_2y_1+x_1y_2,\\
& y_2x_1=-x_1y_1, && y_2x_2=x_1y_1-x_1y_2+x_2y_2,\\
\mathbb V:\quad & y_1x_1=x_2y_1+x_1y_2, && y_1x_2=x_2y_1,\\
& y_2x_1=-x_1y_1+x_2y_1, && y_2x_2=x_2y_2,
\end{align*}
with $y_2y_1=-y_1y_2$ in both, $x_2x_1=-x_1x_2$ in $\mathbb Q$ and $x_2x_1=x_1x_2$ in $\mathbb V$.

\begin{proposition}\label{prop:QV}
$\mathbb Q(h)$ and $\mathbb V(h)$ are PI if and only if $h\in\mi$. For $h=1$ the elements
\begin{align*}
\mathbb Q(1):&\quad y_1^2+y_2^2,\quad y_1^2y_2^2,\quad x_1^2+x_2^2,\quad x_1^2x_2^2+\tfrac12x_2^4,\\
\mathbb V(1):&\quad y_1^2+y_2^2,\quad y_1^2y_2^2,\quad x_2,\quad x_1^2(x_1-x_2)^2
\end{align*}
are central, and the algebra is a finite module over the subalgebra they generate.
\end{proposition}

\begin{proof}
Necessity. In $\mathbb Q(h)$, $y_1x_1=hx_1y_2$ and $y_2x_1=-hx_1y_1$, so $y_1x_1^2=-h^2x_1^2y_1$; in $\mathbb V(h)$, $y_1x_2=hx_2y_1$. Lemma~\ref{lem:qplane} gives $h^2\in\mi$, respectively $h\in\mi$.

Sufficiency. By Theorem~\ref{thm:twist} it suffices to treat $h=1$. Two applications of the relations give, in $\mathbb Q(1)$,
\begin{align*}
y_1^2x_1&=x_1y_2^2, & y_1^2x_2&=x_1y_1^2-2x_1y_1y_2+x_1y_2^2+x_2y_1^2,\\
y_2^2x_1&=x_1y_1^2, & y_2^2x_2&=-x_1y_1^2+2x_1y_1y_2-x_1y_2^2+x_2y_2^2,
\end{align*}
and in $\mathbb V(1)$,
\begin{align*}
y_1^2x_1&=x_1y_2^2+x_2y_1^2+x_2y_1y_2, & y_1^2x_2&=x_2y_1^2,\\
y_2^2x_1&=x_1y_1^2-x_2y_1^2-x_2y_1y_2, & y_2^2x_2&=x_2y_2^2 .
\end{align*}
Adding, $y_1^2+y_2^2$ commutes with $x_1,x_2$; it commutes with $y_1,y_2$ because squares are central in $\kk_{-1}[y_1,y_2]$. For $y_1^2y_2^2$ one uses the same formulas twice; for instance in $\mathbb V(1)$, since $y_1^2$ is central in $A_y$,
\begin{align*}
y_1^2y_2^2x_1&=y_1^2\bigl(x_1y_1^2-x_2y_1^2-x_2y_1y_2\bigr)=(x_1y_2^2+x_2y_1^2+x_2y_1y_2)y_1^2-x_2y_1^4-x_2y_1^3y_2=x_1y_1^2y_2^2 .
\end{align*}

For the $x$-elements let $n=x_1$, $t=x_1^{-1}x_2$ in $\mathbb Q$ and $n=x_2$, $t=x_2^{-1}x_1$ in $\mathbb V$. In \eqref{eq:Mt},
\[
M_{\mathbb Q}(t)=\begin{pmatrix}-1&1-t\\1+t&1\end{pmatrix},\qquad M_{\mathbb V}(t)=\begin{pmatrix}1&t\\1-t&0\end{pmatrix}.
\] By Cayley--Hamilton, $M_{\mathbb Q}^2=(2-t^2)I$ and $M_{\mathbb V}^2-M_{\mathbb V}=(t-t^2)I$. So $F_{\mathbb Q}(t)=t^2(2-t^2)$ and $F_{\mathbb V}(t)=(t^2-t)^2$ satisfy $F(M(t))=F(t)I$, that is, they commute with $y_1,y_2$. In $\mathbb Q$, $t$ anticommutes with $x_1$, $x_1^4$ commutes with $y_1,y_2$, and
\[
x_1^4F_{\mathbb Q}(t)=-2x_1^2x_2^2-x_2^4,
\]
which is central since it lies in $\kk[x_1^2,x_2^2]$. In the same way $1-M_{\mathbb Q}^2=-(1-t^2)I$ and $y_jx_1^2=-x_1^2y_j$, so $x_1^2(1-t^2)=x_1^2+x_2^2$ commutes with $y_1,y_2$, and it is central in $A_x$. In $\mathbb V$, $x_2$ is central and $x_2^4F_{\mathbb V}(t)=x_1^2(x_1-x_2)^2$.

Finiteness: $\kk[x_1^2,x_2^2]$ is integral over $\kk[x_1^2+x_2^2,\ x_1^2x_2^2+\frac12x_2^4]$, since $x_2^2$ is a root of $\frac12T^2-(x_1^2+x_2^2)T+(x_1^2x_2^2+\frac12x_2^4)$, and $\kk[x_1,x_2]$ is integral over $\kk[x_2,x_1^2(x_1-x_2)^2]$. With the $y$-elements and the basis this gives a finite module over a central subalgebra.
\end{proof}

\subsection{The family \texorpdfstring{$\mathbb L$}{L} and its PI degree}

The relations of $\mathbb L(q,f,h)$ are
\begin{align*}
x_2x_1&=q\,x_1x_2, & y_1x_1&=hf\,x_1y_2, & y_1x_2&=h\,x_2y_2,\\
y_2y_1&=-y_1y_2, & y_2x_1&=hf\,x_1y_1, & y_2x_2&=h\,x_2y_1 .
\end{align*}
Both $x_1$ and $x_2$ are normal and exchange $y_1$ and $y_2$ up to scalars.

\begin{theorem}\label{thm:L}
$\mathbb L(q,f,h)$ is PI if and only if $q,f,h\in\mi$. In that case let $N=|\langle-1,q,f,h\rangle|$ and $N_0=|\langle-1,q,f\rangle|$. Then
\[
\PIdeg\mathbb L(q,f,h)=\begin{cases}2N, & N/N_0\ \text{odd},\\ N, & N/N_0\ \text{even}.\end{cases}
\]
\end{theorem}

\begin{proof}
The criterion is Corollary~\ref{cor:mono}. Put $u=x_1^2$ and $w=x_1x_2$. From the relations,
\begin{align}
y_au&=(hf)^2\,uy_a, & y_aw&=h^2f\,wy_a, & wu&=q^2\,uw, & y_2y_1&=-y_1y_2 ,\label{eq:Ltorus}
\end{align}
and the monomials $y_1^cy_2^du^jw^l$ are independent because $u^jw^l$ is a nonzero multiple of $x_1^{2j+l}x_2^l$. So $y_1,y_2,u,w$ generate a quantum affine space; let $T$ be its quantum torus inside $Q(\mathbb L)$. Conjugation by $x_1$ sends $y_1\mapsto(hf)^{-1}y_2$, $y_2\mapsto(hf)^{-1}y_1$, $u\mapsto u$, $w\mapsto q^{-1}w$, so $x_1$ normalizes $Q(T)$, and $x_1^2=u\in T$, $x_2=x_1^{-1}w$. Since $x_1$ normalizes $Q(T)$ and $x_1^2\in T$, the sum $Q(T)+Q(T)x_1$ is a ring; it is a domain of dimension at most two over $Q(T)$, hence a division ring, and it contains $\mathbb L$. So $Q(\mathbb L)=Q(T)+Q(T)x_1$.

The element $y_1^N$ is central in $Q(T)$ by \eqref{eq:Ltorus}, and $x_1y_1^Nx_1^{-1}=y_2^N\ne y_1^N$. So $x_1\notin Q(T)$, the sum is direct, and $\Ad_{x_1}$ acts on $Z(Q(T))$ with order exactly two, its square being $\Ad_u$. If $t_0+t_1x_1$ is central, commuting it with $y_1^N$ gives $t_1(y_2^N-y_1^N)=0$, so $t_1=0$; hence $Z(Q(\mathbb L))=Z(Q(T))^{\langle\Ad_{x_1}\rangle}$ and
\[
[Q(\mathbb L):Z(Q(\mathbb L))]=2\cdot(\PIdeg T)^2\cdot 2 .
\]
It remains to compute $\PIdeg T$. Let $\zeta$ be a primitive $N$th root of unity with $q=\zeta^\theta$, $f=\zeta^\varphi$, $h=\zeta^\eta$, so that $\zeta^{N/2}=-1$ and $\gcd(\theta,\varphi,\eta,N/2)=1$. In the basis $w,y_1,y_1^{-1}y_2,u$ of the lattice of monomials, \eqref{eq:Ltorus} gives the commutation matrix with upper entries
\[
-(2\eta+\varphi),\quad 0,\quad 2\theta,\quad N/2,\quad 2\eta+2\varphi,\quad 0,
\]
whose Pfaffian is $\theta N$. Let $g$ be the gcd of the entries; it divides $N/2$. By the normal form of alternating integer matrices \cite[Chapter IV]{New}, as in \cite[Lemma 2.4]{GV-K}, $T$ has PI degree
\[
\frac Ng\cdot\frac{N}{\gcd(N,\theta N/g)}=\frac Ng\cdot\frac{g}{\gcd(g,\theta)}=\frac{N}{\gcd(\theta,N/2,\varphi,2\eta)} .
\]
The last gcd is $1$ or $2$, and it is $2$ exactly when $\theta$, $\varphi$ and $N/2$ are even, that is, when $-1,q,f$ lie in the subgroup of index two of $\langle-1,q,f,h\rangle$. This is the statement.
\end{proof}

For instance $\PIdeg\mathbb L(1,1,1)=4$, $\PIdeg\mathbb L(1,1,i)=4$ and $\PIdeg\mathbb L(i,1,1)=8$. In $\mathbb K$, with the same parameters, the PI degree is always $2N$ \cite[Theorem 2.5]{GV-K}; in $\mathbb L$ the exchange of $y_1,y_2$ by both $x$-generators can halve it.

\section{Normal subalgebras: the families \texorpdfstring{$\mathbb F$}{F}, \texorpdfstring{$\mathbb I$}{I} and \texorpdfstring{$\mathbb C$}{C}}\label{sec:CFI}

Neither method applies here. Zhang and Zhang proved these three families strongly noetherian with normal elements of degree two and three \cite[Propositions 5.13 and 5.17]{ZZ09}, and the same elements decide the PI property. In \cite[Proposition 5.13]{ZZ09} the family $\mathbb F$ is reduced to $\mathbb I$ through their $\Sigma$-$M$-duality; Lemma~\ref{lem:Fnormal} below is the counterpart for $\mathbb F$ of the computation made there for $\mathbb I$.

\begin{lemma}\label{lem:rank4}
Let $R=\kk_{-1}[a,b]$, $\omega_1=ab$ and $\omega_2=a^2+b^2$. Then $\omega_1$ and $\omega_2$ commute and are algebraically independent, and $R$ is free as a left and as a right $\kk[\omega_1,\omega_2]$-module with basis $1,a,b,a^2$.
\end{lemma}

\begin{proof}
$a^2,b^2$ are central, so $\omega_1$ commutes with $\omega_2$, and $\omega_1^2=-a^2b^2$. Since $\omega_2$ and $\omega_1^2$ are algebraically independent in $\kk[a^2,b^2]$, so are $\omega_1$ and $\omega_2$. From $a^3=a\omega_2+b\omega_1=\omega_2a-\omega_1b$, $b^3=b\omega_2-a\omega_1=\omega_2b+\omega_1a$ and $a^4=\omega_2a^2+\omega_1^2$, together with $ab=\omega_1$ and $ba=-\omega_1$, the four elements generate $R$ on either side, and the surjection $\kk[\omega_1,\omega_2]\oplus\kk[\omega_1,\omega_2](-1)^2\oplus\kk[\omega_1,\omega_2](-2)\to R$ is an isomorphism because both sides have Hilbert series $(1-t)^{-2}$.
\end{proof}

\subsection{The family \texorpdfstring{$\mathbb F$}{F}}

The relations of $\mathbb F(p,h)$, $p^2=-1$, are $x_2x_1=-x_1x_2$, $y_2y_1=py_1y_2$ and
\begin{align*}
y_1x_1&=h(-x_1y_1-px_2y_1+x_1y_2-x_2y_2), & y_1x_2&=h(-px_1y_1+x_2y_1+x_1y_2+x_2y_2),\\
y_2x_1&=h(-px_1y_1+px_2y_1+px_1y_2+x_2y_2), & y_2x_2&=h(-px_1y_1-px_2y_1+x_1y_2-px_2y_2).
\end{align*}

\begin{lemma}\label{lem:Fnormal}
Let $\omega_1=x_1x_2$, $\omega_2=x_1^2+x_2^2$ and $\lambda=-2ph^2$. Then, for $j=1,2$,
\begin{equation}\label{eq:Fnormal}
y_j\omega_1=\lambda\,\omega_1y_j,\qquad y_j\omega_2=\lambda\,\omega_2y_j .
\end{equation}
\end{lemma}

\begin{proof}
Apply the mixed relations twice and use $x_2x_1=-x_1x_2$, $p^2=-1$. For $y_1\omega_1$ the four terms of $y_1x_1$ give
\begin{align*}
h^{-2}y_1x_1x_2&=\bigl(p\,x_1^2y_1-x_1x_2y_1-x_1^2y_2-x_1x_2y_2\bigr)\\
&\quad+\bigl(x_1x_2y_1-p\,x_2^2y_1+p\,x_1x_2y_2-p\,x_2^2y_2\bigr)\\
&\quad+\bigl(-p\,x_1^2y_1-p\,x_1x_2y_1+x_1^2y_2-p\,x_1x_2y_2\bigr)\\
&\quad+\bigl(-p\,x_1x_2y_1+p\,x_2^2y_1+x_1x_2y_2+p\,x_2^2y_2\bigr)\\
&=-2p\,x_1x_2y_1 .
\end{align*}
In the same way
\begin{align*}
h^{-2}y_1x_1^2&=(1-p)x_1^2y_1+(p-1)x_1^2y_2+2(p+1)x_1x_2y_2-(p+1)x_2^2y_1+(p-1)x_2^2y_2,\\
h^{-2}y_1x_2^2&=-(p+1)x_1^2y_1-(p-1)x_1^2y_2-2(p+1)x_1x_2y_2-(p-1)x_2^2y_1-(p-1)x_2^2y_2,\\
h^{-2}y_2x_1^2&=(p+1)x_1^2y_1-(p+1)x_1^2y_2+2(p-1)x_1x_2y_1+(p+1)x_2^2y_1-(p-1)x_2^2y_2,\\
h^{-2}y_2x_2^2&=-(p+1)x_1^2y_1-(p-1)x_1^2y_2-2(p-1)x_1x_2y_1-(p+1)x_2^2y_1-(p+1)x_2^2y_2,\\
h^{-2}y_2x_1x_2&=-2p\,x_1x_2y_2,
\end{align*}
and adding the pairs gives $y_j\omega_2=-2ph^2\,\omega_2y_j$.
\end{proof}

\begin{theorem}\label{thm:F}
$\mathbb F(p,h)$ is PI if and only if $2h^2\in\mi$.
\end{theorem}

\begin{proof}
Let $S=\kk\langle \omega_1,\omega_2,y_1,y_2\rangle$. By \eqref{eq:Fnormal}, $\omega_1\omega_2=\omega_2\omega_1$ and $y_2y_1=py_1y_2$, $S$ is a quotient of the quantum affine space with these relations, and it is that quantum affine space because the monomials $\omega_1^a\omega_2^by_1^cy_2^d$ are independent (Lemma~\ref{lem:rank4} and the basis of $B$). If $\mathbb F$ is PI, Lemma~\ref{lem:qplane} applied to $\omega_1$ and $y_1$ gives $\lambda\in\mi$, that is, $2h^2\in\mi$ because $p\in\mi$. Conversely, Lemma~\ref{lem:rank4} and the basis give
\begin{equation}\label{eq:Fdec}
\mathbb F(p,h)=S\oplus x_1S\oplus x_2S\oplus x_1^2S ,
\end{equation}
a free right $S$-module of rank four. If $\lambda\in\mi$ then $S$ is PI, and Lemma~\ref{lem:finite} applies.
\end{proof}

\subsection{The family \texorpdfstring{$\mathbb I$}{I}}

\begin{theorem}\label{thm:I}
$\mathbb I(q,h)$ is PI if and only if $2h^2\in\mi$.
\end{theorem}

\begin{proof}
Let $\omega_1=y_1y_2$, $\omega_2=y_1^2+y_2^2$ and $\lambda'=-2qh^2$. We claim that
\begin{equation}\label{eq:Inormal}
\omega_kx_i=\lambda'\,x_i\omega_k\qquad(i,k\in\{1,2\}).
\end{equation}
For $h=1$ this is in the proof of \cite[Proposition 5.13]{ZZ09}, with the scalar $-(1+q)^2=-2q$, and \eqref{eq:twistscalar} then gives it for every $h$. We also derive it from Lemma~\ref{lem:Inorm}, which we need later anyway. Take $w_1,w_2$, $Z_1,Z_2$ and $c=h(1-q)$ as in Lemma~\ref{lem:Inorm} below, whose proof does not use this section; note that $c^2=h^2(1-2q+q^2)=\lambda'$. By that lemma, conjugation by $Z_j^{-1}$ sends $y_1,y_2$ to $c\,y_1,-c\,y_2$ for $j=1$ and to $c\,y_2,c\,y_1$ for $j=2$, so it multiplies $\omega_1$ by $-c^2$ and $\omega_2$ by $c^2$. Since $Z_1-Z_2=(1-q)x_2w_2$, we have
\begin{align*}
x_2&=(1-q)^{-1}(Z_1-Z_2)\,w_2^{-1}, & x_1&=(Z_1-x_2w_2)\,w_1^{-1}=\bigl(Z_1-(1-q)^{-1}(Z_1-Z_2)\bigr)w_1^{-1},
\end{align*}
so each $x_i$ is a sum $Z_1d_1+Z_2d_2$ with $d_1,d_2\in\kk w_1^{-1}+\kk w_2^{-1}$. Now $\omega_2$ is central in $A_y$, while $\omega_1$ anticommutes with $y_1,y_2$, hence with $w_1^{-1},w_2^{-1}$. Hence $\omega_1Z_jd_j=-c^2Z_j\omega_1d_j=c^2Z_jd_j\omega_1$ and $\omega_2Z_jd_j=c^2Z_jd_j\omega_2$, which is \eqref{eq:Inormal}. Lemma~\ref{lem:qplane} applied to $x_1$ and $\omega_1$ shows that $\lambda'\in\mi$ if $\mathbb I$ is PI, that is, $2h^2\in\mi$. Now argue as for $\mathbb F$, with $S'=\kk\langle x_1,x_2,\omega_1,\omega_2\rangle$ and $\mathbb I=S'\oplus S'y_1\oplus S'y_2\oplus S'y_1^2$ as a left $S'$-module.
\end{proof}

\subsection{The family \texorpdfstring{$\mathbb C$}{C}}

\begin{theorem}\label{thm:C}
$\mathbb C(p,h)$ is PI if and only if $3h^2\in\mi$.
\end{theorem}

\begin{proof}
Let $a=y_1^2y_2$ and $b=y_1^3-y_2^3$. In the proof of \cite[Proposition 5.17]{ZZ09} it is shown that, for $h=1$, $ax=(p-1)^2(1-p^2)\,xa$ for $x\in\kk x_1+\kk x_2$, $bx_2=(1-p^2)^3x_2b$ and $bx_1=(p-1)(1-p^2)(p^2-p)\,x_1b$. With $p^2+p+1=0$ the first scalar is $3(1-p)$ and the last two are both $3(2p+1)$, and by \eqref{eq:twistscalar}, in $\mathbb C(p,h)$,
\begin{align}
ax_i&=\alpha\,x_ia, & bx_i&=\beta\,x_ib, & \alpha&=3(1-p)h^3, & \beta&=3(2p+1)h^3 .\label{eq:Cnormal}
\end{align}
Then $\alpha=p^2\beta$, $\alpha^2=-p(3h^2)^3$ and $\beta^2=-(3h^2)^3$, so $\alpha,\beta\in\mi$ if and only if $3h^2\in\mi$.

Since $p^3=1$, $Y_1=y_1^3$ and $Y_2=y_2^3$ are central in $A_y$, so $a$ and $b$ commute. We have $b=Y_1-Y_2$ and $a^3=Y_1^2Y_2$, so $Y_1$ is a root of $T^3-bT^2-a^3$, and $A_y$, which is free of rank nine over $\kk[Y_1,Y_2]$, is a finite left module over $\kk[a,b]$; in particular $a,b$ are algebraically independent. Hence $S=\kk\langle x_1,x_2,a,b\rangle$ is a quantum affine space with parameters $p,\alpha,\beta$, and $\mathbb C=\sum_wS\,w$ for a finite set of $w\in A_y$. If $\mathbb C$ is PI, Lemma~\ref{lem:qplane} for $x_1$ and $a$ gives $\alpha\in\mi$. If $3h^2\in\mi$, $S$ is PI and Lemma~\ref{lem:finite} applies.
\end{proof}

\section{The PI degree of the monomial families}\label{sec:monodeg}

Here $B$ is one of $\mathbb B$, $\mathbb E$, $\mathbb J$, $\mathbb M$, $\mathbb N$, $\mathbb O$, $\mathbb P$, $\mathbb S$, $\mathbb T$, $\mathbb U$, $\mathbb W$, $\mathbb Z$ (with $f\ne0$ where it applies), with the datum of Table~\ref{tab:mono}; the family $\mathbb J$ is reduced to $\mathbb E$ by \eqref{eq:EJ}. The same method gives back the PI degrees of $\mathbb K$ \cite{GV-K} and $\mathbb L$ (Theorem~\ref{thm:L}).

\subsection{The three kinds of plane}

Call a basis $v_1,v_2$ of a plane \emph{balanced} if the plane is $\kk_{-1}$ and $v_1^2=v_2^2$; every basis of the form $v_1=w_1+cw_2$, $v_2=w_1-cw_2$ in a coordinate basis $w_1,w_2$ of $\kk_{-1}$ is balanced. Inspecting Table~\ref{tab:mono}, each side of each datum is of one of three kinds:
\begin{center}
\begin{tabular}{|c|c|c|}
\hline
kind & families (side) & rule for $w=v_1^{-1}v_2$\\\hline
balanced & $x$: $\mathbb E,\mathbb M,\mathbb O,\mathbb P,\mathbb S,\mathbb T,\mathbb U,\mathbb Z$;\ \ $y$: $\mathbb M,\mathbb N,\mathbb O,\mathbb S,\mathbb T,\mathbb U,\mathbb W$ & $v_1wv_1^{-1}=w^{-1}$\\
commutative & $x$: $\mathbb W$;\ \ $y$: $\mathbb Z$ & $v_1wv_1^{-1}=w$\\
coordinate & $x$: $\mathbb B$ ($q=p$), $\mathbb N$ ($q=-1$);\ \ $y$: $\mathbb B$, $\mathbb E$ ($q=p$), $\mathbb P$ ($q=-1$) & $v_1wv_1^{-1}=q^{-1}w$\\\hline
\end{tabular}
\end{center}
In the balanced case $v_1^2$ is central in the plane and $v_2v_1^{-1}=w^{-1}$, because $v_1^{-1}=v_1/v_1^2$ and $v_2^{-1}=v_2/v_1^2$.

In $Q(B)$ put
\begin{align}
e&=e_1, & r&=e_1^{-1}e_2, & X&=e^2, & z&=z_1, & \rho&=z_1^{-1}z_2, & Y&=z^2 .\notag
\end{align}
For elements $a,b$ we write $a\triangleright b=aba^{-1}$.

\begin{proposition}\label{prop:crossed}
Let $T'$ be the subalgebra of $Q(B)$ generated by $r^{\pm1},\rho^{\pm1},X^{\pm1},Y^{\pm1}$. Then:
\begin{enumerate}
\item[(i)] $T'$ is a quantum torus, and $e$ and $z$ normalize it.
\item[(ii)] $Q(B)=Q(T')\oplus Q(T')e\oplus Q(T')z\oplus Q(T')ez$.
\item[(iii)] Assume that $B$ is PI. Let
$G\subseteq\operatorname{Aut}(K)$,
$K=Z(Q(T'))$, be the group generated by
$e\triangleright$ and $z\triangleright$,
and let $n'=\PIdeg T'$. If $G$ has order four, then $\PIdeg B=4n'$. If $G$ has order two, then $\PIdeg B=2n'$.
\end{enumerate}
\end{proposition}

\begin{proof}
(i) The relations of the next subsection (proved there), together with Lemmas~\ref{lem:A2} and~\ref{lem:A3}(iv) for $z\triangleright X$ and $e\triangleright Y$, show that $e$ and $z$ send each of $r,\rho,X,Y$ to a scalar multiple of a monomial in them, and that these four elements commute up to scalars. On the $x$-side, $X^cr^a$ is a nonzero multiple of $X^{c-a}(e_1e_2)^a$ in the balanced and commutative cases and of $x_1^{2c-a}x_2^a$ in the coordinate case. These elements of $Q(A_x)$ are linearly independent. In the commutative and coordinate cases they are distinct monomials in $e_1^{\pm1},e_2^{\pm1}$, respectively $x_1^{\pm1},x_2^{\pm1}$. In the balanced case $X$ is central in $A_x$ and $e_2e_1=X^2(e_1e_2)^{-1}$; by Lemma~\ref{lem:plane}, with $u=v=X$, the elements $X$ and $e_1e_2+e_2e_1$ are algebraically independent, hence so are the commuting elements $X$ and $e_1e_2$. The same holds on the $y$-side, and Lemma~\ref{lem:fracind} shows that the elements $r^aX^c\rho^bY^d$ are linearly independent. Up to scalars these are the monomials of $T'$.

(ii) The sum $R$ of the four spaces is closed under multiplication, because $e^2=X$, $z^2=Y$, and $ze$ is a multiple of $e\,z$ times a monomial of $T'$ by \eqref{eq:mono}. It contains $e_1$, $e_2=er$, $z_1$ and $z_2=z\rho$, so $B\subseteq R$; it is a domain of dimension at most four over $Q(T')$, hence a division ring, and $R=Q(B)$. The sum is direct: the automorphisms of $B$ given by $x_i\mapsto-x_i$ and by $y_j\mapsto-y_j$ extend to $Q(B)$ and fix $r,\rho,X,Y$, hence $Q(T')$, while they multiply $e$, $z$, $ez$ by the pairs of signs $(-1,1)$, $(1,-1)$, $(-1,-1)$. So the four summands lie in different joint eigenspaces.

(iii) This is the usual computation for crossed products. If $\sum_ga_gw_g$, with $w_g\in\{1,e,z,ez\}$, commutes with $Q(T')$, then $a_g\,g(c)=c\,a_g$ for $c\in K$, so $a_g=0$ unless $g$ acts trivially on $K$. If $G$ has order four, the centralizer of $Q(T')$ is $K$, the center of $Q(B)$ is $K^G$, and $[Q(B):K^G]=4\cdot n'^2\cdot4$. If $G$ has order two, let $w$ be the element of $\{e,z,ez\}$ acting trivially on $K$, and $w'$ one of the other two. By Skolem--Noether $w\triangleright$ agrees on $Q(T')$ with conjugation by some $t\in Q(T')$, so $u_0=t^{-1}w$ centralizes $Q(T')$. Since $e^2=X$, $z^2=Y$ and $(ez)^2\in Q(T')$, we have $w^2\in Q(T')$, so $u_0^2\in Q(T')$ centralizes $Q(T')$ and lies in $K$; and $u_0\notin Q(T')$ by (ii). Hence $K'=K(u_0)$ is a field of degree two over $K$. It is the centralizer of $Q(T')$: by the computation above, such an element lies in $Q(T')+Q(T')u_0$, and then in $K+Ku_0$. So $Q(T')K'\cong Q(T')\otimes_KK'$ is central simple of degree $n'$ over $K'$, and $Q(B)=Q(T')K'\oplus Q(T')K'w'$. Conjugation by $w'$ normalizes $Q(T')$, so it preserves $K'$, and it is nontrivial on $K\subseteq K'$; its square is conjugation by an element of $Q(T')K'$. Thus $Z(Q(B))$ is the fixed field of $w'\triangleright$ in $K'$, and $[Q(B):Z(Q(B))]=2\cdot n'^2\cdot2$.
\end{proof}

\subsection{The relations}

All formulas below follow from \eqref{eq:mono} and the rules of the table above; the proofs come after the lists. We write $\kappa$ for the scalar with $\rho\triangleright r=\kappa r$; $X$ commutes with $r$ except in $\mathbb B$, and $Y$ with $\rho$ except in $\mathbb B$ and $\mathbb E$.

\emph{Type I, both sides balanced} ($\mathbb M$, $\mathbb O$):
\begin{align*}
e\triangleright r&=r^{-1}, & e\triangleright\rho&=\tfrac{\lambda_{21}}{\lambda_{11}}\rho^{-1}, & z\triangleright r&=\tfrac{\lambda_{22}}{\lambda_{21}}r^{-1}, & z\triangleright\rho&=\rho^{-1},\\
\kappa&=\tfrac{\lambda_{12}\lambda_{21}}{\lambda_{11}\lambda_{22}}, & z\triangleright X&=\lambda_{11}\lambda_{21}X, & e\triangleright Y&=\tfrac{1}{\lambda_{21}\lambda_{22}}Y, & Y\triangleright X&=(\lambda_{11}\lambda_{21})^2X,
\end{align*}
and $X\rho=\rho X$, $Yr=rY$.

\emph{Type I, both sides coordinate with $q=p$} ($\mathbb B$): $e\triangleright r=p^{-1}r$, $z\triangleright\rho=p^{-1}\rho$, $z\triangleright r=-r^{-1}$, $e\triangleright\rho=-\rho^{-1}$, $\kappa=-1$, $z\triangleright X=-ph^2Xr^2$, $e\triangleright Y=-ph^{-2}Y\rho^2$, $X\triangleright r=-r$, $Y\triangleright\rho=-\rho$, $Y\triangleright X=-h^4X$, and $X\rho=\rho X$, $Yr=rY$.

\emph{Type I, $x$ balanced, $y$ coordinate} ($\mathbb E$): the same as for $\mathbb M$, $\mathbb O$, except $z\triangleright\rho=p^{-1}\rho$, $Y\triangleright\rho=p^{-2}\rho=-\rho$ and $e\triangleright Y=\frac{p}{\lambda_{21}\lambda_{22}}Y\rho^2$.

\emph{Type $0$, both balanced} ($\mathbb S$): $e\triangleright r=r^{-1}$, $z\triangleright\rho=\rho^{-1}$, $z\triangleright r=\frac{\lambda_{12}}{\lambda_{11}}r$, $e\triangleright\rho=\frac{\lambda_{11}}{\lambda_{21}}\rho$, $\kappa=\frac{\lambda_{11}\lambda_{22}}{\lambda_{12}\lambda_{21}}$, $X\triangleright\rho=(\frac{\lambda_{11}}{\lambda_{21}})^2\rho$, $Y\triangleright r=(\frac{\lambda_{12}}{\lambda_{11}})^2r$, $Y\triangleright X=\lambda_{11}^4X$.

\emph{Type II} ($\mathbb P$, $\mathbb T$, $\mathbb Z$; $x$ balanced): $e\triangleright r=r^{-1}$, $z\triangleright r=\frac{\lambda_{22}}{\lambda_{21}}r$, $e\triangleright\rho=\frac{\lambda_{21}}{\lambda_{11}}\rho^{-1}$, $\kappa=\frac{\lambda_{12}\lambda_{21}}{\lambda_{11}\lambda_{22}}$, $z\triangleright X=\lambda_{11}\lambda_{21}X$, $Y\triangleright r=(\frac{\lambda_{22}}{\lambda_{21}})^2r$, $Y\triangleright X=(\lambda_{11}\lambda_{21})^2X$, $X\rho=\rho X$; and $z\triangleright\rho$ is $\rho^{-1}$, $-\rho$ or $\rho$ for $\mathbb T$, $\mathbb P$, $\mathbb Z$.

\emph{Type III} ($\mathbb N$, $\mathbb U$, $\mathbb W$; $y$ balanced): $z\triangleright\rho=\rho^{-1}$, $z\triangleright r=\frac{\lambda_{12}}{\lambda_{11}}r^{-1}$, $\kappa=\frac{\lambda_{11}\lambda_{22}}{\lambda_{12}\lambda_{21}}$, $z\triangleright X$ equal to $\lambda_{11}^2X$, $-\lambda_{11}^2Xr^2$, $\lambda_{11}^2Xr^2$ for $\mathbb U$, $\mathbb N$, $\mathbb W$, $Y\triangleright X=(\lambda_{11}\lambda_{12})^2X$ for $\mathbb N$, $\mathbb W$ and $\lambda_{11}^4X$ for $\mathbb U$, $e\triangleright Y=(\lambda_{11}\lambda_{12})^{-1}Y$, $e\triangleright\rho=\frac{\lambda_{11}}{\lambda_{21}}\rho$; and
\[
e\triangleright r=\begin{cases}r^{-1}&(\mathbb U)\\ -r&(\mathbb N)\\ r&(\mathbb W)\end{cases}.
\]

We now prove these lists, from \eqref{eq:mono} and the rules for the three kinds of plane. Since $\sigma$ and $\tau$ are involutions, \eqref{eq:mono} can be solved for the products of $z_a^{\pm1}$ and $e_b^{\pm1}$ in any order.

\begin{lemma}\label{lem:A1}
For all $a,b,c,d\in\{1,2\}$,
\begin{align}
z_ae_c^{-1}&=\lambda_{\tau(a),c}^{-1}\,e_{\sigma(c)}^{-1}z_{\tau(a)}, & e_bz_c^{-1}&=\lambda_{\tau(c),b}\,z_{\tau(c)}^{-1}e_{\sigma(b)},\label{eq:A1}
\end{align}
and consequently
\begin{align}
z_a\triangleright(e_b^{-1}e_c)&=\frac{\lambda_{\tau(a),c}}{\lambda_{\tau(a),b}}\,e_{\sigma(b)}^{-1}e_{\sigma(c)}, & e_b\triangleright(z_c^{-1}z_d)&=\frac{\lambda_{\tau(c),b}}{\lambda_{\tau(d),b}}\,z_{\tau(c)}^{-1}z_{\tau(d)},\label{eq:A2}\\
z_a\triangleright e_b^2&=\lambda_{ab}\lambda_{\tau(a),b}\,e_{\sigma(b)}^2, & e_b\triangleright z_c^2&=\bigl(\lambda_{\tau(c),b}\lambda_{\tau(c),\sigma(b)}\bigr)^{-1}z_{\tau(c)}^2 .\label{eq:A3}
\end{align}
\end{lemma}

\begin{proof}
By \eqref{eq:mono}, $z_{\tau(a)}e_c=\lambda_{\tau(a),c}\,e_{\sigma(c)}z_a$. Multiplying by $e_{\sigma(c)}^{-1}$ on the left and by $e_c^{-1}$ on the right gives the first formula of \eqref{eq:A1}, and inverting both sides of it, with $(a,c)$ renamed $(c,b)$, gives the second.

For \eqref{eq:A2}, the first formula of \eqref{eq:A1} and \eqref{eq:mono} give
\[
z_ae_b^{-1}e_cz_a^{-1}=\lambda_{\tau(a),b}^{-1}\,e_{\sigma(b)}^{-1}\bigl(z_{\tau(a)}e_c\bigr)z_a^{-1}=\frac{\lambda_{\tau(a),c}}{\lambda_{\tau(a),b}}\,e_{\sigma(b)}^{-1}e_{\sigma(c)},
\]
and both formulas of \eqref{eq:A1} give
\[
e_bz_c^{-1}z_de_b^{-1}=\lambda_{\tau(c),b}\,z_{\tau(c)}^{-1}e_{\sigma(b)}\bigl(z_de_b^{-1}\bigr)=\frac{\lambda_{\tau(c),b}}{\lambda_{\tau(d),b}}\,z_{\tau(c)}^{-1}e_{\sigma(b)}e_{\sigma(b)}^{-1}z_{\tau(d)} .
\]
For \eqref{eq:A3}, $z_ae_be_b=\lambda_{ab}\,e_{\sigma(b)}z_{\tau(a)}e_b=\lambda_{ab}\lambda_{\tau(a),b}\,e_{\sigma(b)}^2z_a$. Applying the second formula of \eqref{eq:A1} twice,
\[
e_bz_c^{-2}=\lambda_{\tau(c),b}\,z_{\tau(c)}^{-1}e_{\sigma(b)}z_c^{-1}=\lambda_{\tau(c),b}\lambda_{\tau(c),\sigma(b)}\,z_{\tau(c)}^{-2}e_b,
\]
and inverting gives the last formula.
\end{proof}

\begin{lemma}\label{lem:A2}
Let $v_1,v_2$ be one of the two bases of a datum, $w=v_1^{-1}v_2$ and $V=v_1^2$. Then $v_2=v_1w$ and $v_2^{-1}v_1=w^{-1}$. Moreover:
\begin{enumerate}
\item[(i)] if the basis is balanced, $v_1\triangleright w=w^{-1}$, $v_2^2=V$, and $V$ commutes with $w$;
\item[(ii)] if the plane is commutative, $v_1\triangleright w=w$, $v_2^2=Vw^2$, and $V$ commutes with $w$;
\item[(iii)] if the basis is a coordinate basis with $v_2v_1=q\,v_1v_2$, then $v_1\triangleright w=q^{-1}w$, $V\triangleright w=q^{-2}w$ and $v_2^2=q\,Vw^2$.
\end{enumerate}
\end{lemma}

\begin{proof}
In (i), $v_1^2=v_2^2$ is central in the plane, as is every square of a linear form in $\kk_{-1}[x_1,x_2]$, and
\[
v_1wv_1^{-1}=v_2v_1^{-1}=v_2^{-1}v_2^2v_1^{-1}=v_2^{-1}v_1=w^{-1}.
\]
Part (ii) is clear. In (iii), $v_1^{-1}v_2v_1=q\,v_2$ gives $v_2v_1^{-1}=q^{-1}v_1^{-1}v_2$, that is, $v_1\triangleright w=q^{-1}w$. Applying this twice, $V\triangleright w=q^{-2}w$, so $wV=q^2Vw$ and
\[
v_2^2=v_1wv_1w=(v_1\triangleright w)\,Vw=q^{-1}wVw=q\,Vw^2. \qedhere
\]
\end{proof}

Conjugation is multiplicative, so $Y\triangleright{}=(z\triangleright{})^2$, $X\triangleright{}=(e\triangleright{})^2$ and $\rho\triangleright r=z_1^{-1}\triangleright(z_2\triangleright r)$. Put
\[
\mu_a=\frac{\lambda_{\tau(a),2}}{\lambda_{\tau(a),1}}\qquad(a=1,2).
\]

\begin{lemma}\label{lem:A3}
For every monomial datum:
\begin{enumerate}
\item[(i)] $z_a\triangleright r=\mu_a\,r$ if $\sigma$ is the identity, and $z_a\triangleright r=\mu_a\,r^{-1}$ otherwise;
\item[(ii)] $\kappa=\mu_2/\mu_1$, where $\rho\triangleright r=\kappa r$;
\item[(iii)] $e\triangleright\rho=\dfrac{\lambda_{\tau(1),1}}{\lambda_{\tau(2),1}}\,\rho^{\pm1}$, with the sign $+$ exactly when $\tau$ is the identity;
\item[(iv)] $z\triangleright X=\lambda_{11}\lambda_{\tau(1),1}\,e_{\sigma(1)}^2$ and $e\triangleright Y=\bigl(\lambda_{\tau(1),1}\lambda_{\tau(1),\sigma(1)}\bigr)^{-1}z_{\tau(1)}^2$.
\end{enumerate}
\end{lemma}

\begin{proof}
Part (i) is \eqref{eq:A2} with $b=1$, $c=2$, together with $e_2^{-1}e_1=r^{-1}$. For (ii), if $\sigma$ is the identity then $z_1^{-1}\triangleright r=\mu_1^{-1}r$ and $\rho\triangleright r=z_1^{-1}\triangleright(\mu_2r)=\mu_2\mu_1^{-1}r$. Otherwise $z_1\triangleright r^{-1}=\mu_1^{-1}r$, so $z_1^{-1}\triangleright r=\mu_1r^{-1}$ and
\[
\rho\triangleright r=z_1^{-1}\triangleright(\mu_2r^{-1})=\mu_2\,(\mu_1r^{-1})^{-1}=\mu_2\mu_1^{-1}r .
\]
Part (iii) is \eqref{eq:A2} with $b=c=1$, $d=2$, since $z_{\tau(1)}^{-1}z_{\tau(2)}$ is $\rho$ or $z_2^{-1}z_1=\rho^{-1}$. Part (iv) is \eqref{eq:A3} with $a=b=1$ and $b=c=1$.
\end{proof}

\begin{proof}[Proof of the lists above]
When $z\triangleright r$ is a multiple of $r^{-1}$, say $\mu r^{-1}$, then $Y\triangleright r=z\triangleright(\mu r^{-1})=\mu(\mu r^{-1})^{-1}=r$; in the same way $X\triangleright\rho=\rho$ when $e\triangleright\rho$ is a multiple of $\rho^{-1}$. The relations $e\triangleright r$ and $z\triangleright\rho$ are Lemma~\ref{lem:A2} applied to the two sides, and so is the action of $X$ on $r$ (in $\mathbb N$, Lemma~\ref{lem:A2}(iii) gives $X\triangleright r=q^{-2}r=r$ because $q=-1$).

\emph{Type I} ($\sigma$ and $\tau$ both the transposition). Here $\mu_1=\lambda_{22}/\lambda_{21}$ and $\mu_2=\lambda_{12}/\lambda_{11}$, and Lemma~\ref{lem:A3} gives
\begin{align*}
z\triangleright r&=\frac{\lambda_{22}}{\lambda_{21}}\,r^{-1}, & \kappa&=\frac{\lambda_{12}\lambda_{21}}{\lambda_{11}\lambda_{22}}, & e\triangleright\rho&=\frac{\lambda_{21}}{\lambda_{11}}\,\rho^{-1}, & z\triangleright X&=\lambda_{11}\lambda_{21}\,e_2^2, & e\triangleright Y&=\frac{z_2^2}{\lambda_{21}\lambda_{22}} .
\end{align*}
In $\mathbb M$, $\mathbb O$ and $\mathbb E$ the $x$-side is balanced, so $e_2^2=X$ and $Y\triangleright X=(\lambda_{11}\lambda_{21})^2X$. In $\mathbb M$ and $\mathbb O$ the $y$-side is balanced as well, so $z_2^2=Y$. In $\mathbb E$ it is a coordinate basis with $q=p$, so $z\triangleright\rho=p^{-1}\rho$, $Y\triangleright\rho=p^{-2}\rho=-\rho$ and $z_2^2=pY\rho^2$.

In $\mathbb B$ both bases are coordinate bases with $q=p$ and $(\lambda_{ab})=h(1,1,-1,1)$, so $\mu_1=-1$, $\mu_2=1$. Lemma~\ref{lem:A2}(iii) gives $e\triangleright r=p^{-1}r$, $z\triangleright\rho=p^{-1}\rho$, $X\triangleright r=-r$, $Y\triangleright\rho=-\rho$, $e_2^2=pXr^2$ and $z_2^2=pY\rho^2$. Hence $z\triangleright X=-ph^2Xr^2$, $e\triangleright Y=-ph^{-2}Y\rho^2$ and
\[
Y\triangleright X=-ph^2\,(z\triangleright X)(z\triangleright r)^2=p^2h^4\,Xr^2r^{-2}=-h^4X .
\]

\emph{Type $0$} ($\mathbb S$). Here $\mu_1=\lambda_{12}/\lambda_{11}$ and $\mu_2=\lambda_{22}/\lambda_{21}$, so
\begin{align*}
z\triangleright r&=\frac{\lambda_{12}}{\lambda_{11}}\,r, & \kappa&=\frac{\lambda_{11}\lambda_{22}}{\lambda_{12}\lambda_{21}}, & e\triangleright\rho&=\frac{\lambda_{11}}{\lambda_{21}}\,\rho, & z\triangleright X&=\lambda_{11}^2X,
\end{align*}
and squaring, $Y\triangleright r=(\lambda_{12}/\lambda_{11})^2r$, $X\triangleright\rho=(\lambda_{11}/\lambda_{21})^2\rho$ and $Y\triangleright X=\lambda_{11}^4X$.

\emph{Type II} ($\tau$ the transposition, $\sigma$ the identity). Here $\mu_1=\lambda_{22}/\lambda_{21}$ and $\mu_2=\lambda_{12}/\lambda_{11}$, so
\begin{align*}
z\triangleright r&=\frac{\lambda_{22}}{\lambda_{21}}\,r, & \kappa&=\frac{\lambda_{12}\lambda_{21}}{\lambda_{11}\lambda_{22}}, & e\triangleright\rho&=\frac{\lambda_{21}}{\lambda_{11}}\,\rho^{-1}, & z\triangleright X&=\lambda_{11}\lambda_{21}X,
\end{align*}
and squaring, $Y\triangleright r=(\lambda_{22}/\lambda_{21})^2r$ and $Y\triangleright X=(\lambda_{11}\lambda_{21})^2X$. The $y$-side is balanced in $\mathbb T$, a coordinate basis with $q=-1$ in $\mathbb P$ and commutative in $\mathbb Z$, which gives $z\triangleright\rho=\rho^{-1}$, $-\rho$, $\rho$.

\emph{Type III} ($\sigma$ the transposition, $\tau$ the identity). Here $\mu_1=\lambda_{12}/\lambda_{11}$ and $\mu_2=\lambda_{22}/\lambda_{21}$, so
\begin{align*}
z\triangleright r&=\frac{\lambda_{12}}{\lambda_{11}}\,r^{-1}, & \kappa&=\frac{\lambda_{11}\lambda_{22}}{\lambda_{12}\lambda_{21}}, & e\triangleright\rho&=\frac{\lambda_{11}}{\lambda_{21}}\,\rho, & z\triangleright X&=\lambda_{11}^2\,e_2^2, & e\triangleright Y&=\frac{Y}{\lambda_{11}\lambda_{12}} .
\end{align*}
The $y$-side is balanced in all three families. On the $x$-side, Lemma~\ref{lem:A2} gives $e_2^2=X$ in $\mathbb U$, $e_2^2=-Xr^2$ in $\mathbb N$ (coordinate, $q=-1$) and $e_2^2=Xr^2$ in $\mathbb W$ (commutative). This gives the three values of $z\triangleright X$. In $\mathbb U$, $Y\triangleright X=\lambda_{11}^4X$. In $\mathbb N$ and $\mathbb W$, with $z\triangleright X=\pm\lambda_{11}^2Xr^2$,
\[
Y\triangleright X=\pm\lambda_{11}^2\,(z\triangleright X)(z\triangleright r)^2=\lambda_{11}^4\,Xr^2\Bigl(\frac{\lambda_{12}}{\lambda_{11}}\Bigr)^2r^{-2}=(\lambda_{11}\lambda_{12})^2X. \qedhere
\]
\end{proof}

\subsection{The formulas}

\begin{theorem}\label{thm:monodeg}
Let $\varepsilon(u,v)=2$ if $u$ and $v$ both have odd order, and $\varepsilon(u,v)=4$ otherwise. For the PI members of the families below,
\[
\begin{array}{|c|c|c|}
\hline
\text{family} & \PIdeg & u,\ v\\\hline
\mathbb B & 8\ord(h^4) & \\
\mathbb S,\ \mathbb T,\ \mathbb U & 8\ord(16h^4) & \\
\mathbb E,\ \mathbb J & 8\ord(4h^4) & \\
\mathbb M & 4\ord(u^2)\ord(v^2) & u=\frac{s+1}{s-1},\ v=h^2(f-1)\\
\mathbb N & 4\ord(u^2)\ord(v^2) & u=\frac{f-g}{f+g},\ v=f^2-g^2\\
\mathbb O & \varepsilon(u,v)\ord(u^2)\ord(v^2) & u=\frac{1+s}{1-s},\ v=h^2(1-f)\\
\mathbb P & 4\ord(u^2)\ord(v^2) & u=\frac{1+s}{1-s},\ v=h^2(1-f)\\
\mathbb W & \varepsilon(u,v)\ord(u^2)\ord(v^2) & u=\frac{s-i}{s+i},\ v=h^2(f+1)\\
\mathbb Z & \varepsilon(u,v)\ord(u^2)\ord(v^2) & u=\frac{i-s}{i+s},\ v=h^2(f+1)\\\hline
\end{array}
\]
\end{theorem}

\begin{proof}
In each case we first split $T'$ as a tensor product of two quantum tori of rank two, then decide whether $G$ has order four or two, and apply Proposition~\ref{prop:crossed}.

$\mathbb O$: here $(\lambda_{ab})=(hb',ha',ha',hb')$ with $a'=1+s$, $b'=1-s$, so $\kappa=u^2$ and $Y\triangleright X=v^2X$; $T'$ is the tensor product of the tori on $r,\rho$ and on $X,Y$, and $n'=\ord(u^2)\ord(v^2)$. The center $K$ is generated by $r^{\pm m},\rho^{\pm m},X^{\pm m'},Y^{\pm m'}$, $m=\ord(u^2)$, $m'=\ord(v^2)$. Both $e\triangleright$ and $z\triangleright$ invert $r^m$ or $\rho^m$, so they are not trivial on $K$, while $ez$ acts by
\[
r\mapsto u^{-1}r,\qquad \rho\mapsto u^{-1}\rho,\qquad X\mapsto vX,\qquad Y\mapsto v^{-1}Y ,
\]
which is trivial on $K$ exactly when $u^m=v^{m'}=1$, that is, when $u$ and $v$ have odd order.

$\mathbb M$: now $(\lambda_{ab})=h(i(s+1),-i(s-1),i(s-1),-i(s+1))$, $\kappa=u^{-2}$ and $(\lambda_{11}\lambda_{21})^2=v^2$, so $n'$ is as before. But $ez$ sends $r\mapsto-ur$ and $\rho\mapsto u\rho$; triviality on $r^m$ and $\rho^m$ would need $(-1)^m=1$ and $u^m=1$, and $u^m=1$ with $m=\ord(u^2)$ forces $m$ odd. So $G$ has order four.

$\mathbb S$, $\mathbb T$, $\mathbb U$: $(\lambda_{ab})$ is $2h(1,-1,-1,-1)$, $h(2i,-2,-2i,-2)$, $h(2i,-2i,-2,-2)$. In all three $\kappa=-1$, and the remaining scalars are $Y\triangleright X=16h^4X$ together with $Y\triangleright r=-r$ in $\mathbb T$ and $X\triangleright\rho=-\rho$ in $\mathbb U$. Replacing $Y$ by $Y\rho$ in $\mathbb T$ and $X$ by $Xr$ in $\mathbb U$ gives the tensor product of the tori on $r,\rho$ (parameter $-1$) and on two generators with parameter $16h^4$, so $n'=2\ord(16h^4)$. Each nontrivial element of $G$ inverts $r$ or $\rho$, so $G$ has order four.

$\mathbb B$: here $(\lambda_{ab})=h(1,1,-1,1)$ and both bases are coordinate bases of $\kk_p$. By the list above, $\kappa=-1$, $X\triangleright r=-r$, $Y\triangleright\rho=-\rho$ and $Y\triangleright X=-h^4X$, while $X$ commutes with $\rho$ and $Y$ with $r$. With $X'=X\rho$ and $Y'=Yr$ the torus splits into the tori on $r,\rho$ (parameter $-1$) and on $X',Y'$, where $Y'\triangleright X'=h^4X'$; so $n'=2\ord(h^4)$. The element $e$ inverts $\rho$, and $z$ and $ez$ invert $r$, so $G$ has order four.

$\mathbb E$: $\kappa=-1$, $Y\triangleright\rho=-\rho$, $Y\triangleright X=(\lambda_{11}\lambda_{21})^2X=4h^4X$. With $Y'=Yr$ the torus splits and $n'=2\ord(4h^4)$. $e$ and $z$ invert $r$, and $ez$ inverts $\rho$. For $\mathbb J$ use \eqref{eq:EJ}: the opposite ring has the same PI degree, and the formula for $\mathbb E$ does not depend on $p$.

$\mathbb P$, $\mathbb Z$ (type II): $(\lambda_{ab})$ is $h(a',b',b',a')$ for $\mathbb P$ and $(A,B,B,A)$ with $A=ih(s-i)$, $B=-ih(s+i)$ for $\mathbb Z$, so $\kappa=u^{-2}$, $Y\triangleright r=u^2r$, $Y\triangleright X=v^2X$, and $Y'=Y\rho$ commutes with $r$ and $\rho$; hence $n'=\ord(u^2)\ord(v^2)$. The elements $e$ and $ez$ invert $r$. The element $z$ acts by $r\mapsto ur$, $X\mapsto vX$, and by $\rho\mapsto-\rho$ in $\mathbb P$, $\rho\mapsto\rho$ in $\mathbb Z$. In $\mathbb Z$ it is trivial on $K$ exactly when $u,v$ have odd order. In $\mathbb P$ it would need $u^m=1$ and $(-1)^m=1$, which is impossible as for $\mathbb M$.

$\mathbb N$, $\mathbb W$ (type III): $(\lambda_{ab})$ is $(a,b,-b,-a)$ with $a=f-g$, $b=f+g$ for $\mathbb N$, and $(A,B,-B,-A)$ for $\mathbb W$; so $\kappa=(\lambda_{11}/\lambda_{12})^2$, $X\triangleright\rho=\kappa\rho$ and $Y\triangleright X=(\lambda_{11}\lambda_{12})^2X$; hence $X'=Xr$ commutes with $\rho$. Hence $n'=\ord(u^2)\ord(v^2)$, where for $\mathbb W$ we note $(A/B)^2=u^2$ and $AB=v$. The elements $z$ and $ez$ invert $\rho$. The element $e$ acts by $\rho\mapsto(\lambda_{11}/\lambda_{21})\rho$, $Y\mapsto v^{-1}Y$, and $r\mapsto-r$, $X'\mapsto-X'$ in $\mathbb N$, $r\mapsto r$, $X'\mapsto X'$ in $\mathbb W$. In $\mathbb W$, $\lambda_{11}/\lambda_{21}=-A/B=u$ and $e$ is trivial on $K$ exactly when $u$ and $v$ have odd order. In $\mathbb N$, triviality needs $(-1)^m=1$ and $(-u)^m=1$ with $m=\ord(u^2)$ even, so $\ord u=2m$ and $u^m=-1$, a contradiction.
\end{proof}

\begin{remark}\label{rem:monodeg}
$\mathbb O$ and $\mathbb P$ share the criterion and the scalars $u,v$, yet for $u,v$ of odd order the PI degree of $\mathbb O$ is half that of $\mathbb P$. The algebras $\mathbb W$ and $\mathbb Z$ are $\Sigma$-$M$-dual \cite[Subcase 4.4.4]{ZZ09}, but the duality moves the parameters. With $s^2=f$, the assignment
\[
x_1\mapsto sy_1+y_2,\qquad x_2\mapsto s(sy_1-y_2),\qquad y_1\mapsto x_1,\qquad y_2\mapsto x_2
\]
defines an isomorphism $\mathbb W(f^{-1},sh)\cong\mathbb Z(f,h)^{\mathrm{op}}$, as one checks on the four mixed relations. Taking $s^{-1}$ as the square root of $f^{-1}$, the scalars $u,v$ of $\mathbb W(f^{-1},sh)$ are $u^{-1},v$ for $\mathbb Z(f,h)$, so the two formulas of Theorem~\ref{thm:monodeg} agree. At equal values of the parameters they differ: for $f=3$ and $h=-\tfrac12$, $\PIdeg\mathbb Z=6$ and $\PIdeg\mathbb W=12$, while the dual member $\mathbb W(\tfrac13,-\tfrac{\sqrt3}2)$ has PI degree $6$. When $\alpha$ and $\beta$ have odd order, the formula for $\mathbb O$ agrees with \cite[Theorem~4.3]{Rub-ozone}: there $\alpha=h(1+s)$ and $\beta=h(1-s)$, so $u=\alpha/\beta$, $v=\alpha\beta$, and $\PIdeg\mathbb O=2\ord(u)\ord(v)$.
\end{remark}

\section{The families \texorpdfstring{$\mathbb F$}{F}, \texorpdfstring{$\mathbb I$}{I}, \texorpdfstring{$\mathbb Q$}{Q}, \texorpdfstring{$\mathbb R$}{R} and \texorpdfstring{$\mathbb V$}{V}}\label{sec:norm}

For these five families we use elements of bidegree $(1,1)$
or $(1,2)$ normalizing one of the planes, and $Q(B)$ is again
a crossed product over a quantum torus. Here a \emph{quantum torus} in a division ring $D$ is a subalgebra $T'$ generated by units $v_1^{\pm1},\dots,v_4^{\pm1}$ that commute up to scalars and whose monomials are linearly independent; $Q(T')\subseteq D$ is its quotient division ring, $K=Z(Q(T'))$, and $a\triangleright b=aba^{-1}$.

\subsection{Two lemmas}

\begin{lemma}\label{lem:cross}
Let $D$ be a division ring containing $\kk$, $T'\subseteq D$ a quantum torus and $w_1=1,w_2,\dots,w_N\in D^\times$ such that
\begin{enumerate}
\item[(i)] each $w_g\triangleright$ restricts to an automorphism of $T'$ that sends every generator to a scalar multiple of a monomial of $T'$;
\item[(ii)] for all $g,g'$ there is $g''$ with $w_gw_{g'}\in Q(T')\,w_{g''}$;
\item[(iii)] the automorphisms of $K$ induced by $w_1\triangleright,\dots,w_N\triangleright$ are pairwise distinct;
\item[(iv)] $D$ is generated as a division ring by $T'$ and the $w_g$.
\end{enumerate}
Then $D=\bigoplus_gQ(T')w_g$, these automorphisms form a group $G$ of order $N$, and $Z(D)=K^G$. If $T'$ is PI of PI degree $n'$, then $D$ has dimension $(Nn')^2$ over its center; if $T'$ is not PI, neither is any subalgebra of $D$ with quotient division ring $D$.
\end{lemma}

\begin{proof}
By (i), conjugation by each $w_g$ is an automorphism of $Q(T')$, so it induces an automorphism of $K$ and (iii) makes sense. By (ii) the sum $R=\sum_gQ(T')w_g$ is a subring, and the $N$ automorphisms of (iii) are closed under composition; being distinct, they form a group. Suppose $\sum_ga_gw_g=0$ with $a_g\in Q(T')$, not all zero and with the fewest nonzero terms, and let $a_{g_0}\ne0$. For $c\in K$, multiplying by $c$ on the right and by $g_0(c)$ on the left and subtracting gives $\sum_ga_g\bigl(g(c)-g_0(c)\bigr)w_g=0$, a shorter relation. So $g(c)=g_0(c)$ for all $c$ whenever $a_g\ne0$, and by (iii) only $a_{g_0}$ is nonzero, which is absurd. Hence the sum is direct, $R$ is a domain of finite dimension over $Q(T')$, so a division ring, and $R=D$ by (iv).

If $\sum_ga_gw_g$ commutes with $K$, then $a_g(g(c)-c)=0$ for all $c\in K$, so $a_g=0$ for $g\ne1$. Thus a central element of $D$ lies in $Q(T')$, hence in $K$, and it must be fixed by $G$; conversely $K^G$ is central. If $T'$ is PI, $[D:K^G]=N\cdot n'^2\cdot N$. The last claim is clear, since $Q(T')\subseteq D$.
\end{proof}

Condition (i) only has to be checked on the generators when $w_g$ is a product of elements $w$ with $w^2\in T'$. Indeed, if $w\triangleright$ sends the generators to scalar multiples of monomials, it maps $T'$ into itself, and its square is conjugation by the unit $w^2$ of $T'$, so it is onto. This is the situation in all our applications.

When the lattice $\Gamma$ of exponents of central monomials of $T'$ has finite index, $T'$ is free over its center $Z(T')=\kk[v^{\mathbf a}:\mathbf a\in\Gamma]$ (Lemma~\ref{lem:centerQT}) on monomials representing the cosets of $\Gamma$, so $n'^2=[\mathbb Z^4:\Gamma]$.

\begin{lemma}\label{lem:words}
Let $E=Q(A_x)\subseteq Q(B)$ and let $Z_1,Z_2\in Ey_1+Ey_2$ be left linearly independent over $E$ with $Z_iE=EZ_i$. Then for every $k\ge0$ the words of length $k$ in $Z_1,Z_2$ span a left $E$-space of dimension $k+1$. In particular, $k+1$ words of length $k$ that span it are left independent over $E$. The same holds with the two planes and left and right exchanged.
\end{lemma}

\begin{proof}
Since $B$ is a free left $A_x$-module on the $y_1^cy_2^d$ and a domain, $E(A_y)_k=\sum_{c+d=k}Ey_1^cy_2^d$ has left dimension $k+1$. Put $M=Ey_1+Ey_2=EZ_1+EZ_2$. Then $ME=M$, because $Z_iE=EZ_i$. So $y_jE\subseteq M$, and $M^k$ is at the same time $\sum E\,y_{j_1}\cdots y_{j_k}=E(A_y)_k$ and the $E$-span of the words of length $k$ in $Z_1,Z_2$.
\end{proof}

\subsection{The families \texorpdfstring{$\mathbb Q$}{Q} and \texorpdfstring{$\mathbb V$}{V}}

We keep the relations displayed before Proposition~\ref{prop:QV}, multiplied by $h$. In $A_y=\kk_{-1}[y_1,y_2]$ the squares $y_1^2,y_2^2$ are central, and we put
\begin{align}
\Theta_{\mathbb Q}&=x_2(y_2^2-y_1^2)-x_1(y_1^2-2y_1y_2+y_2^2), & \Theta_{\mathbb V}&=x_1(y_1^2-y_2^2)-x_2(y_1^2+y_1y_2).\label{eq:Theta}
\end{align}

\begin{lemma}\label{lem:Theta}
In $\mathbb Q(h)$ and $\mathbb V(h)$,
\begin{align*}
\mathbb Q:&\quad y_1\Theta=h\,\Theta y_1, & y_2\Theta&=h\,\Theta y_2, & \Theta x_1&=h^2x_1\Theta,\\
\mathbb V:&\quad y_1\Theta=h\,\Theta y_2, & y_2\Theta&=-h\,\Theta y_1, & \Theta x_2&=h^2x_2\Theta .
\end{align*}
\end{lemma}

\begin{proof}
Write $a=y_1^2$, $b=y_2^2$, $w=y_1y_2$, so that $y_1w=ay_2$, $y_2w=-by_1$, $wy_1=-ay_2$, $wy_2=by_1$. It is enough to take $h=1$ in the first two columns, because $\Theta$ has $x$-degree one. In $\mathbb Q$,
\begin{align*}
y_1\Theta&=(x_1y_1+x_2y_1+x_1y_2)(b-a)-x_1y_2(a+b-2w)=x_2y_1(b-a)-x_1\bigl(y_1(a+b)+2ay_2\bigr),\\
\Theta y_1&=x_2(b-a)y_1-x_1\bigl((a+b)y_1-2wy_1\bigr)=x_2y_1(b-a)-x_1\bigl((a+b)y_1+2ay_2\bigr),\\
y_2\Theta&=(x_1y_1-x_1y_2+x_2y_2)(b-a)+x_1y_1(a+b-2w)=x_2y_2(b-a)+x_1\bigl(2by_1-(a+b)y_2\bigr),\\
\Theta y_2&=x_2(b-a)y_2-x_1\bigl((a+b)y_2-2by_1\bigr).
\end{align*}
In $\mathbb V$,
\begin{align*}
y_1\Theta&=(x_2y_1+x_1y_2)(a-b)-x_2y_1(a+w)=x_1y_2(a-b)-x_2(by_1+ay_2)=\Theta y_2,\\
y_2\Theta&=(-x_1y_1+x_2y_1)(a-b)-x_2y_2(a+w)=-x_1y_1(a-b)+x_2(ay_1-ay_2)=-\Theta y_1 .
\end{align*}
For the last column, in $\mathbb Q$ one has $ax_1=h^2x_1b$, $bx_1=h^2x_1a$ and $wx_1=h^2x_1w$, so
\begin{align*}
\Theta x_1&=x_2(b-a)x_1-x_1(a+b-2w)x_1=-h^2x_2x_1(b-a)-h^2x_1^2(a+b-2w)=h^2x_1\Theta .
\end{align*}
In $\mathbb V$, $x_2$ commutes with $x_1$ and $cx_2=h^2x_2c$ for every $c\in(A_y)_2$.
\end{proof}

Let $z_\pm=y_1\pm iy_2$, a balanced basis of $A_y$, and
\begin{align*}
z&=z_+, & \rho&=z_+^{-1}z_-, & Y&=z_+^2=y_1^2-y_2^2 ,
\end{align*}
so that $z\triangleright\rho=\rho^{-1}$, and $Y$ is central in $A_y$. From the relations and Lemma~\ref{lem:Theta},
\begin{align}
\mathbb Q:&\quad z_\pm x_1=\mp ih\,x_1z_\pm, & z_\pm\Theta&=h\,\Theta z_\pm, & \Theta x_1&=h^2x_1\Theta,\label{eq:Qmodel}\\
\mathbb V:&\quad z_\pm x_2=h\,x_2z_\pm, & z_\pm\Theta&=\mp ih\,\Theta z_\pm, & \Theta x_2&=h^2x_2\Theta .\label{eq:Vmodel}
\end{align}

\begin{theorem}\label{thm:QVdeg}
Let $n$ be the normal generator ($x_1$ in $\mathbb Q$, $x_2$ in $\mathbb V$). The elements $n,\Theta,\rho,Y$ generate a quantum torus $T'$ in $Q(B)$, and $Q(B)=Q(T')\oplus Q(T')z$. If $h\in\mi$, then
\[
\PIdeg\mathbb Q(h)=\PIdeg\mathbb V(h)=4\ord(h^2).
\]
\end{theorem}

\begin{proof}
From \eqref{eq:Qmodel} and \eqref{eq:Vmodel},
\begin{align*}
\mathbb Q:&\quad \rho x_1=-x_1\rho, & \rho\Theta&=\Theta\rho, & Yx_1&=-h^2x_1Y, & Y\Theta&=h^2\Theta Y,\\
\mathbb V:&\quad \rho x_2=x_2\rho, & \rho\Theta&=-\Theta\rho, & Yx_2&=h^2x_2Y, & Y\Theta&=-h^2\Theta Y,
\end{align*}
and $Y\rho=\rho Y$. The monomials are independent. Indeed, let $x'$ be the generator other than $n$. Moving an element of $A_y$ past $x'$ produces $x'$ times an element of $A_y$ of the same degree, plus terms in $nA_y$. So by \eqref{eq:Theta} and induction, $\Theta^b\in x'^bc_b+\sum_{j<b}n^{b-j}x'^jA_y$ with $0\ne c_b\in A_y$. The elements $n^ax'^j$, $a\in\mathbb Z$, $j\ge0$, are linearly independent in $Q(A_x)$, hence right independent over $Q(A_y)$ by Lemma~\ref{lem:fracind}. Looking at the largest $b$ in a relation, the $n^a\Theta^b$, $a\in\mathbb Z$, $b\ge0$, are right independent over $Q(A_y)$ as well, and the $\rho^cY^d$ are independent in $Q(A_y)$ as in Proposition~\ref{prop:crossed}. Negative powers of $\Theta$ are removed by multiplying a relation by a power of $\Theta$.

Let $m=\ord(h^2)$. A monomial $n^a\Theta^b\rho^cY^d$ is central in $T'$ if and only if
\begin{align*}
\mathbb Q:&\quad a\ \text{even},\quad c\equiv d\ (2),\quad b\equiv-a,\ d\equiv a\ (m);\\
\mathbb V:&\quad b\ \text{even},\quad c\equiv d\ (2),\quad a\equiv-b,\ d\equiv-b\ (m).
\end{align*}
In both cases the lattice has index $4m^2$, so $\PIdeg T'=2m$.

We apply Lemma~\ref{lem:cross} with $w\in\{1,z\}$. Conjugation by $z$ multiplies $n$ and $\Theta$ by scalars, fixes $Y$ and inverts $\rho$; and $z^2=Y$. Since $\rho^2\in K$ and $\rho^2\ne\rho^{-2}$, $z$ acts nontrivially on $K$. Finally $T'$ and $z$ generate $Q(B)$: $z_-=z\rho$, the elements $y_1^2,y_2^2,y_1y_2$ lie in the span of $Y$, $Y\rho$, $Y\rho^{-1}$, and $x'$ is recovered from \eqref{eq:Theta}. So $\PIdeg B=2\cdot2m$.
\end{proof}

By Lemma~\ref{lem:cross}, the same model gives back the criterion of Proposition~\ref{prop:QV}.

\subsection{The family \texorpdfstring{$\mathbb R$}{R}}

The relations of $\mathbb R(h)$ are $x_2x_1=-x_1x_2$, $y_2y_1=-y_1y_2$ and
\begin{align*}
y_1x_1&=h(x_1y_1+x_2y_1+x_1y_2), & y_1x_2&=h\,x_1y_2,\\
y_2x_1&=h\,x_2y_1, & y_2x_2&=h(-x_2y_1-x_1y_2+x_2y_2).
\end{align*}
Let $e_\pm=x_1\pm x_2$ and $z_\pm=y_1\pm y_2$. Adding and subtracting the relations,
\begin{align}
z_+x_1&=h(e_+z_++x_2z_-), & z_+x_2&=-h\,x_2z_-, & z_-x_1&=h\,x_1z_+, & z_-x_2&=h(x_1z_+-e_-z_-).\label{eq:Rz}
\end{align}

\begin{lemma}\label{lem:Rnorm}
In $\mathbb R(h)$ let
\begin{align*}
Z_1&=x_1z_++x_2z_-=(x_1+x_2)y_1+(x_1-x_2)y_2, & Z_2&=x_1y_2+x_2y_1 .
\end{align*}
Then
\begin{align*}
Z_1x_1&=h\,x_1Z_1, & Z_1x_2&=h\,x_2Z_1, & Z_2x_1&=h\,x_2Z_2, & Z_2x_2&=-h\,x_1Z_2, & Z_2Z_1&=Z_1Z_2 .
\end{align*}
\end{lemma}

\begin{proof}
Write $a=x_1^2$, $b=x_2^2$, $u=x_1x_2$. By \eqref{eq:Rz},
\begin{align*}
Z_1x_1&=h\bigl(x_1e_+z_++x_1x_2z_-+x_2x_1z_+\bigr)=h\bigl(a\,z_++u\,z_-\bigr)=h\,x_1Z_1,\\
Z_1x_2&=h\bigl(-x_1x_2z_-+x_2x_1z_+-x_2e_-z_-\bigr)=h\bigl(-u\,z_++b\,z_-\bigr)=h\,x_2Z_1 .
\end{align*}
Note that $2Z_2=e_+z_+-e_-z_-$. Since $e_+^2-e_-x_1=x_2e_+$, $e_+x_2=-x_2e_-$, $e_-x_1=x_1e_+$ and $e_-^2-e_+x_2=x_1e_-$,
\begin{align*}
2Z_2x_1&=h\bigl(e_+^2z_++e_+x_2z_--e_-x_1z_+\bigr)=h\,x_2(e_+z_+-e_-z_-),\\
2Z_2x_2&=h\bigl(-e_+x_2z_--e_-x_1z_++e_-^2z_-\bigr)=-h\,x_1(e_+z_+-e_-z_-).
\end{align*}
For the last relation put $Y=z_+^2=z_-^2$ and $v=a+b$. By \eqref{eq:Rz}, $z_+e_+=he_+z_+$, $z_-e_-=he_-z_-$, $z_+e_-=h(e_+z_++2x_2z_-)$ and $z_-e_+=h(2x_1z_+-e_-z_-)$. Hence
\begin{align*}
2Z_1Z_2&=x_1z_+e_+z_+-x_1z_+e_-z_-+x_2z_-e_+z_+-x_2z_-e_-z_-\\
&=h\bigl((x_1e_+-2x_1x_2+2x_2x_1-x_2e_-)Y-x_1e_+\,z_+z_--x_2e_-\,z_-z_+\bigr)\\
&=h\bigl((v-2u)Y-(a+u)z_+z_-+(b+u)z_-z_+\bigr),\\
2Z_2Z_1&=e_+z_+x_1z_++e_+z_+x_2z_--e_-z_-x_1z_+-e_-z_-x_2z_-\\
&=h\bigl((e_+^2-e_+x_2-e_-x_1+e_-^2)Y+e_+x_2\,z_-z_+-e_-x_1\,z_+z_-\bigr)\\
&=h\bigl((v-2u)Y+(b+u)z_-z_+-(a+u)z_+z_-\bigr),
\end{align*}
because $x_1e_+=e_-x_1=a+u$, $-x_2e_-=e_+x_2=b+u$ and $e_\pm^2=v$.
\end{proof}

\begin{theorem}\label{thm:Rdeg}
In $Q(\mathbb R(h))$ let $X=x_1^2$, $r=x_1^{-1}x_2$ and $W=Z_2^2$. Then $X,r,Z_1,W$ generate a quantum torus $T'$ whose only nontrivial relations are $Z_1X=h^2XZ_1$ and $WX=h^4XW$, and
\[
Q(\mathbb R)=Q(T')\oplus Q(T')x_1\oplus Q(T')Z_2\oplus Q(T')x_1Z_2 .
\]
$\mathbb R(h)$ is PI if and only if $h\in\mi$, and then $\PIdeg\mathbb R(h)=4\ord(h^2)$.
\end{theorem}

\begin{proof}
$Z_1$ and $Z_2$ normalize $E=Q(A_x)$ by Lemma~\ref{lem:Rnorm}, and they are left independent over $E$, because they induce different automorphisms of the center $\kk(x_1^2,x_2^2)$ of $E$. So Lemma~\ref{lem:words} applies: the words $Z_1^aZ_2^b$ of each length are left independent over $E$. Words of different lengths lie in the spaces $E(A_y)_k$, whose sum is direct by Lemma~\ref{lem:fracind}, and negative powers of $Z_1$ and $W$ are removed by multiplying a relation by a power of $Z_1W$. Since the $X^ir^j$ are independent in $E$, the monomials of $T'$ are independent. The relations of $T'$ follow from Lemma~\ref{lem:Rnorm}: $Z_1\triangleright$ multiplies $x_1,x_2$ by $h$, and $Z_2^2\triangleright$ multiplies them by $-h^2$. A monomial $X^ir^jZ_1^kW^l$ is central if and only if $m\mid i$ and $m\mid k+2l$, where $m=\ord(h^2)$. So $\PIdeg T'=m$ when $h\in\mi$, and $T'$ is not PI otherwise.

The elements $x_1$ and $Z_2$ act on the generators by
\[
\begin{array}{c|cccc}
 & r & X & Z_1 & W\\\hline
x_1\triangleright & -r & X & h^{-1}Z_1 & -h^{-2}W\\
Z_2\triangleright & -r^{-1} & -h^2r^2X & Z_1 & W
\end{array}
\]
Here $x_1\triangleright Z_2=h^{-1}Z_2r$ and $(Z_2r)^2=-W$. Moreover $x_1^2=X$, $Z_2^2=W$ and $Z_2x_1=h\,x_1rZ_2$. The element $r$ lies in $K$, and $x_1$, $Z_2$, $x_1Z_2$ send it to $-r$, $-r^{-1}$, $r^{-1}$. Since $x_2=x_1r$ and $y_1,y_2\in EZ_1+EZ_2$, Lemma~\ref{lem:cross} applies with $N=4$.
\end{proof}

This recovers \cite[Corollary 4.8]{ZZ08} and adds the PI degree. Zhang and Zhang proved the criterion with the normal elements $x_1x_2$, $x_1^2+x_2^2$, $y_1y_2$, $y_1^2+y_2^2$ \cite[Corollary 4.7(a)]{ZZ08}; for the PI degree we need, in addition, the elements $Z_1,Z_2$ of Lemma~\ref{lem:Rnorm}, which normalize the whole plane $A_x$.

\subsection{The families \texorpdfstring{$\mathbb F$}{F} and \texorpdfstring{$\mathbb I$}{I}}

In $\mathbb F(p,h)$, with the relations of Section~\ref{sec:CFI}, let
\begin{align*}
g&=x_1+px_2, & \bar g&=x_1-px_2, & c&=h(1-p),
\end{align*}
so that $g^2=\bar g^2=x_1^2-x_2^2$, and $c^2=-2ph^2$ is the scalar $\lambda$ of Lemma~\ref{lem:Fnormal}. Using $(p-1)/(p+1)=p$, the mixed relations become
\begin{align}
y_1g&=h(1+p)\,gy_2, & y_1\bar g&=h\bigl(-2gy_1+(1-p)\bar gy_2\bigr),\notag\\
y_2\bar g&=-h(1+p)\,\bar gy_1, & y_2g&=h\bigl((1-p)gy_1+2p\,\bar gy_2\bigr).\label{eq:Fg}
\end{align}

\begin{lemma}\label{lem:Fnorm}
In $\mathbb F(p,h)$ let $Z_1=gy_1-\bar gy_2$ and $Z_2=gy_1+p\,\bar gy_2$. Then
\begin{align*}
Z_1x_1&=-c\,x_1Z_1, & Z_1x_2&=c\,x_2Z_1, & Z_2x_1&=c\,x_2Z_2, & Z_2x_2&=c\,x_1Z_2, & Z_2^2&=-p\,Z_1^2 .
\end{align*}
\end{lemma}

\begin{proof}
Let $\omega_0=g^2=\bar g^2$, which is central in $A_x$. By \eqref{eq:Fg},
\begin{align*}
Z_1g&=h(1+p)\omega_0 y_2-h\bar g\bigl((1-p)gy_1+2p\bar gy_2\bigr)=-h(1-p)\,\bar g(gy_1-\bar gy_2)=-c\,\bar gZ_1,\\
Z_1\bar g&=hg\bigl(-2gy_1+(1-p)\bar gy_2\bigr)+h(1+p)\omega_0 y_1=-h(1-p)\,g(gy_1-\bar gy_2)=-c\,gZ_1,\\
Z_2g&=h(1+p)\omega_0 y_2+hp\,\bar g\bigl((1-p)gy_1+2p\bar gy_2\bigr)=hp(1-p)\,\bar g(gy_1+p\bar gy_2)=cp\,\bar gZ_2,\\
Z_2\bar g&=hg\bigl(-2gy_1+(1-p)\bar gy_2\bigr)-hp(1+p)\omega_0 y_1=-hp(1-p)\,g(gy_1+p\bar gy_2)=-cp\,gZ_2 .
\end{align*}
The first four relations follow, since $2x_1=g+\bar g$ and $2px_2=g-\bar g$. For the last one put $G_1=gy_1$, $G_2=\bar gy_2$. By \eqref{eq:Fg} and $y_2y_1=py_1y_2$,
\begin{align*}
G_1^2&=h(p-1)\,\omega_0\,y_1y_2, & G_1G_2&=h\bigl(-2\omega_0\,y_1y_2+(1-p)g\bar g\,y_2^2\bigr),\\
G_2^2&=-h(1+p)\,\omega_0\,y_1y_2, & G_2G_1&=h\bigl((1-p)\bar gg\,y_1^2-2\omega_0\,y_1y_2\bigr).
\end{align*}
So $G_1^2-G_2^2=2hp\,\omega_0 y_1y_2=-p(G_1^2+G_2^2)$, and with $S=G_1G_2+G_2G_1$,
\begin{align*}
Z_2^2&=G_1^2-G_2^2+pS=-p\bigl(G_1^2+G_2^2-S\bigr)=-p\,Z_1^2 .\qedhere
\end{align*}
\end{proof}

The relations of $\mathbb I(q,h)$, $q^2=-1$, are $x_2x_1=q\,x_1x_2$, $y_2y_1=-y_1y_2$ and
\begin{align*}
y_1x_1&=h(-qx_1y_1-qx_2y_1+x_1y_2-qx_2y_2), & y_1x_2&=h(x_1y_1+x_2y_1+x_1y_2-qx_2y_2),\\
y_2x_1&=h(x_1y_1+qx_2y_1+qx_1y_2-qx_2y_2), & y_2x_2&=h(-x_1y_1-qx_2y_1+x_1y_2-x_2y_2).
\end{align*}
With $w_1=y_2-qy_1$ and $w_2=-(y_2+qy_1)$, so that $w_1^2=w_2^2=y_2^2-y_1^2$, they become
\begin{align}
w_1x_1&=-h(1+q)\,x_2w_1, & w_1x_2&=h\bigl((1-q)x_1w_1+2x_2w_2\bigr),\notag\\
w_2x_2&=h(1+q)\,x_1w_2, & w_2x_1&=h\bigl(-2q\,x_1w_1+(1-q)x_2w_2\bigr).\label{eq:Iw}
\end{align}

\begin{lemma}\label{lem:Inorm}
In $\mathbb I(q,h)$ let $Z_1=x_1w_1+x_2w_2$, $Z_2=x_1w_1+q\,x_2w_2$ and $c=h(1-q)$. Then
\begin{align*}
y_1Z_1&=c\,Z_1y_1, & y_2Z_1&=-c\,Z_1y_2, & y_1Z_2&=c\,Z_2y_2, & y_2Z_2&=c\,Z_2y_1, & Z_2^2&=q\,Z_1^2 .
\end{align*}
\end{lemma}

\begin{proof}
This is the computation of Lemma~\ref{lem:Fnorm} with the two planes exchanged. Let $\omega_0=w_1^2=w_2^2$, $u=x_1x_2$, $H_1=x_1w_1$ and $H_2=x_2w_2$. From \eqref{eq:Iw},
\begin{align*}
w_1Z_1&=-h(1+q)x_2\omega_0+h\bigl((1-q)x_1w_1+2x_2w_2\bigr)w_2=h(1-q)(x_1w_1+x_2w_2)w_2=c\,Z_1w_2,
\end{align*}
and in the same way $w_2Z_1=c\,Z_1w_1$, $w_1Z_2=cq\,Z_2w_2$ and $w_2Z_2=-cq\,Z_2w_1$. These give the first four relations, because $2y_2=w_1-w_2$ and $2qy_1=-(w_1+w_2)$. For the last one,
\begin{align*}
H_1^2&=-h(1+q)\,u\omega_0, & H_2^2&=hq(1+q)\,u\omega_0,
\end{align*}
so $H_1^2+q^2H_2^2=q(H_1^2+H_2^2)$, and $Z_2^2=H_1^2+q^2H_2^2+q(H_1H_2+H_2H_1)=q\,Z_1^2$.
\end{proof}

\begin{theorem}\label{thm:FIdeg}
Let $B$ be $\mathbb F(p,h)$ and $\delta^2=-p$, or $\mathbb I(q,h)$ and $\delta^2=q$, and put $Z_2'=\delta^{-1}Z_2$, so that $Z_2'^2=Z_1^2$. In $\mathbb F$ let
\begin{align*}
X&=x_1^2, & r&=x_1^{-1}x_2, & Y&=Z_1^2, & \rho&=Z_1^{-1}Z_2', & P&=\rho^2 .
\end{align*}
Then $X,r,Y,P$ generate a quantum torus $T'$ whose only nontrivial relation is $YX=c^4XY$, and
\[
Q(\mathbb F)=\bigoplus_{i,j,k\in\{0,1\}}Q(T')\,x_1^iZ_1^j\rho^k .
\]
The same holds for $\mathbb I$ with $y_2,y_1$ in place of $x_1,x_2$ and products taken in the opposite ring; the eight induced automorphisms of $K$ are again pairwise distinct. If $2h^2\in\mi$, then
\[
\PIdeg\mathbb F(p,h)=\PIdeg\mathbb I(q,h)=8\ord(-4h^4).
\]
\end{theorem}

\begin{proof}
We first treat $\mathbb F$. As in Theorem~\ref{thm:Rdeg}, $Z_1,Z_2$ normalize $E=Q(A_x)$ and are left independent over $E$. Since $Y$ commutes with $Z_1$ and $Z_2'$, every word of length $k$ in $Z_1,Z_2'$ is a power of $Y$ times an alternating word, and there are $k+1$ such products of length $k$. By Lemma~\ref{lem:words} they are left independent over $E$, and products of different lengths are independent because the spaces $E(A_y)_k$ form a direct sum (Lemma~\ref{lem:fracind}); negative powers of $Y$ and $P$ are removed by multiplying a relation by a power of $Y$, which commutes with $Z_1,Z_2'$. Now $Y^jP^l=Y^{j-2l}(Z_1Z_2')^{2l}$, because $\rho=Y^{-1}Z_1Z_2'$, so the monomials of $T'$ are independent.

By Lemma~\ref{lem:Fnorm}, $Z_1\triangleright\rho=\rho^{-1}$ and
\[
\begin{array}{c|cccc}
 & r & X & Y & P\\\hline
x_1\triangleright & -r & X & c^{-2}Y & -P\\
Z_1\triangleright & -r & c^2X & Y & P^{-1}\\
\rho\triangleright & -r^{-1} & -r^2X & Y & P
\end{array}
\]
For the first row, $x_1\triangleright Z_1=-c^{-1}Z_1$ and $x_1\triangleright Z_2'=-c^{-1}Z_2'r$, so $x_1\triangleright\rho=\rho r$, and $(\rho r)^2=-P$ because $\rho r=-r^{-1}\rho$. For the third, $\rho\triangleright x_1=x_2$ and $\rho\triangleright x_2=-x_1$. In particular $P\triangleright$ fixes $r$ and $X$, while $Y\triangleright$ fixes $r$ and multiplies $X$ by $c^4$. These are the relations of $T'$, and since $c^4=-4h^4$ we get $\PIdeg T'=\ord(-4h^4)$.

The eight elements $x_1^iZ_1^j\rho^k$ satisfy (i) and (ii) of Lemma~\ref{lem:cross}: $x_1^2=X$, $Z_1^2=Y$, $\rho^2=P$, and they commute with each other up to units of $T'$. The elements $r,P$ lie in $K$, and
\[
\begin{array}{c|ccccccc}
 & x_1 & Z_1 & \rho & x_1Z_1 & x_1\rho & Z_1\rho & x_1Z_1\rho\\\hline
r\mapsto & -r & -r & -r^{-1} & r & r^{-1} & r^{-1} & -r^{-1}\\
P\mapsto & -P & P^{-1} & P & -P^{-1} & -P & P^{-1} & -P^{-1}
\end{array}
\]
so the eight automorphisms of $K$ are distinct. Since $x_2=x_1r$, $Z_2=\delta Z_1\rho$ and $y_1,y_2\in EZ_1+EZ_2$, Lemma~\ref{lem:cross} gives $\PIdeg\mathbb F=8\ord(-4h^4)$.

For $\mathbb I$, work in the opposite ring and put $u_1=y_2$,
$u_2=y_1$. With products now taken in that ring,
Lemma~\ref{lem:Inorm} gives
\[
Z_1u_1=-c\,u_1Z_1,\quad Z_1u_2=c\,u_2Z_1,\quad
Z_2u_1=c\,u_2Z_2,\quad Z_2u_2=c\,u_1Z_2.
\]
Also $u_2u_1=-u_1u_2$ and $Z_2^2=qZ_1^2$.
After normalizing by $\delta^2=q$, these are exactly the relations
used in the calculation for $\mathbb F$.
The right-hand version of Lemma~\ref{lem:words} supplies the same
independence argument, and the torus and the eight center actions
therefore give $\PIdeg\mathbb I=8\ord(-4h^4)$ as well.
\end{proof}

Since $-4h^4=-(2h^2)^2$, Lemma~\ref{lem:cross} also recovers Theorems~\ref{thm:F} and~\ref{thm:I}.

\begin{theorem}\label{thm:FIQRV}
By Theorems~\ref{thm:QVdeg}, \ref{thm:Rdeg} and~\ref{thm:FIdeg}, for the PI members
\begin{align*}
\PIdeg\mathbb Q(h)&=\PIdeg\mathbb R(h)=\PIdeg\mathbb V(h)=4\ord(h^2),\\
\PIdeg\mathbb F(p,h)&=\PIdeg\mathbb I(q,h)=8\ord(-4h^4).
\end{align*}
\end{theorem}

\begin{remark}\label{rem:bounds}
The quantum affine space $S$ of the proof of Theorem~\ref{thm:F} gives less. In its quantum torus $\omega_1\omega_2^{-1}$ is central, and the remaining generators $\omega_2,y_1,y_2$ have commutation scalars $\lambda$ and $p$; so $\PIdeg S=|\langle\lambda,p\rangle|=\operatorname{lcm}(4,\ord(2h^2))$, and Lemma~\ref{lem:finite} with \eqref{eq:Fdec} gives $\PIdeg\mathbb F\le4\operatorname{lcm}(4,\ord(2h^2))$. This bound is attained for $2h^2=1$, where both sides are $16$, and not for $2h^2=\pm i$, where $\PIdeg\mathbb F=8$.
\end{remark}

\section{The family \texorpdfstring{$\mathbb C$}{C} and a summary of PI degrees}\label{sec:degrees}

\begin{proposition}\label{prop:bounds}
Let $\mathbb C(p,h)$ be PI, and put $\beta=3(2p+1)h^3$ and $n=\operatorname{lcm}\bigl(3,\ord\beta\bigr)$. Then
\[
\PIdeg\mathbb C(p,h)\in\{n,\ 3n,\ 9n\}.
\]
\end{proposition}

\begin{proof}
Keep $a=y_1^2y_2$, $b=y_1^3-y_2^3$ and the quantum affine space $S=\kk\langle x_1,x_2,a,b\rangle$ of the proof of Theorem~\ref{thm:C}. Write $p=\zeta^{k}$ and $\beta=\zeta^{l}$ for a root of unity $\zeta$, so that $\alpha=p^2\beta=\zeta^{l+2k}$. The commutation matrix of $S$ in the order $x_1,x_2,a,b$ has entries $-k,\ -(l+2k),\ -l,\ -(l+2k),\ -l,\ 0$ and Pfaffian $0$, so $\PIdeg S=|\langle p,\beta\rangle|=n$.

Next, $A_y$ is a free left $\kk[a,b]$-module of rank nine. The element $a$ is normal and regular in $A_y$, and the monomials $y_1^iy_2^j$ with $i\le1$ or $j=0$ give a basis of $A_y/aA_y$. Let $v=\sum c_{ij}y_1^iy_2^j$ be a combination of these monomials with $bv\in aA_y$. Since $y_1^{i+3}y_2^j\in aA_y$ for $j\ge1$ and $y_2^3$ commutes with $y_1$,
\[
bv\equiv\sum_ic_{i0}\,y_1^{i+3}-\sum_{i\le1,\ j\ge0}c_{ij}\,y_1^iy_2^{j+3}\pmod{aA_y}.
\]
The monomials of the first sum have $y_2$-exponent $0$, those of the second have $y_2$-exponent at least $3$, and all of them belong to the basis of $A_y/aA_y$. So every $c_{ij}$ vanishes, and $b$ is regular on $A_y/aA_y$. As $b$ is central, $A_y/(aA_y+bA_y)$ has Hilbert series
\[
\frac{(1-t^3)^2}{(1-t)^2}=(1+t+t^2)^2 .
\]
Lifting a homogeneous basis gives $w_1,\dots,w_9\in A_y$ which generate $A_y$ as a left $\kk[a,b]$-module by the graded Nakayama lemma. The free module on them has Hilbert series $(1+t+t^2)^2(1-t^3)^{-2}=(1-t)^{-2}$, the series of $A_y$, so the $w_i$ are a basis.

Since $S=A_x\,\kk[a,b]$ and $A_x\otimes A_y\to\mathbb C$ is bijective, $\mathbb C=\bigoplus_{i=1}^9Sw_i$. Clearing a common left denominator shows that the $w_i$ remain left independent over $Q(S)$, so $V=\sum_iQ(S)w_i=Q(S)\,\mathbb C$ is a left $Q(S)$-space of dimension nine, stable under right multiplication by $\mathbb C$. For $0\ne c\in\mathbb C$, right multiplication by $c$ is an injective $Q(S)$-linear map of $V$, hence bijective, so $V$ is also stable under right multiplication by $c^{-1}$. Since $1\in V$, we get $V=Q(\mathbb C)$.

Put $\Delta=Q(\mathbb C)$, $D=Q(S)$, $K=Z(D)$, $F=Z(\Delta)$ and $N=\PIdeg\mathbb C$, so $[D:K]=n^2$ and $[\Delta:F]=N^2$. The subring $L=FK$ of $\Delta$ is a field, finite over $F$. The algebra $D\otimes_KL$ is simple, so $DL\cong D\otimes_KL$ is central simple over $L$ of dimension $n^2$. By the double centralizer theorem its centralizer $\Delta'$ in $\Delta$ is a division algebra with center $L$, say $[\Delta':L]=\ell^2$, and $[\Delta:F]=[DL:F]\,[\Delta':F]$. Counting dimensions over $F$ and over $K$,
\begin{align*}
N^2&=n^2[L:F]\cdot\ell^2[L:F], & 9n^2&=\dim_K\Delta=[L:K]\,\frac{N^2}{[L:F]} ,
\end{align*}
where $\dim_K\Delta=\dim_D\Delta\cdot[D:K]$ is the dimension of $\Delta$ as a left $K$-space. Hence
\begin{equation}\label{eq:Cindex}
N=n\,\ell\,[L:F],\qquad 9=[L:K]\,[L:F]\,\ell^2 .
\end{equation}
If $\ell=1$ then $[L:F]$ divides $9$, and if $\ell=3$ then $[L:F]=1$. In both cases $N/n\in\{1,3,9\}$.
\end{proof}

In particular $n\le\PIdeg\mathbb C(p,h)\le9n$; the lower bound is just $\PIdeg S$, and the upper bound also follows from Lemma~\ref{lem:finite}.

Table~\ref{tab:degrees} collects what is known. In the formulas for $\mathbb D$ and $\mathbb G$, $\zeta$ is a primitive $N$th root of unity, $p=\zeta^k$ and $h=\zeta^{k'}$; the PI degree of $\mathbb G$ does not depend on $f$.

\begin{longtable}{|c|l|l|}
\caption{PI degrees of the PI members.}\label{tab:degrees}\\
\hline
Family & $\PIdeg$ & Reference\\\hline
\endfirsthead
\hline
Family & $\PIdeg$ & Reference\\\hline
\endhead
$\mathbb D$ & $N^2/\gcd\bigl(N,2k(k+k')\bigr)$, $N=|\langle-1,p,h\rangle|$ & Prop.~\ref{prop:DGdeg}\\
$\mathbb G$ & $N^2/\gcd\bigl(N,2k(k+k')\bigr)$, $N=|\langle p,h\rangle|$ & Prop.~\ref{prop:DGdeg}, \cite{GV-G}\\
$\mathbb K$ & $2|\langle-1,q,f,h\rangle|$ & \cite{GV-K}, \cite[Thm.~3.2]{Rub-ozone}\\
$\mathbb L$ & $2N$ or $N$, $N=|\langle-1,q,f,h\rangle|$ & Thm.~\ref{thm:L}\\
$\mathbb B$, $\mathbb E$, $\mathbb J$, $\mathbb M$, $\mathbb N$, $\mathbb O$, $\mathbb P$ & see Thm.~\ref{thm:monodeg} & Thm.~\ref{thm:monodeg}\\
$\mathbb S$, $\mathbb T$, $\mathbb U$, $\mathbb W$, $\mathbb Z$ & see Thm.~\ref{thm:monodeg} & Thm.~\ref{thm:monodeg}\\
$\mathbb F$, $\mathbb I$ & $8\ord(-4h^4)$ & Thm.~\ref{thm:FIQRV}\\
$\mathbb Q$, $\mathbb R$, $\mathbb V$ & $4\ord(h^2)$ & Thm.~\ref{thm:FIQRV}\\
$\mathbb C$ & $n$, $3n$ or $9n$, $n=\operatorname{lcm}(3,\ord\beta)$ & Prop.~\ref{prop:bounds}\\
\hline
\end{longtable}

\begin{question}\label{q:degrees}
Which of $n$, $3n$, $9n$ is the PI degree of $\mathbb C(p,h)$? Equivalently, by \eqref{eq:Cindex}, how do the centers of $Q(S)$ and of $Q(\mathbb C)$ sit inside $Q(\mathbb C)$?
\end{question}

As far as we can see, the method of Section~\ref{sec:norm} does not reach $\mathbb C$: we know no element of low bidegree normalizing one of its planes, and both planes are $\kk_p$ with $p$ of order three, so their quotient rings have degree three, not two.

\section{The determinant of \texorpdfstring{$\sigma$}{sigma}}\label{sec:det}

Zhang and Zhang attach to a double extension the automorphism $\det\sigma=-p_{11}\sigma_{12}\sigma_{11}+\sigma_{22}\sigma_{11}-p_{12}\sigma_{12}\sigma_{21}$ of $A_x$ \cite[Section 2]{ZZ08}, and for family $\mathbb R$ they prove that $B$ is PI if and only if $\det\sigma$ has finite order \cite[Corollary 4.8]{ZZ08}. The published version of \cite{ZZ09} computes $\det\sigma$ only for $\mathbb A$. On degree one, the matrix for every family is listed in Section 4 of the arXiv version of \cite{ZZ09}, and we reproduce that list in Table~\ref{tab:det}.

\begin{longtable}{|c|c||c|c|}
\caption{The automorphism $\det\sigma$ in degree one (Section 4 of the arXiv version of \cite{ZZ09}).}\label{tab:det}\\
\hline
Family & $\det\sigma$ & Family & $\det\sigma$\\\hline
\endfirsthead
\hline
Family & $\det\sigma$ & Family & $\det\sigma$\\\hline
\endhead
$\mathbb A$ & $h^2I$ & $\mathbb N$ & $(f^2-g^2)I$\\
$\mathbb B$ & $ph^2\operatorname{diag}(1,-1)$ & $\mathbb O$ & $(f-1)h^2I$\\
$\mathbb C$ & $-3h^2\operatorname{diag}(p,1)$ & $\mathbb P$ & $(1-f)h^2I$\\
$\mathbb D$ & $-p^2h^2I$ & $\mathbb Q$ & $h^2\operatorname{diag}(-1,1)$\\
$\mathbb E$ & $2ph^2\begin{psmallmatrix}0&1\\-1&0\end{psmallmatrix}$ & $\mathbb R$ & $h^2\begin{psmallmatrix}0&1\\-1&0\end{psmallmatrix}$\\
$\mathbb F$ & $-2ph^2I$ & $\mathbb S$ & $4h^2I$\\
$\mathbb G$ & $p^2h^2I$ & $\mathbb T$ & $4h^2\begin{psmallmatrix}0&1\\1&0\end{psmallmatrix}$\\
$\mathbb H$ & $h^2\begin{psmallmatrix}1&0\\2f&1\end{psmallmatrix}$ & $\mathbb U$ & $4h^2I$\\
$\mathbb I$ & $2h^2\operatorname{diag}(1,-1)$ & $\mathbb V$ & $h^2\begin{psmallmatrix}-1&1\\0&1\end{psmallmatrix}$\\
$\mathbb J$ & $2h^2I$ & $\mathbb W$ & $(f+1)h^2I$\\
$\mathbb K$ & $h^2\operatorname{diag}(1,f)$ & $\mathbb X$ & $h^2\begin{psmallmatrix}1&0\\2&1\end{psmallmatrix}$\\
$\mathbb L$ & $h^2\operatorname{diag}(f^2,1)$ & $\mathbb Y$ & $h^2I$\\
$\mathbb M$ & $(1-f)h^2I$ & $\mathbb Z$ & $-(f+1)h^2I$\\
\hline
\end{longtable}

\begin{corollary}\label{cor:det}
If a member of the list is PI, then $\det\sigma$ has finite order.
\end{corollary}

\begin{proof}
Compare Tables~\ref{tab:main} and~\ref{tab:det}. In every family that has PI members the matrix is a scalar $\gamma$ times a matrix of finite order, or a diagonal matrix, and the scalars involved are roots of unity when the PI criterion holds: $\gamma$ is $\pm h^2$, $\pm p^2h^2$ or $ph^2$ in $\mathbb B$, $\mathbb D$, $\mathbb G$, $\mathbb Q$, $\mathbb R$, $\mathbb V$; it is $\pm2h^2$ or $\pm2ph^2$ in $\mathbb E$, $\mathbb F$, $\mathbb I$, $\mathbb J$; it is $4h^2=(2h)^2$ in $\mathbb S$, $\mathbb T$, $\mathbb U$; the entries are $-3h^2p$, $-3h^2$ in $\mathbb C$ and $h^2$, $h^2f$, $h^2f^2$ in $\mathbb K$, $\mathbb L$; and
\begin{align*}
(f-1)h^2&=-h(1-s)\cdot h(1+s), & (f+1)h^2&=h(s+i)\cdot h(s-i), & f^2-g^2&=(f-g)(f+g)
\end{align*}
in $\mathbb M$, $\mathbb O$, $\mathbb P$, in $\mathbb W$, $\mathbb Z$ and in $\mathbb N$.
\end{proof}

\begin{remark}\label{rem:detsigma}
The converse fails, so \cite[Corollary 4.8]{ZZ08} does not extend to the list. For $h\in\mi$, $\det\sigma=h^2I$ in $\mathbb Y(f,h)$ and $\det\sigma=-h^2I$ in $\mathbb D(\pm1,h)$, but these algebras are not PI (Proposition~\ref{prop:noPI}). For a member that is not degenerate, take $\mathbb O$ with $s=3$ and $h=(2\sqrt2)^{-1}$: then $\det\sigma=h^2(f-1)I=I$, while $h(1+s)=\sqrt2\notin\mi$. The twist is not to blame: in $\mathbb N$, where there is no $h$, the member with $f=2$, $g=\sqrt3$ has $\det\sigma=(f^2-g^2)I=I$ and is not PI, because $f+g=2+\sqrt3\notin\mi$. In $\mathbb N$ the determinant records the products in $\mathcal I_{\rm III}$ but not the ratios. In $\mathbb O$ it does not even record the two elements of $\mathcal I_{\rm I}$: in the example above they are $\lambda_{11}\lambda_{22}=\frac12$ and $\lambda_{12}\lambda_{21}=2$, neither of finite order, and $\det\sigma$ only sees their product.
\end{remark}

\section{Primitivity}\label{sec:primitive}

\begin{proposition}\label{prop:notprimitive}
Let $B=\bigoplus_{n\ge0}B_n$ be a connected graded PI algebra with $B_1\ne0$. Then $B$ is not primitive. In particular no PI member of the list is primitive.
\end{proposition}

\begin{proof}
By Kaplansky's theorem a primitive PI ring is simple artinian \cite[Theorem 13.3.8]{MR}, while $B_{\ge1}$ is a proper nonzero ideal.
\end{proof}

Connected graded algebras can be primitive, so the non-PI members need another argument. By \cite[Proposition 3.8]{ZZ09} the opposite ring of a member of the list is again a connected graded double extension of a regular algebra of dimension two, so it is also strongly noetherian \cite[Theorem 0.1]{ZZ09}. Every argument below therefore applies on both sides, and \emph{primitive} means left and right primitive. The obstruction we use is this.

\begin{proposition}\label{prop:rational}
Let $B$ be a member of the list. If the center of $Q(B)$ is not algebraic over $\kk$, then $B$ is not primitive.
\end{proposition}

\begin{proof}
$B$ is a noetherian domain, strongly noetherian by \cite[Theorem 0.1]{ZZ09}, and connected graded with finite-dimensional components. So it is a Jacobson algebra satisfying the Nullstellensatz \cite[Theorem 0.4]{ASZ}. For such algebras a primitive ideal is rational \cite[Lemma II.7.15]{BG}, see also \cite{GL}; in particular, if $(0)$ is primitive, then $Z(Q(B))$ is algebraic over $\kk$.
\end{proof}

\subsection{The family \texorpdfstring{$\mathbb B$}{B}}

The relations of $\mathbb B(p,h)$, $p^2=-1$, are $x_2x_1=p\,x_1x_2$, $y_2y_1=p\,y_1y_2$ and
\[
y_1x_1=h\,x_2y_2,\qquad y_1x_2=h\,x_1y_2,\qquad y_2x_1=-h\,x_2y_1,\qquad y_2x_2=h\,x_1y_1 .
\]

\begin{proposition}\label{prop:Bnotprim}
For every $h\in\kk^\times$ the element
\[
w=(x_1^4+x_2^4)^2\,(x_1^4x_2^4)^{-1}\in Q(\mathbb B)
\]
is central and transcendental over $\kk$. In particular no member of $\mathbb B$ is primitive.
\end{proposition}

\begin{proof}
Applying the mixed relations twice,
\begin{align*}
y_1x_1^2&=h\,x_2y_2x_1=-h^2x_2^2y_1, & y_2x_1^2&=-h\,x_2y_1x_1=-h^2x_2^2y_2,\\
y_1x_2^2&=h\,x_1y_2x_2=h^2x_1^2y_1, & y_2x_2^2&=h\,x_1y_1x_2=h^2x_1^2y_2,
\end{align*}
so $yx_1^2=-h^2x_2^2y$ and $yx_2^2=h^2x_1^2y$ for $y\in\{y_1,y_2\}$. Applying each rule twice,
\begin{equation}\label{eq:Bswap}
y\,x_1^4=h^4\,x_2^4\,y,\qquad y\,x_2^4=h^4\,x_1^4\,y .
\end{equation}
Since $p^4=1$, $x_1^4$ and $x_2^4$ are central in $A_x$. By \eqref{eq:Bswap}, $y(x_1^4+x_2^4)^2=h^8(x_1^4+x_2^4)^2y$ and $y\,x_1^4x_2^4=h^8\,x_2^4x_1^4y=h^8\,x_1^4x_2^4y$, so $w$ commutes with $y_1,y_2$, and also with $x_1,x_2$. By the PBW basis $x_1^4$ and $x_2^4$ are algebraically independent, so $w$ is a nonconstant element of $\kk(x_1^4,x_2^4)$, hence transcendental over $\kk$. Proposition~\ref{prop:rational} finishes the proof.
\end{proof}

\subsection{The families \texorpdfstring{$\mathbb C$}{C}, \texorpdfstring{$\mathbb F$}{F}, \texorpdfstring{$\mathbb I$}{I}, \texorpdfstring{$\mathbb Q$}{Q}, \texorpdfstring{$\mathbb R$}{R}, \texorpdfstring{$\mathbb V$}{V}}

\begin{proposition}\label{prop:CFIQRV}
For every value of the twist the following elements are central in $Q(B)$ and transcendental over $\kk$:
\begin{align*}
\mathbb F,\ \mathbb R:&\quad (x_1x_2)^2(x_1^2+x_2^2)^{-2}, & \mathbb I:&\quad (y_1y_2)^2(y_1^2+y_2^2)^{-2}, & \mathbb C:&\quad (y_1^2y_2)^3(y_1^3-y_2^3)^{-3},\\
\mathbb Q:&\quad t^2(2-t^2),\ \ t=x_1^{-1}x_2, & \mathbb V:&\quad (t^2-t)^2,\ \ t=x_2^{-1}x_1 . &&
\end{align*}
In particular no member of these six families is primitive.
\end{proposition}

\begin{proof}
$\mathbb F$: by Lemma~\ref{lem:Fnormal}, $\xi=\omega_1\omega_2^{-1}$ commutes with $y_1,y_2$. Now $\omega_1^2=-x_1^2x_2^2$ and $\omega_2$ are central in $A_x$, so $\xi^2$ is central in $Q(\mathbb F)$. For $\mathbb I$ use the elements $\omega_1,\omega_2$ of the proof of Theorem~\ref{thm:I}, and for $\mathbb R$ the relations $y_j\omega_1=-h^2\omega_1y_j$, $y_j\omega_2=h^2\omega_2y_j$ for $\omega_1=x_1x_2$, $\omega_2=x_1^2+x_2^2$, which are in the proof of \cite[Corollary 4.7(a)]{ZZ08}.

$\mathbb C$: by \eqref{eq:Cnormal} and $\alpha^3=\beta^3$, the element $a^3b^{-3}$ commutes with $x_1,x_2$; and $a^3=y_1^6y_2^3$ and $b$ are central in $A_y$.

$\mathbb Q$, $\mathbb V$: the proof of Proposition~\ref{prop:QV} shows that $F_{\mathbb Q}(t)$ and $F_{\mathbb V}(t)$ commute with $y_1,y_2$, and \eqref{eq:Mt} does not involve $h$. In $\mathbb Q$, $t$ anticommutes with $x_1$ and $F_{\mathbb Q}$ is even; in $\mathbb V$, $t$ commutes with $x_1,x_2$.

Each element is a nonconstant element of a purely transcendental extension of $\kk$, namely $\kk(x_1^2,x_2^2)$, $\kk(y_1^2,y_2^2)$, $\kk(y_1^3,y_2^3)$ or $\kk(t)$. Proposition~\ref{prop:rational} finishes the proof.
\end{proof}

\subsection{The families \texorpdfstring{$\mathbb D$}{D} and \texorpdfstring{$\mathbb G$}{G}}

We use the presentations of Section~\ref{sec:normal}; for $\mathbb G$ this is the presentation of \cite{GV-G}. Let $p^2\ne1$ and $q=hp$; throughout this subsection and the next,
$q$ denotes $hp$. The relations of $\mathbb D(p,h)$ are
\begin{align*}
x_2x_1&=-x_1x_2, & y_2y_1&=p\,y_1y_2, & y_1x_1&=-q\,x_1y_1,\\
y_2x_1&=q\,x_1y_2, & y_1x_2&=h(-p^2x_2y_1+x_1y_2), & y_2x_2&=h(x_2y_2+x_1y_1),
\end{align*}
and those of $\mathbb G(p,f,h)$ are
\begin{align*}
x_2x_1&=x_1x_2, & y_ix_1&=q\,x_1y_i, & y_2y_1&=p\,y_1y_2,\\
y_1x_2&=h(p^2x_2y_1+x_1y_2), & y_2x_2&=h(x_2y_2+fx_1y_1). &&
\end{align*}
For all $p,h$ with $p^2\ne1$ (and all admissible $f$), $Q(\mathbb D)$ and $Q(\mathbb G)$ are the quotient division rings of the quantum tori on $x_1,\vartheta,y_1,y_2$ with
\begin{align}
\mathbb D:&\quad \vartheta x_1=-x_1\vartheta,\ \ \vartheta y_1=p^{-1}y_1\vartheta,\ \ \vartheta y_2=py_2\vartheta,\label{eq:Dtorus}\\
&\quad y_1x_1=-q\,x_1y_1,\ \ y_2x_1=q\,x_1y_2,\ \ y_2y_1=p\,y_1y_2,\notag\\
\mathbb G:&\quad \vartheta x_1=x_1\vartheta,\ \ \vartheta y_1=p^{-1}y_1\vartheta,\ \ \vartheta y_2=py_2\vartheta,\ \ y_ix_1=q\,x_1y_i,\ \ y_2y_1=p\,y_1y_2 .\label{eq:Gtorus}
\end{align}
Here $\vartheta$ is $x_1^{-1}x_2$ minus an explicit Laurent polynomial in $y_1y_2^{-1}$; for $\mathbb G$ it is called $T$ in \cite{GV-G}, where this model was found \cite[Theorem 3.2]{GV-G}. A proof is in Lemma~\ref{lem:DGmodel} below.

\begin{lemma}\label{lem:centerQT}
Let $\mathcal T$ be a quantum torus on $v_1^{\pm1},\dots,v_n^{\pm1}$, with $v_iv_j=c_{ij}v_jv_i$, and let $\Gamma\subseteq\mathbb Z^n$ be the set of exponents of central monomials. Then $Z(Q(\mathcal T))$ is the fraction field of $Z(\mathcal T)=\kk[v^{\mathbf a}:\mathbf a\in\Gamma]$, and its transcendence degree over $\kk$ is the rank of $\Gamma$.
\end{lemma}

\begin{proof}
For $\mathbf b\in\mathbb Z^n$ let $\chi_{\mathbf b}$ be the character with $v^{\mathbf b}v^{\mathbf a}v^{-\mathbf b}=\chi_{\mathbf b}(\mathbf a)v^{\mathbf a}$. Two monomials $v^{\mathbf a},v^{\mathbf a'}$ have the same characters if and only if $\mathbf a-\mathbf a'\in\Gamma$. So an element $\phi\in\mathcal T$ with $v^{\mathbf b}\phi v^{-\mathbf b}=\chi_{\mathbf b}(\mathbf a_0)\phi$ for all $\mathbf b$ is of the form $v^{\mathbf a_0}c$ with $c\in Z(\mathcal T)$.

Let $0\ne z\in Z(Q(\mathcal T))$ and $I=\{\phi\in\mathcal T:z\phi\in\mathcal T\}$. Since $z$ is central, $I$ is a two-sided ideal, and it is nonzero because $z=a\phi^{-1}$ for some $a,\phi\in\mathcal T$. It is stable under conjugation by monomials, since $z\,v^{\mathbf b}\phi v^{-\mathbf b}=v^{\mathbf b}(z\phi)v^{-\mathbf b}$. Pick $\phi\in I$ nonzero with the fewest monomials, and $\mathbf a_0$ in its support. Then $v^{\mathbf b}\phi v^{-\mathbf b}-\chi_{\mathbf b}(\mathbf a_0)\phi$ lies in $I$ and has fewer monomials, so it vanishes, and $\phi=v^{\mathbf a_0}c$ with $c\in Z(\mathcal T)$. The element $z\phi\in\mathcal T$ satisfies the same rule, so $z\phi=v^{\mathbf a_0}c'$ with $c'\in Z(\mathcal T)$. Hence $z=v^{\mathbf a_0}c'c^{-1}v^{-\mathbf a_0}=c'c^{-1}$. Finally $Z(\mathcal T)$ is a Laurent polynomial ring in $\operatorname{rk}\Gamma$ variables, by choosing a basis of $\Gamma$.
\end{proof}

\begin{theorem}\label{thm:DG}
Let $B$ be $\mathbb D(p,h)$ with $p^2\ne1$, or $\mathbb G(p,f,h)$, and $q=hp$.
\begin{enumerate}
\item[(i)] The transcendence degree of $Z(Q(B))$ over $\kk$ is $4$ if $p,q\in\mi$, $2$ if exactly one of $p,q$ is in $\mi$, and $0$ if $p,q\notin\mi$. In the last case $Z(Q(B))=\kk$.
\item[(ii)] If $q\in\mi$, then $x_1^{m}$, $m=\operatorname{lcm}(2,\ord q)$, is central in $\mathbb D(p,h)$, and $x_1^{\ord q}$ is central in $\mathbb G(p,f,h)$.
\item[(iii)] If $B$ is not PI and $p\in\mi$ or $q\in\mi$, then $B$ is not primitive.
\end{enumerate}
\end{theorem}

\begin{proof}
(i) Commuting $x_1^a\vartheta^by_1^cy_2^d$ past $x_1,\vartheta,y_1,y_2$ with \eqref{eq:Dtorus} multiplies it by the inverses of
\begin{equation}\notag
(-1)^{b+c}q^{-(c+d)},\qquad (-1)^ap^{d-c},\qquad (-q)^ap^{b-d},\qquad q^ap^{c-b},
\end{equation}
and doing the same for $x_1^a\vartheta^by_1^cy_2^d$ with \eqref{eq:Gtorus} gives
\begin{equation}\notag
q^{-(c+d)},\qquad p^{d-c},\qquad q^ap^{b-d},\qquad q^ap^{c-b}.
\end{equation}
The lattice $\Gamma$ is cut out by setting these four scalars equal to $1$. If $p,q\in\mi$, it has finite index and rank four. Suppose $p\notin\mi$. The second condition gives $d=c$, and in $\mathbb D$ also $a$ even. The last two give $q^ap^{c-b}=q^ap^{b-c}=1$, hence $p^{2(b-c)}=1$, $b=c$ and $q^a=1$. The first then reads $q^{2c}=1$. So
\[
\Gamma=\{(a,c,c,c):\ q^a=q^{2c}=1,\ a\text{ even in }\mathbb D\},
\]
which has rank two if $q\in\mi$ and is zero if $q\notin\mi$. Suppose $p\in\mi$ and $q\notin\mi$. The last condition forces $q^a\in\mi$, so $a=0$; the first forces $c+d=0$, and in $\mathbb D$ also $b+c$ even. What remains is $p^{b+c}=p^{c-b}=1$, so
\[
\Gamma=\{(0,b,c,-c):\ p^{b+c}=p^{c-b}=1,\ b+c\text{ even in }\mathbb D\},
\]
which contains $(0,2n,0,0)$ and $(0,0,2n,-2n)$, $n=\ord p$, and has rank two. Lemma~\ref{lem:centerQT} gives (i).

(ii) In $\mathbb D$, $m$ is even and $q^m=1$, so $x_1^m$ commutes with $x_2$, $y_1$ and $y_2$ by the first, third and fourth relations. In $\mathbb G$, $x_1$ commutes with $x_2$ and $y_ix_1^{\ord q}=q^{\ord q}x_1^{\ord q}y_i$.

(iii) If $B$ is not PI, then $p$ and $q$ are not both roots of unity, so by (i) $Z(Q(B))$ has transcendence degree two. Proposition~\ref{prop:rational} applies.
\end{proof}

When $q\in\mi$ and $p\notin\mi$, part (ii) gives an algebra that is neither PI nor primitive and still has central elements of positive degree.

\subsection{Primitivity of \texorpdfstring{$\mathbb D$}{D} and \texorpdfstring{$\mathbb G$}{G} when \texorpdfstring{$p,hp\notin\mi$}{p, hp not roots of unity}}

Now let $p,q\notin\mi$. Then $Z(Q(B))=\kk$ and Proposition~\ref{prop:rational} gives nothing. We show that $B$ is primitive, treating $\mathbb D$ and $\mathbb G$ together through the part of $B$ over $\kk[x_1^{-1}x_2]$. For the construction of $A$ and the quantum torus model of $Q(B)$, we temporarily allow all admissible parameters with $p^2\ne1$, including the case in which $p$ and $h$ are roots of unity.

Let $\beta_1,\beta_2\in\kk^\times$ and $p^2\ne1$, and let $A$ be the algebra generated by $t,y_1,y_2$ with
\begin{equation}\label{eq:Arel}
y_2y_1=p\,y_1y_2,\qquad ty_1=p^{-1}y_1t+\beta_1y_2,\qquad ty_2=py_2t+\beta_2y_1 .
\end{equation}
The defining relations identify $A$ as an Ore extension
$\kk_p[y_1,y_2][t;\sigma,\delta]$, with $\sigma$ an
automorphism; hence $A$ is noetherian.

Put
\begin{equation}\notag
b=-\frac{\beta_2}{p^2-1},\qquad c=\frac{p^2\beta_1}{p^2-1},\qquad a=bc,\qquad \lambda(w)=w+aw^{-1},
\end{equation}
and
\[
\Omega=ty_1y_2-p^{-1}b\,y_1^2-c\,y_2^2 .
\]
The constants satisfy
\begin{equation}\label{eq:aident}
a(p-p^{-1})^2=-\beta_1\beta_2,\qquad a(p^2-1)=p^2\beta_1b .
\end{equation}
In both families let $t=x_1^{-1}x_2$. For $\mathbb D$ take $\beta_1=p^{-2}$, $\beta_2=-1$; then $b=c=\kappa_0:=(p^2-1)^{-1}$ and $\Omega=x_1^{-1}\Xi$ with $\Xi$ as in Remark~\ref{rem:DGnormal} below. For $\mathbb G$ take $\beta_1=-p^{-2}$, $\beta_2=-f$; then $b,c$ are the constants of \cite{GV-G} and $\Omega=x_1^{-1}\Xi$ with $\Xi=x_2y_1y_2-x_1(cy_2^2+p^{-1}by_1^2)$. In both cases \eqref{eq:Arel} are the relations \eqref{eq:Mt} between $t$ and the $y$'s, with the matrices \eqref{eq:MDG}, and they do not involve $h$. By the PBW bases of $\mathbb D$ and $\mathbb G$, $A$ is the subalgebra of $B[x_1^{-1}]$ generated by $t,y_1,y_2$, it has basis $t^iy_1^jy_2^l$, and
\begin{equation}\label{eq:Bx}
B[x_1^{-1}]=A[x_1^{\pm1};\alpha],\qquad \alpha=x_1(-)x_1^{-1},
\end{equation}
We grade $A$ by $\deg t=0$, $\deg y_i=1$; the relations are homogeneous, and $A_n$ is a free left $\kk[t]$-module on the $y_1^jy_2^l$, $j+l=n$.

\begin{lemma}\label{lem:Omega}
$\Omega$ is central in $A$. In $A$ we also have $y_1t=p\,ty_1-p\beta_1y_2$ and $y_2t=p^{-1}ty_2-p^{-1}\beta_2y_1$.
\end{lemma}

\begin{proof}
The last two relations are \eqref{eq:Arel} solved for $y_it$. In $Q(A)$ put $\varpi=y_1y_2^{-1}$ and $\vartheta=t-b\varpi-c\varpi^{-1}$. From $y_2^{-1}y_1=p^{-1}y_1y_2^{-1}$ one gets $\varpi y_1=p^{-1}y_1\varpi$, $\varpi^{-1}y_1=y_2$, $y_1\varpi^{-1}=p^{-1}y_2$, $\varpi y_2=y_1$, $y_2\varpi=py_1$, hence
\begin{align*}
(b\varpi+c\varpi^{-1})y_1-p^{-1}y_1(b\varpi+c\varpi^{-1})&=c(1-p^{-2})y_2=\beta_1y_2,\\
(b\varpi+c\varpi^{-1})y_2-py_2(b\varpi+c\varpi^{-1})&=b(1-p^2)y_1=\beta_2y_1,
\end{align*}
and subtracting from \eqref{eq:Arel} gives $\vartheta y_1=p^{-1}y_1\vartheta$ and $\vartheta y_2=py_2\vartheta$. Since $\varpi y_1y_2=p^{-1}y_1^2$ and $\varpi^{-1}y_1y_2=y_2^2$, we have $\Omega=\vartheta y_1y_2$. Now $(y_1y_2)y_1=p\,y_1(y_1y_2)$ and $(y_1y_2)y_2=p^{-1}y_2(y_1y_2)$, so $\vartheta y_1y_2$ commutes with $y_1$ and $y_2$; it commutes with $\vartheta$, and with $\varpi$ because $(y_1y_2)\varpi=p^2\varpi(y_1y_2)$ and $\vartheta\varpi=p^{-2}\varpi\vartheta$. Hence it commutes with $t=\vartheta+b\varpi+c\varpi^{-1}$.
\end{proof}

\begin{lemma}\label{lem:DGmodel}
Let $B$ be $\mathbb D(p,h)$ with $p^2\ne1$, or $\mathbb G(p,f,h)$, and $\vartheta=t-b\varpi-c\varpi^{-1}\in Q(B)$ as in the proof of Lemma~\ref{lem:Omega}. Then $x_1,\vartheta,y_1,y_2$ satisfy \eqref{eq:Dtorus}, respectively \eqref{eq:Gtorus}, their monomials are linearly independent, and they generate $Q(B)$ as a division ring.
\end{lemma}

\begin{proof}
The relations of $\vartheta$ with $y_1,y_2$ are in the proof of Lemma~\ref{lem:Omega}. In $\mathbb D$,
\begin{align*}
x_1tx_1^{-1}&=-t, & x_1y_1x_1^{-1}&=-q^{-1}y_1, & x_1y_2x_1^{-1}&=q^{-1}y_2,
\end{align*}
so $x_1\varpi x_1^{-1}=-\varpi$ and $x_1\vartheta x_1^{-1}=-\vartheta$. In $\mathbb G$, $x_1$ commutes with $t$ and $x_1y_ix_1^{-1}=q^{-1}y_i$, so $x_1$ commutes with $\varpi$ and with $\vartheta$. Since $t=\vartheta+b\varpi+c\varpi^{-1}$ and $x_2=x_1t$, the four elements generate $Q(B)$.

For the independence, note that by Lemma~\ref{lem:Omega} there are maps $\sigma',\delta'$ of $A_y$ with $st=t\sigma'(s)+\delta'(s)$ for $s\in A_y$; then $s^{-1}t=t\sigma'(s)^{-1}-s^{-1}\delta'(s)\sigma'(s)^{-1}$ for $s\ne0$, so $Q(A_y)\,t\subseteq t\,Q(A_y)+Q(A_y)$, so $\vartheta^k\in t^k+\sum_{i<k}t^iQ(A_y)$. The elements $x_1^at^i$, $a\in\mathbb Z$, $i\ge0$, are nonzero multiples of the distinct monomials $x_1^{a-i}x_2^i$ of $Q(A_x)$, so they are right independent over $Q(A_y)$ by Lemma~\ref{lem:fracind}. Looking at the largest power of $\vartheta$ in a relation, the same holds for the $x_1^a\vartheta^k$, and therefore the $x_1^a\vartheta^ky_1^jy_2^l$, $a,j,l\in\mathbb Z$, $k\ge0$, are linearly independent. Negative powers of $\vartheta$ are handled by multiplying a relation by a power of $\vartheta$, which commutes with the other three generators up to scalars.
\end{proof}

\begin{proposition}\label{prop:DGdeg}
Let $B$ be $\mathbb D(p,h)$ with $p\ne\pm1$, or $\mathbb G(p,f,h)$, and assume $p,h\in\mi$. Put $N=|\langle-1,p,h\rangle|$ for $\mathbb D$ and $N=|\langle p,h\rangle|$ for $\mathbb G$, let $\zeta$ be a primitive $N$th root of unity, and write $p=\zeta^k$, $h=\zeta^{k'}$. Then
\[
\PIdeg B=\frac{N^2}{\gcd\bigl(N,\,2k(k+k')\bigr)} .
\]
\end{proposition}

\begin{proof}
By Lemma~\ref{lem:DGmodel}, $Q(B)$ is the quotient division ring of the quantum torus on $x_1,\vartheta,y_1,y_2$, so $\PIdeg B$ is the PI degree of that torus. Here $q=\zeta^{k+k'}$ and, in $\mathbb D$, $-1=\zeta^{N/2}$. Write $v_iv_j=\zeta^{a_{ij}}v_jv_i$ in the order $x_1,\vartheta,y_1,y_2$. By \eqref{eq:Dtorus} and \eqref{eq:Gtorus}, the entries $(a_{12},a_{13},a_{14},a_{23},a_{24},a_{34})$ are
\begin{align*}
\mathbb D:&\quad \Bigl(\frac N2,\ -\frac N2-k-k',\ -k-k',\ -k,\ k,\ -k\Bigr),\\
\mathbb G:&\quad \bigl(0,\ -k-k',\ -k-k',\ -k,\ k,\ -k\bigr),
\end{align*}
and in both cases
\[
a_{12}a_{34}-a_{13}a_{24}+a_{14}a_{23}=2k(k+k').
\]
The gcd $d_1$ of the entries is $\gcd(N/2,k,k')$ in $\mathbb D$ and $\gcd(k,k')$ in $\mathbb G$, and $\gcd(N,d_1)=1$ because $\zeta^{N/2}$, $\zeta^k$, $\zeta^{k'}$, respectively $\zeta^k$, $\zeta^{k'}$, generate $\langle\zeta\rangle$. By the normal form of alternating integer matrices \cite[Chapter IV]{New}, the matrix is equivalent to two blocks with invariants $d_1\mid d_2$ and $d_1d_2=2k(k+k')$, and the torus has PI degree
\[
\frac{N}{\gcd(N,d_1)}\cdot\frac{N}{\gcd(N,d_2)}=\frac{N^2}{\gcd\bigl(N,2k(k+k')\bigr)},
\]
as in the proof of Theorem~\ref{thm:L}.
\end{proof}

From now on $p\notin\mi$.

\begin{lemma}\label{lem:hw}
Let $M$ be a left $A$-module, $w\in\kk^\times$ and $0\ne v\in M$ with $tv=\lambda(w)v$. Then
\[
\bigl(t-\lambda(pw)\bigr)\bigl(t-\lambda(w/p)\bigr)y_iv=0\qquad(i=1,2).
\]
If $t-\lambda(pw)$ is injective on $M$, then $\Omega v=0$.
\end{lemma}

\begin{proof}
Write $\lambda=\lambda(w)$. By \eqref{eq:Arel},
\[
t\,y_1v=p^{-1}\lambda\,y_1v+\beta_1\,y_2v,\qquad t\,y_2v=p\lambda\,y_2v+\beta_2\,y_1v .
\]
The $2\times2$ matrix of this rule has trace $(p+p^{-1})\lambda$ and determinant $\lambda^2-\beta_1\beta_2$, and by \eqref{eq:aident}
\begin{align*}
\lambda(pw)+\lambda(w/p)&=(p+p^{-1})\lambda(w),\\
\lambda(pw)\lambda(w/p)&=w^2+a(p^2+p^{-2})+a^2w^{-2}=\lambda(w)^2-\beta_1\beta_2 .
\end{align*}
Cayley--Hamilton gives the first claim. If $t-\lambda(pw)$ is injective, then $t\,y_iv=\lambda(w/p)\,y_iv$ for $i=1,2$, and the first displayed rule gives $y_2v=\gamma\,y_1v$ with
\[
\gamma=\frac{\lambda(w/p)-p^{-1}\lambda(w)}{\beta_1}=\frac{a(p^2-1)}{p\beta_1w}=\frac{pb}{w}.
\]
Then $y_1y_2v=\gamma y_1^2v$, $y_2^2v=\gamma\,y_2y_1v=p\gamma^2y_1^2v$, and
\[
t\,y_1^2v=p^{-1}y_1\,t\,y_1v+\beta_1y_2y_1v=\bigl(p^{-1}\lambda(w/p)+p\beta_1\gamma\bigr)y_1^2v .
\]
Since $p^2\beta_1-c\,p^2=-c$, we get
\begin{align*}
\Omega v&=\bigl(\gamma p^{-1}\lambda(w/p)+p\beta_1\gamma^2-p^{-1}b-cp\gamma^2\bigr)y_1^2v\\
&=p^{-1}\bigl(\gamma\lambda(w/p)-c\gamma^2-b\bigr)y_1^2v\\
&=p^{-1}\Bigl(b+\frac{ap^2b}{w^2}-\frac{cp^2b^2}{w^2}-b\Bigr)y_1^2v=0,
\end{align*}
because $a=bc$.
\end{proof}

\begin{corollary}\label{cor:tfinite}
Let $M\ne0$ a left $A$-module with $g(t)M=0$ for some nonzero $g\in\kk[t]$. Then $\Omega v=0$ for some $0\ne v\in M$. Consequently, if $Q$ is a prime ideal of $A$ containing a nonzero element of $\kk[t]$, then $\Omega\in Q$.
\end{corollary}

\begin{proof}
Let $\mathcal E$ be the set of $\mu\in\kk$ with $\ker(t-\mu)\ne0$ on $M$. It is a nonempty set of roots of $g$, and $t-\mu$ is injective for $\mu\notin\mathcal E$. The set $W=\{w\in\kk^\times:\lambda(w)\in\mathcal E\}$ is finite and nonempty, and since $p\notin\mi$ there is $w_0\in W$ with $pw_0\notin W$. Take $v\ne0$ with $tv=\lambda(w_0)v$ and apply Lemma~\ref{lem:hw}. For the second claim apply this to $M=A/Q$: there is $r\notin Q$ with $\Omega r\in Q$, so $(A\Omega A)(ArA)=A\Omega rA\subseteq Q$ and $\Omega\in Q$.
\end{proof}

Next, the weights of the right action of $t$. Let $F=\kk(t)$ and $F(w)$ the extension by a root of $w^2-tw+a$, so that $t=\lambda(w)$; since $a\ne0$ this polynomial is irreducible, and $w$ is transcendental over $\kk$. Put $V=F(w)\otimes_{\kk[t]}A$, using the left $\kk[t]$-structure of $A$. Right multiplication by any nonzero element of $A$ is left $\kk[t]$-linear, so it acts $F(w)$-linearly on $V=\bigoplus_nV_n$, and it is injective because $A$ is a domain and $F(w)$ is flat over $\kk[t]$. Also $\dim_{F(w)}V_n=n+1$. The scalars $\lambda(p^{\nu}w)$, $\nu\in\mathbb Z$, are pairwise distinct: $\lambda(p^{\nu}w)=\lambda(p^{\nu'}w)$ forces $\nu=\nu'$ or $w^2\in\kk$.

\begin{lemma}\label{lem:lower}
Let $\eta\in F(w)^\times$ with $\lambda(p\eta)\ne\lambda(\eta/p)$, and $0\ne m\in V_n$ with $mt=\lambda(\eta)m$. Then $my_1,my_2\in V_{n+1}$ lie in the sum of the $\lambda(p\eta)$- and $\lambda(\eta/p)$-eigenspaces of right multiplication by $t$. The $\lambda(\eta/p)$-components of $my_1$ and $my_2$ are not both zero, and neither are their $\lambda(p\eta)$-components.
\end{lemma}

\begin{proof}
By Lemma~\ref{lem:Omega}, with $\mu=\lambda(\eta)$,
\[
(my_1)t=p\mu\,my_1-p\beta_1\,my_2,\qquad (my_2)t=p^{-1}\mu\,my_2-p^{-1}\beta_2\,my_1 ,
\]
a matrix with trace $(p+p^{-1})\mu$ and determinant $\mu^2-\beta_1\beta_2$, hence with the distinct eigenvalues $\lambda(p\eta)$ and $\lambda(\eta/p)$, as in Lemma~\ref{lem:hw}. Suppose the $\lambda(\eta/p)$-components of $my_1$ and $my_2$ vanish. Then the first equation gives $my_2=\gamma'my_1$ with
\[
\gamma'=\frac{p\lambda(\eta)-\lambda(p\eta)}{p\beta_1}=\frac{a(p^2-1)}{p^2\beta_1\eta}=\frac b\eta ,
\]
and as before, using $y_1y_2=p^{-1}y_2y_1$,
\[
m\Omega=p^{-1}\bigl(\lambda(\eta)\gamma'-b-c\gamma'^2\bigr)my_1^2=p^{-1}\,\frac{b(a-bc)}{\eta^2}\,my_1^2=0,
\]
which is impossible because right multiplication by $\Omega$ is injective on $V$. For the $\lambda(p\eta)$-components apply the same argument to $\eta'=a/\eta$, noting that $\lambda(\eta')=\lambda(\eta)$ and $\lambda(\eta'/p)=\lambda(p\eta)$.
\end{proof}

\begin{proposition}\label{prop:weights}
\begin{enumerate}
\item[(i)] Right multiplication by $t$ on $V_n$ is diagonalizable, with the $n+1$ distinct eigenvalues $\lambda(p^{\nu}w)$, $\nu\in\{-n,-n+2,\dots,n\}$. For $n=2k$ the $\lambda(w)$-eigenspace is $F(w)\,\Omega^k$.
\item[(ii)] Every nonzero graded ideal of $A$ contains a power of $\Omega$.
\end{enumerate}
\end{proposition}

\begin{proof}
(i) Induction on $n$, with $V_0=F(w)$ and $1\cdot t=\lambda(w)$. Since $A_n=A_{n-1}y_1+A_{n-1}y_2$, Lemma~\ref{lem:lower}, applied to $\eta=p^{\nu}w$, shows that $V_n$ is spanned by eigenvectors with eigenvalues $\lambda(p^{\nu}w)$, $|\nu|\le n$, $\nu\equiv n\pmod2$, and that each of these $n+1$ values occurs. As $\dim V_n=n+1$, each eigenspace is a line. Finally $\Omega^k\in A_{2k}$ commutes with $t$, so $\Omega^kt=\lambda(w)\Omega^k$.

(ii) Let $J\ne0$ be a graded ideal and $J_n\ne0$. Then $U_n=F(w)\otimes J_n\subseteq V_n$ is nonzero and stable under right multiplication by $t$, so it contains an eigenvector $m$, say with eigenvalue $\lambda(p^{\nu}w)$. Suppose $\nu>0$. By Lemma~\ref{lem:lower}, one of $my_1,my_2\in U_{n+1}$ has a nonzero $\lambda(p^{\nu-1}w)$-component, and this component lies in $U_{n+1}$ because the eigenprojections are polynomials in right multiplication by $t$. After $|\nu|$ steps (lowering if $\nu>0$, raising if $\nu<0$) we find $\Omega^k\in U_{2k}$ for some $k$, by (i). So the class of $\Omega^k$ in $A_{2k}/J_{2k}$ dies in $F(w)\otimes(A_{2k}/J_{2k})$, hence is $\kk[t]$-torsion: $g(t)\Omega^k\in J$ for some nonzero $g$.

Now let $Q_1,\dots,Q_\ell$ be prime ideals containing $J$ with $Q_1\cdots Q_\ell\subseteq J$. Each $Q_i$ contains $g(t)\Omega^k$, hence $(Ag(t)A)(A\Omega^kA)$, so it contains $g(t)$ or $\Omega^k$; in both cases $\Omega\in Q_i$, by Corollary~\ref{cor:tfinite} in the first. Therefore $\Omega^\ell\in J$.
\end{proof}

\begin{theorem}\label{thm:DGprim}
Let $B$ be $\mathbb D(p,h)$ with $p^2\ne1$, or $\mathbb G(p,f,h)$, and suppose that $B$ is not PI. Then $B$ is primitive if and only if neither $p$ nor $hp$ is a root of unity, that is, if and only if $Z(Q(B))=\kk$.
\end{theorem}

\begin{proof}
One direction is Theorem~\ref{thm:DG}(iii). Let $p,q\notin\mi$, $q=hp$. The elements $x_1$ and $\Xi=x_1\Omega$ are normal in $B$ (Remark~\ref{rem:DGnormal}). We show that every nonzero prime $P$ of $B$ contains $x_1\Xi$; then $(0)$ is locally closed, and $B$ is primitive by \cite[Lemma II.7.15]{BG} and \cite[Theorem 0.4]{ASZ}, as in the proof of Theorem~\ref{thm:primmono} below.

Let $x_1\notin P$. Then $P'=P\,B[x_1^{-1}]$ is a nonzero ideal of $A[x_1^{\pm1};\alpha]$ with $P'\cap B=P$. Conjugation by $x_1$ multiplies $\Omega$ by $q^{-2}$, so $\Omega x_1^i=q^{2i}x_1^i\Omega$. Take $0\ne g=\sum_{i\in I}\theta_ix_1^i\in P'$, $\theta_i\in A$, with $|I|$ minimal, and $j\in I$. Then
\[
\Omega g-q^{2j}g\Omega=\sum_{i\in I}(q^{2i}-q^{2j})\,\theta_ix_1^i\,\Omega\in P'
\]
has fewer terms, so it vanishes, and since $q\notin\mi$ and $A$ is a domain, $g=\theta_jx_1^j$. Hence $J=P'\cap A\ne0$. It is stable under $\alpha$, which acts on the basis by
\[
\alpha(t^iy_1^jy_2^l)=\varepsilon\,q^{-(j+l)}\,t^iy_1^jy_2^l,\qquad \varepsilon\in\{\pm1\},
\]
($\varepsilon=(-1)^{i+j}$ for $\mathbb D$, $\varepsilon=1$ for $\mathbb G$). Since $q\notin\mi$, elements of different degree have different eigenvalues, and a Vandermonde argument shows that $J$ is graded. By Proposition~\ref{prop:weights}(ii), $\Omega^\ell\in J$ for some $\ell$. Then $\Xi^\ell$, a nonzero multiple of $x_1^\ell\Omega^\ell$, lies in $P'\cap B=P$, and since $\Xi$ is normal and $P$ is prime, $\Xi\in P$.
\end{proof}

\begin{remark}\label{rem:DGnormal}
In $\mathbb D$ the element
\[
\Xi=x_2y_1y_2-\kappa_0x_1\bigl(p^{-1}y_1^2+y_2^2\bigr)=x_1\Omega
\]
is normal: $\Xi x_1=q^2x_1\Xi$, $\Xi y_1=-q^{-1}y_1\Xi$, $\Xi y_2=q^{-1}y_2\Xi$ and $\Xi x_2=-q^2x_2\Xi$, because $\Omega$ is central in $A$ and $x_1$ normalizes $A$. In $\mathbb G$, $\Xi=x_1\Omega$ satisfies $\Xi x_i=q^2x_i\Xi$ and $\Xi y_i=q^{-1}y_i\Xi$ \cite{GV-G}. For $\mathbb G$ and $h=1$ this is the element called $C$ in \cite{GV-G}. Unlike in the monomial families, the normal elements we use only see the $x$-side, and Proposition~\ref{prop:weights}(ii) plays the role of a second normal element. 
\end{remark}

\subsection{The monomial families}

For the families of Section~\ref{sec:monodeg} the crossed product settles primitivity for all parameters: Proposition~\ref{prop:crossed}(i)--(ii), the relations of Section~\ref{sec:monodeg} and the splitting of $T'$ in the proof of Theorem~\ref{thm:monodeg} do not need the PI hypothesis. The two factors have parameters
\begin{equation}\label{eq:twoparams}
\begin{array}{c|c|c}
 & \text{first factor} & \text{second factor}\\\hline
\mathbb S,\ \mathbb T,\ \mathbb U & -1 & 16h^4\\
\mathbb E,\ \mathbb J & -1 & 4h^4\\
\mathbb M,\ \mathbb N,\ \mathbb O,\ \mathbb P,\ \mathbb W,\ \mathbb Z & u^{\pm2} & v^2
\end{array}
\end{equation}
with $u,v$ as in Theorem~\ref{thm:monodeg}.

\begin{lemma}\label{lem:loc}
Let $N_x\in A_x$ and $N_y\in A_y$ be the elements
\[
\begin{array}{c|ccccc}
 & \mathbb M,\ \mathbb O & \mathbb P & \mathbb Z & \mathbb W & \mathbb N\\\hline
N_x & e_1^2 & e_1^2 & e_1^2 & e_1e_2 & x_1x_2\\
N_y & z_1^2 & y_1y_2 & z_1z_2 & z_1^2 & z_1^2
\end{array}
\]
They are normal in $B$, and the localization of $B$ at the multiplicative set they generate is the ring
\[
\Lambda=T'\oplus T'e\oplus T'z\oplus T'ez\subseteq Q(B).
\]
\end{lemma}

\begin{proof}
In each case $N_x$ is normal in $A_x$: $e_1^2$ is central for a balanced basis, and $e_1e_2$, $x_1x_2$ are normal in a commutative or coordinate plane. By \eqref{eq:mono}, $z_aN_x$ is a scalar multiple of $N_x'z_a$ with $N_x'\in\{e_1^2,e_2^2\}$ or $\{e_1e_2,e_2e_1\}$, and $N_x'$ is a scalar multiple of $N_x$ in the cases listed. So $N_x$ is normal in $B$, and the same argument applies to $N_y$. In $\Lambda$ both are units, so $B[N_x^{-1},N_y^{-1}]\subseteq\Lambda$. Conversely, inverting $N_x$ inverts $e_1$ and $e_2$ (for instance $e_1^{-1}=e_1N_x^{-1}$ when $N_x=e_1^2$ and the basis is balanced, and $e_1^{-1}=e_2N_x^{-1}$ when $N_x=e_1e_2$), so it contains $r^{\pm1}$ and $X^{\pm1}$; the same holds on the $y$-side, and $\Lambda\subseteq B[N_x^{-1},N_y^{-1}]$.
\end{proof}

\begin{theorem}\label{thm:primmono}
\begin{enumerate}
\item[(i)] The members of $\mathbb E$, $\mathbb J$, $\mathbb S$, $\mathbb T$, $\mathbb U$, $\mathbb K$ and $\mathbb L$ are never primitive.
\item[(ii)] Let $B$ be a member of $\mathbb M$, $\mathbb N$, $\mathbb O$, $\mathbb P$, $\mathbb W$ or $\mathbb Z$, with $f\ne0$ in $\mathbb M,\mathbb O,\mathbb P,\mathbb W$, and let $u,v$ be as in Theorem~\ref{thm:monodeg}. Then $B$ is primitive if and only if neither $u$ nor $v$ is a root of unity.
\end{enumerate}
\end{theorem}

\begin{proof}
Suppose one of the two factors in \eqref{eq:twoparams} has a parameter in $\mi$; this always happens in the first two rows. Then $T'$ has a central monomial $c$ other than a scalar, for instance $r^2$ when the first factor has parameter $-1$. Put
\[
\prod_{g\in G}\bigl(T-g(c)\bigr)\in K[T],
\]
where $G$ is the group of automorphisms of $K$ induced by $e$ and $z$. It is a quotient of $(\mathbb Z/2)^2$, since $e^2,z^2\in T'$ and $ze\in T'^{\times}ez$, and $K^G$ commutes with $Q(T')$, $e$ and $z$, so it is central in $Q(B)$. The coefficients lie in $K^G$, and they are not all in $\kk$, because $c$ is transcendental over $\kk$ and is a root of this polynomial. So the center of $Q(B)$ is not algebraic over $\kk$, and $B$ is not primitive by Proposition~\ref{prop:rational}. For $\mathbb S$, for instance, $r^2+r^{-2}$ is already central in $Q(B)$; in terms of $x_1,x_2$,
\[
r^2+r^{-2}=\frac{2\bigl(x_1^4-6x_1^2x_2^2+x_2^4\bigr)}{(x_1^2+x_2^2)^2},
\]
for every $h$. For $\mathbb K$ and $\mathbb L$ the elements
\[
\frac{f^2y_1^4+y_2^4}{y_1^2y_2^2}\in Q(\mathbb K),\qquad \frac{y_1^4+y_2^4}{y_1^2y_2^2}\in Q(\mathbb L)
\]
are central for all $q,f,h$. Indeed, the mixed relations of $\mathbb K$ are $y_1x_1=hx_1y_1$, $y_2x_1=hx_1y_2$, $y_1x_2=hx_2y_2$, $y_2x_2=hfx_2y_1$, and those of $\mathbb L$ are in Section~\ref{sec:normal}. Applying them twice,
\begin{align*}
\mathbb K:&\quad y_j^2x_1=h^2x_1y_j^2, & y_1^2x_2&=h^2x_2y_2^2, & y_2^2x_2&=h^2f^2x_2y_1^2,\\
\mathbb L:&\quad y_1^2x_1=h^2f^2x_1y_2^2, & y_2^2x_1&=h^2f^2x_1y_1^2, & y_1^2x_2&=h^2x_2y_2^2,\quad y_2^2x_2=h^2x_2y_1^2 .
\end{align*}
So in $\mathbb K$ both $f^2y_1^4+y_2^4$ and $y_1^2y_2^2$ are multiplied by $h^4$ when moved past $x_1$ and by $h^4f^2$ past $x_2$, and in $\mathbb L$ both $y_1^4+y_2^4$ and $y_1^2y_2^2$ are multiplied by $h^4f^4$ past $x_1$ and by $h^4$ past $x_2$. Moreover $y_1^2,y_2^2$ are central in $\kk_{-1}[y_1,y_2]$. For $\mathbb J$ use \eqref{eq:EJ}: $Q(\mathbb J)$ is anti-isomorphic to $Q(\mathbb E)$, so it has the same center. This proves (i) and one direction of (ii).

Now let $u,v\notin\mi$ in (ii). The quantum torus $T'$ is a tensor product of two rank-two tori with parameters that are not roots of unity, so its center is $\kk$, and $T'$ is simple \cite[Proposition 1.3]{MP}. Its units are the nonzero scalar multiples of monomials, so an inner automorphism of $T'$ multiplies each monomial by a scalar. Here $e$, $z$ and $ez$ act on $T'$ by conjugation, and $\Lambda$ is a crossed product of $T'$ by the group $(\mathbb Z/2)^2$ of these actions modulo inner automorphisms (recall $e^2=X$, $z^2=Y$ and $ze\in T'^{\times}ez$). We claim that the automorphisms of $T'$ induced by $e$, $z$ and $ez$ are outer. The elements that invert $r$ or $\rho$ are not of that form. The remaining nontrivial element is $ez$ in $\mathbb M$ and $\mathbb O$, $z$ in $\mathbb P$ and $\mathbb Z$, and $e$ in $\mathbb N$ and $\mathbb W$, and it multiplies $r$ or $\rho$ by one of $u^{\pm1}$, $-u$. Conjugation by a monomial of $T'$ multiplies $r$ and $\rho$ by elements of $\langle u^2\rangle$, by the relations of Section~\ref{sec:monodeg}. If $\pm u^{\pm1}\in\langle u^2\rangle$ then $u^{2j\mp1}=\pm1$ for some $j$, and $u\in\mi$. This proves the claim.

Now $\Lambda$ is simple. Let $I\ne0$ be an ideal and $x=\sum_ga_gw_g\in I$ a nonzero element with the fewest nonzero terms, where $w_g\in\{1,e,z,ez\}$ and $a_g\in T'$. Multiplying by some $w_g^{-1}$ we may assume $a_1\ne0$. The coefficients $a_1$ of the elements of $I$ with support contained in that of $x$ form a nonzero ideal of $T'$, so we may take $a_1=1$. For $t\in T'$ the element $tx-xt=\sum_g\bigl(ta_g-a_g(w_g\triangleright t)\bigr)w_g$ has fewer terms, so $ta_g=a_g(w_g\triangleright t)$ for all $t$. If $a_g\ne0$ for some $g\ne1$, then $a_gT'=T'a_g$ is a nonzero ideal, so $a_g$ is a unit and $w_g\triangleright$ is conjugation by $a_g^{-1}$, against the claim. Hence $x=1$.

Let $P\ne0$ be a prime ideal of $B$. If $P$ missed the multiplicative set generated by $N_x$ and $N_y$, then $P\Lambda$ would be a proper nonzero ideal of $\Lambda$ by Lemma~\ref{lem:loc}. So $P$ contains a product of copies of $N_x$ and $N_y$, and since these are normal and $P$ is prime, $P$ contains $N_x$ or $N_y$, hence $N_xN_y$. Therefore the intersection of the nonzero primes of $B$ contains $N_xN_y\ne0$, and $(0)$ is locally closed. For a noetherian algebra satisfying the Nullstellensatz, locally closed primes are primitive \cite[Lemma II.7.15]{BG}, and $B$ satisfies the Nullstellensatz \cite[Theorem 0.4]{ASZ}.
\end{proof}

\begin{remark}
In (ii), $u,v\in\mi$ gives a PI algebra that is not primitive, $u,v\notin\mi$ a primitive algebra that is not PI, and if exactly one of $u,v$ is a root of unity $B$ is neither. In the families of (i), and in $\mathbb B$, the non-PI members are neither PI nor primitive.
\end{remark}

\subsection{Further nonprimitive members}\label{subsec:HXtorsion}

The following proposition treats the families $\mathbb H$ and
$\mathbb X$ for all twists and several remaining specializations
at roots of unity.

\begin{proposition}\label{prop:HXprim}
\begin{enumerate}
\item[(i)]
Let $B$ be $\mathbb H(f,h)$ or $\mathbb X(h)$, and put
\[
S_+=y_1^2+y_2^2,\qquad S_-=y_1^2-y_2^2,
\qquad \varsigma=S_-S_+^{-1}\in Q(B).
\]
Then $\varsigma$ anticommutes with $x_1,x_2$ and commutes with $y_1,y_2$.
Consequently, $\varsigma^2$ is central and transcendental over $\kk$,
and no member of $\mathbb H$ or $\mathbb X$ is primitive.

\item[(ii)]
Let $h\in\mi$, $n=\ord(h)$ and $m=\operatorname{lcm}(2,n)$.
Then $x_1^n$ is central in $\mathbb A(h)$ and $\mathbb Y(f,h)$,
$x_1^m$ is central in $\mathbb D(\pm1,h)$, $\mathbb O(0,h)$,
$\mathbb P(0,h)$ and $\mathbb W(0,h)$, and $x_2^m$ is central
in $\mathbb M(0,h)$. In particular, these members are not primitive.
\end{enumerate}
\end{proposition}

\begin{proof}
(i) In both families $A_y=\kk_{-1}[y_1,y_2]$, so $S_+,S_-$
are central in $A_y$. The mixed relations are
\[
y_1x_1=hx_1y_2,\qquad y_2x_1=hx_1y_1,
\qquad y_1x_2=hL y_2,\qquad y_2x_2=hL y_1,
\qquad L=fx_1+x_2,
\]
where $f=1$ in $\mathbb X$. Let $\psi$ be the automorphism
of $A_x$ defined by
\[
\psi(x_1)=x_1,\qquad \psi(x_2)=2fx_1+x_2.
\]
This triangular substitution preserves both the Jordan-plane
relation in $\mathbb H$ and the commutative-plane relation in
$\mathbb X$. Applying the mixed relations twice gives
\[
y_1^2x=h^2\psi(x)y_2^2,\qquad
y_2^2x=h^2\psi(x)y_1^2
\qquad(x\in\kk x_1+\kk x_2).
\]
Hence $S_\pm x=\pm h^2\psi(x)S_\pm$ and
\[
S_+^{-1}x=h^{-2}\psi^{-1}(x)S_+^{-1},
\qquad
\varsigma x=h^{-2}S_-\psi^{-1}(x)S_+^{-1}=-x\varsigma.
\]
Since $\varsigma\in\kk(y_1^2,y_2^2)$, it commutes with $y_1,y_2$,
so $\varsigma^2$ is central in $Q(B)$. By the PBW basis,
$y_1^2,y_2^2$ are algebraically independent, and $\varsigma^2$ is a
nonconstant element of the purely transcendental field
$\kk(y_1^2,y_2^2)$. Thus $\varsigma^2$ is transcendental over $\kk$,
and Proposition~\ref{prop:rational} applies.

(ii) In $\mathbb A$ and $\mathbb Y$, $x_1$ commutes with $x_2$
and $y_jx_1=hx_1y_j$, so $x_1^n$ is central.
For $\mathbb D(\pm1,h)$, the relations
\[
x_2x_1=-x_1x_2,\qquad
y_1x_1=-hp\,x_1y_1,\qquad
y_2x_1=hp\,x_1y_2
\]
show that $x_1^m$ is central, since $m$ is even, $p=\pm1$
and $h^m=1$. The case $\mathbb O(0,h)$ follows from
$\mathbb O(0,h)=\mathbb D(-1,h)$
(Proposition~\ref{prop:noPI}).
In $\mathbb M(0,h)$, $\mathbb P(0,h)$ and $\mathbb W(0,h)$,
let $\nu=x_2,x_1,x_1$, respectively.
By the relations used in Proposition~\ref{prop:noPI},
\[
y_1\nu=\varepsilon h\nu y_2,\qquad
y_2\nu=\varepsilon h\nu y_1,
\]
with the same sign $\varepsilon\in\{1,-1\}$.
Thus $y_j\nu^2=h^2\nu^2y_j$, and $\nu^2$ is central in $A_x$.
Since $m$ is even and $h^m=1$, $\nu^m$ is central in $B$.
All the displayed central powers are nonzero homogeneous elements
of positive degree, hence transcendental over $\kk$.
Proposition~\ref{prop:rational} applies again.
\end{proof}

Together, these results decide primitivity for every non-PI member, except in the following cases.

\begin{question}\label{q:remainingprim}
Let $h\notin\mi$. Which members of $\mathbb A(h)$ and
$\mathbb Y(f,h)$ are primitive? Which of the degenerate members
$\mathbb D(\pm1,h)$, $\mathbb M(0,h)$, $\mathbb O(0,h)$,
$\mathbb P(0,h)$ and $\mathbb W(0,h)$ are primitive?
Recall that $\mathbb O(0,h)=\mathbb D(-1,h)$.
None of these algebras is PI.
\end{question}

\section*{Declarations}
\textbf{Funding.} This research received no specific grant from any funding
agency in the public, commercial or not-for-profit sectors.

\textbf{Competing interests.} The authors have no competing interests to
declare that are relevant to the content of this article.

\textbf{Author contributions.} Both authors contributed equally to this work
and read and approved the final manuscript.

\textbf{Declaration of AI use.} During the preparation of this work,
the authors used ChatGPT (OpenAI) to assist in developing a Python
notebook for Google Colab to simplify products in the algebras considered
and to check selected computations carried out by hand. The notebook
and verification code are publicly available in the accompanying
GitHub repository \cite{GV-code}. Claude (Anthropic) and ChatGPT were
also used for literature searches, the development and preliminary
checking of mathematical arguments and proofs, and editing the English
text and LaTeX code. The authors have reviewed and verified all
mathematical arguments, proofs, computations and references and take
full responsibility for the final content.

\end{document}